\documentclass[leqno,twoside]{amsart}
\usepackage{latexsym,amsmath,amssymb, url}
\usepackage{graphicx}
\usepackage{pgf,tikz}
\usetikzlibrary{arrows}

\title [rigidity of Sobolev  isometric immersions]{Rigidity and regularity of co-dimension one
Sobolev isometric immersions}

\author{Zhuomin Liu and Mohammad Reza Pakzad}
\address{Zhuomin Liu, Department of Mathematical Analysis, Charles University, Czech Republic} 
\address{Mohammad Reza Pakzad, Department of Mathematics, 
University of Pittsburgh, 301, Thackeray Hall, Pittsburgh, PA 15260, USA}

\thanks{}

\def\ds{\displaystyle}

\def\R{{\mathbb R}}

\newtheorem{theorem}{Theorem}
\newtheorem{lemma}[theorem]{Lemma}

\newtheorem{corollary}[theorem]{Corollary}
\newtheorem{proposition}[theorem]{Proposition}

\theoremstyle{definition}
\newtheorem{remark}[theorem]{Remark}
\newtheorem{definition}[theorem]{Definition}

\newcommand{\bbbr}{\mathbb R}

\numberwithin{theorem}{section}
\numberwithin{equation}{section}

\begin{document}

\definecolor{qqttzz}{rgb}{0,0.2,0.6}
\definecolor{qqqqqq}{rgb}{0,0,0}
\definecolor{qqqqqq}{rgb}{0,0,0}
\definecolor{qqffff}{rgb}{0,1,1}
\definecolor{qqqqff}{rgb}{0,0,1}
\definecolor{zzttqq}{rgb}{0.6,0.2,0}
\definecolor{ccffcc}{rgb}{0.8,1,0.8}
\definecolor{yqyqyq}{rgb}{0.502,0.502,0.502}
\definecolor{dqfqcq}{rgb}{0.816,0.941,0.753}
\definecolor{cqcqcq}{rgb}{0.753,0.753,0.753}
\definecolor{eqeqeq}{rgb}{0.878,0.878,0.878}
\definecolor{ffffqq}{rgb}{1,1,0}
\definecolor{qqzzff}{rgb}{0,0.6,1}
\definecolor{zzffqq}{rgb}{0.6,1,0}
\definecolor{ttttff}{rgb}{0.2,0.2,1}
\definecolor{ttzzqq}{rgb}{0.2,0.6,0}
\definecolor{qqzzqq}{rgb}{0,0.6,0}
\definecolor{zzffzz}{rgb}{0.6,1,0.6}
\definecolor{ffffzz}{rgb}{1,1,0.6}

\email{liuzhuomin@hotmail.com, pakzad@pitt.edu}




\begin{abstract}
We prove the developability and $C_{\rm loc}^{1,1/2}$ regularity of $W^{2,2}$
isometric immersions of $n$-dimensional domains into $\bbbr^{n+1}$. As a conclusion  
we show that any such Sobolev isometry can be approximated by smooth isometries in the $W^{2,2}$ 
strong norm, provided the domain is $C^1$ and convex. Both results fail to be true if the Sobolev 
regularity is weaker than $W^{2,2}$.  
 \end{abstract}

\maketitle
\setcounter{tocdepth}{1}
\tableofcontents
\section{Introduction}

It has been known since at least the 19th century that any smooth surface with zero Gaussian curvature is
locally ruled, i.e.  passing through any point of the surface  is a straight segment lying on the surface. Such surfaces were 
called developable surfaces. This terminology was used as an indication that any such surface is in isometric equivalence 
with the plane, i.e. any piece of it can be {\it developed} on the flat plane without any stretching or compressing. Meanwhile, it was 
already suspected that there exist 
somewhat regular surfaces applicable to the plane, but yet not developable (See \cite{cajori} for a review of this question). 
Nevertheless, it was not until the work of John Nash at the zenith of the last century that the existence of such unintuitive phenomena was rigorously 
established.

In his pioneering work, Nash settled several questions. He established that any Riemannian manifold can be isometrically embedded in an 
Euclidean space \cite{Nash1}.
Moreover, if the dimension of the space is large enough, this embedding can be done in a manner so that the diameter of the image is 
as small as one wishes. 
As for the lower dimensional embeddings, Nash \cite{Nash2} and Kuiper \cite{kuiper}, established 
the existence of a $C^1$ isometric embedding of 
any Riemannian manifold into another manifold of one higher dimension. Their method, which is now famously re-cast in the framework of convex 
integration 
\cite{gromov}, involved iterated perturbations of a given short mapping of the manifold towards realizing an isometry.

A surprising corollary of these results is the existence of a $C^1$ flat torus in $\R^3$ \cite{borrelli}. 
Another one is that there are $C^1$ isometric 
embeddings of 
the two dimensional unit sphere into three dimensional space with arbitrarily small diameter. By contrast, it was established by Hartman 
and Nirenberg that any flat $C^2$ surface in $\R^3$ must be developable \cite{HartmanNirenberg}, while Hilbert had already shown that 
any $C^2$ isometric 
immersion of the sphere  must be a rigid motion. This latter result is a special case of a similar statement for any closed convex surface 
in $\R^3$, 
see \cite[Chapter 12]{spivak}. On the other hand, the former result was generalized by Pogorelov's  for $C^1$ isometries with
total zero curvature in \cite[Chapter II]{Po 56} and \cite[Chapter IX]{Po 73}.  
 
A natural question arises in this context for the analyst: What about isometric immersions of intermediate regularity, say of H\"older or 
Sobolev type? 
Regarding H\"older regularity, rigidity of $C^{1,\alpha}$ isometries of 2 dimensional flat domains has been established for $\alpha\ge 2/3$ 
\cite{bori1,bori2}, 
while their flexibility in the sense of Nash and Kuiper
is known for $\alpha<1/7$ \cite{bori2, CDS}. The critical value for $\alpha$ is conjectured to be $1/2$ in this case.  
As for the regularity of Sobolev isometries, following the results of Kirchheim in 
\cite{kirchheimthesis} on $W^{2,\infty}$ solutions to degenerate Monge-Amp\`ere equations (see Proposition \ref{0pakzad2}),  
the rigidity of $W^{2,2}$ isometries of a flat domain was  established  in \cite{Pak}. More precisely, it was established that 
such mappings are developable in the classical sense, i.e.

\begin{theorem}[Pakzad \cite{Pak}] \label{0pakzad1}
Let $v\in W^{2,2}(\Sigma, \bbbr^3)$ be an isometric immersion, where $\Sigma$ is a bounded Lipschitz domain in
$\bbbr^2$. Then $v\in C_{\rm loc}^{1,1/2}(\Sigma, \bbbr^2)$. Furthermore, for
every point of $x$, either there exists a neighborhood of $x$, or a unique
segment passing through $x$ and joining $\partial \Sigma$ at both ends, on which
$\nabla v$ is constant.
\end{theorem}

\begin{remark}
It can be shown that this statement is actually valid for all bounded open sets $\Sigma \subset \R^2$, i.e.
without any assumption on the regularity of the boundary. All one must prove is that the constancy segments, 
whose existence are locally established, can be extended all the way to the boundary one step at a time. Assuming the existence of 
any supposedly maximal constancy segment which does not reach the boundary, a contradiction could be achieved by 
creating a Lipschitz domain $\Sigma' \subset \Sigma$ including the closure of that segment and applying Theorem {\ref{0pakzad1}} to $\Sigma'$. 
In the same manner, the regularity assumption on $\partial \Omega$ in Theorem \ref{0deve-regu} can be removed.
\end{remark}  

To put this result in context, it is noteworthy 
that a $W^{2,2}$ function on a two dimensional domain fails barely to be $C^1$, but there is information available about second weak derivatives, 
and e.g.   the Gaussian curvature of the image of a $W^{2,2}$ isometric immersion of a flat domain  is identically zero as an $L^1$ function.  
This indicates that these isometries are far from the highly oscillatory solutions of Nash and Kuiper and hence possibly should behave in a rigid manner. 
Note that only the $C^1$ regularity result was stated in \cite{Pak} and was a major ingredient of the proof, but the higher H\"older regularity announced here is 
an immediate consequence of the developability. In \cite{MuPa2} it was established that the $C^1$ regularity can be extended up to the boundary 
if the domain is of class $C^{1,\alpha}$. This does not hold true anymore for merely $C^1$ regular domains. Finally, the following proposition is a key 
step in establishing the above rigidity result and will be instrumental in proving Theorem \ref{0deve-regu}.

\begin{proposition}[Kirchheim \cite{kirchheimthesis}, Pakzad \cite{Pak}] \label{0pakzad2} 
Let $\Sigma$ be as above and let $f\in W^{1,2}(\Sigma, \bbbr^3)$ be a map with almost everywhere symmetric and singular (i.e. of zero determinant) gradient. 
Then $f \in C^0(\Sigma)$ and for every point $x\in \Sigma$, there exists either a neighborhood $U$ of $x$, 
or a segment passing through it and joining $\partial \Sigma$ at its both ends, 
on which $f$ is constant.
\end{proposition}

It was proved furthermore in \cite{Pak} 
that any $W^{2,2}$ isometry on a convex 2d domain can  be approximated in strong norm by smooth isometries. 
This is a nontrivial result, since the usual regularization techniques fail due to the non-linearity of the isometry 
constraint. The idea was to make use of the developability structure of these mappings and reduce the approximation problem to 
the one about mollifying the expressions $R^TR'$ for the Darboux moving frames $R(t)$ along the curves orthogonal to the rulings. The convexity 
assumption is a technical one, 
and as shown by Hornung \cite{Ho0, Ho2}, can be replaced by e.g. piece-wise $C^1$ regularity of the boundary, see also \cite{Ho1}.

It is natural to ask whether these results can be generalized to higher dimensions. In \cite{VWKG}, the authors showed the generalized developability of 
smooth isometric immersions of Euclidean domains into Euclidean spaces.  We would like to pose the same problem for the same class of  
isometric immersions but only considered under sufficient Sobolev regularity assumptions.  A first main result in this direction, presented in this paper, 
is the developability of $W^{2,2}$ 
co-dimension one isometries 
of flat domains in ${\mathbb R}^n$: 

\begin{theorem}\label{0deve-regu}
Let $u\in W^{2,2}(\Omega, \bbbr^{n+1})$ be an isometric immersion, where $\Omega$ is a bounded Lipschitz domain in
$\bbbr^n$. Then $u\in C_{\rm loc}^{1,1/2}(\Omega, \bbbr^{n+1})$. Moreover, for
every $x\in \Omega$, either $\nabla u$ is constant in a neighborhood of $x$, or
there exists a unique $(n-1)$-dimensional hyperplane $P\ni x$ of $\bbbr^n$ such
that $\nabla u$ is constant on the connected component of $x$ in $P\cap \Omega$.
\end{theorem}

The interesting feature of this new result is that the Sobolev regularity $W^{2,2}$ is much below the
required $W^{2,n+\varepsilon}$ for obtaining $C^1$ regularity. An extra difficulty which comes in the way of the proof in dimensions 
higher than 2 is that the  argument used in Lemma 2.1 of \cite{Pak} to show the continuity of the derivatives of the given Sobolev isometry is 
no more generalizable to our case. Indeed, in \cite{Pak}, a very important first step of the proof of developability  is to show the $C^1$ regularity. 
Here, on the other hand, we first show  the developability of the mapping without having the $C^1$ regularity at hand. 
Our  proof is based an induction on the dimension of slices of the domain and careful and detailed geometric arguments. Having established developability, the 
$C^1$ regularity (and better) follows in a straightforward manner. 
 
The problem of regularity and developability of Sobolev isometric immersions of  co-dimension $k>1$ is more involved and could not be tackled 
through the methods discussed in this paper. 
In a forthcoming paper  by Jerrard and the second author \cite{jp},  another approach, more analytical in nature, is adapted to study this problem. 
It is based on the fact that the Hessian rank inequality
\begin{equation}\label{analytic}
\displaystyle \mbox{rank} (\nabla^2 v) \le k \quad \mbox{a.e. in} \,\, \Omega
\end{equation} is satisfied by the components $v=u^j$ of such isometry. Note that this equation becomes the degenerate Monge-Amp\`ere equation when $k=n-1$. 
Similar as in \cite {Pak}, regularity and developability of the Sobolev solutions to (\ref{analytic}) directly implies the same results for the corresponding 
isometries. However, one loses some natural advantages when working with (\ref{analytic}) rather than with the isometries themselves as
done in the present paper: the solution $v$  is no more Lipschitz and being just a scalar function, one loses the extra information  derived from the length 
preserving properties of isometries.  Methods of geometric measure theory applied to the class of Monge-Amp\`ere functions developed by Jerrard in 
\cite{Je1, Je2} are used to overcome these obstacles.
 
The second main result of this paper concerns approximation of $W^{2,2}$ isometries by smooth ones: 

\begin{theorem}\label{0density}
Assume $\Omega \subset \R^n$ is a $C^1$ bounded convex domain and that  $u\in W^{2,2}(\Omega, \bbbr^{n+1})$ is an isometric immersion. 
Then there is a sequence of isometric immersions $u_m\in 
C^\infty(\overline \Omega, \bbbr^{n+1})$ converging to $u$ in $W^{2,2}$ norm.
\end{theorem}
 
The main idea of the proof, similar as in the 2d case, is to mollify the curves which pass orthogonally through the constancy hyperplanes of Theorem \ref{0deve-regu}
both in the domain and on the image. This latter problem, framed within the general isometry mollification problem, is still nonlinear.  
However, identifying these curves with suitable orthonormal moving 
Darboux frames $R(t) \in SO(n)$ and $ \tilde R(t) = [(\nabla u) R(t), {\bf n}(t)] \in SO(n+1)$, where ${\bf n}$ is the unit normal to the image 
of the isometry in $\R^{n+1}$,  we could linearize the problem by considering the curvature matrices  $R^T R'(t) \in so(n)$ and $\tilde R^T \tilde R'(t) \in so(n+1)$ and 
recover an approximating sequence of moving frames through their regularization. Many technical details must nevertheless be taken care of in this process; in particular 
one must make sure that the mollified curves can be used to define new smooth isometries. Also, the mapping as a whole cannot be described by one single
couple of such curves and the domain must be partitioned into suitable subdomains.  

\begin{remark} 
Neither the $C^1$ regularity nor the convexity of the boundary 
seems to be absolutely necessary for the density result to hold true 
(see e.g. \cite{Ho2} for finer results in 2d), but omitting these assumptions  goes beyond the scope of our paper. However, 
both of the results in Theorems \ref{0deve-regu} and \ref{0density} are sharp in the sense that  
they fail to be true if the isometric immersion is only of class $W^{2,p}$ for $p<2$. An immediate counterexample 
is the following isometric immersion $u:B^2\times (0,1)^{n-2} \to \R^{n+1}$, whose image can be visualized as a family of cones over a hyperplane 
of dimension $n-2$:
$$
\ds u(r\cos \theta, r \sin \theta, x_3, \cdots, x_n):= (\frac r2 \cos(2\theta), \frac r2 \sin (2\theta), \frac {\sqrt 3}{2} r, x_3, \cdots, x_n). 
$$
\end{remark}
 
The paper is organized as follows: In Section \ref{spreliminary} we will review some basic analytic properties of isometric 
immersions with second order derivatives. Section  \ref{developsection} is dedicated to the proof of Theorem \ref{0deve-regu}.
In Section \ref{sdensity} we will show that smooth isometric immersions are strongly dense in the space of $W^{2,2}$ 
isometric immersions from a domain of $\bbbr^n$ into $\bbbr^{n+1}$.  The proof of Lemma \ref{0difficult}, 
which is a crucial and difficult step in establishing the density result is postponed to the Appendix for the convenience of the reader.
 
\bigskip

\noindent{\bf Acknowledgments.} 
The authors would like to thank Piotr Haj{\l}asz for his interest in this problem and for many 
fruitful discussions. Z.L. was partially supported by NSF grants DMS-0900871 and DMS-0907844.
M.R.P. was supported by the NSF grants DMS-0907844 and DMS-1210258. The authors thank the referee for  
a careful reading of the original manuscript and for useful comments and suggestions about the quality 
of presentation.

\section{Preliminaries}\label{spreliminary}
Let $\Omega$ be a bounded Lipschitz domain of $\bbbr^n, n\geq
2$. We define the class of Sobolev isometric immersions from $\Omega$ to
$\bbbr^{n+1}$ as,
\begin{equation}\label{0isometric}
I^{2,2}(\Omega, \bbbr^{n+1}):=\{u\in W^{2,2}(\Omega, \bbbr^{n+1}): (\nabla u)^T
\nabla u = \rm{I} \textrm{ a.e.}\}.
\end{equation}

Note that the condition $(\nabla u)^T \nabla u = \rm{I} $
implies that $u$ is Lipschitz continuous, thus,
\begin{equation}\label{0lip}
I^{2,2}(\Omega, \bbbr^{n+1}) \subset  W^{2,2}(\Omega, \bbbr^{n+1}) \cap W^{1,\infty} (\Omega, \bbbr^{n+1}).
\end{equation}

Given $u\in I^{2,2}(\Omega, \bbbr^{n+1})$, let $u^j, 1\leq j \leq n+1,$ be the
$j$-th component of $u$ and let $u_{,i}=\partial u /\partial x_i, 1\leq i\leq
n,$ be the partial derivative of $u$ in the $\mathbf{e}_i$ direction. Throughout the paper 
we will use the same notation for all functions.

For a.e. $x\in \Omega$, consider the cross product
$$ \mathbf{n}(x)=u_{,1}(x)\times \cdots \times u_{,n}(x).$$
That is, $\mathbf{n}(x)$ is the unique unit vector orthogonal to $u_{,i}(x)$ for
all $1\leq i \leq n$ such that
$$ u_{,1}(x), \cdots, u_{,n}(x), \mathbf{n}(x) $$ form a positive basis of
$\bbbr^{n+1}$.

Note that $\mathbf{n}$ can also be identified as differential forms: consider
the $1$-form,
$$ \omega_i=\sum_{j=1}^{n+1} u^j_{,i}dx_j.$$
Then
\begin{equation}\label{0normal}
\mathbf{n}=\ast (\omega_1\wedge \cdots \wedge \omega_n).
\end{equation}
because for any $\xi\in \bigwedge^1(\bbbr^{n+1})$,
$$\langle \xi, \mathbf{n} \rangle=\langle \xi, \ast (\omega_1\wedge \cdots
\wedge \omega_n)\rangle =(-1)^n\xi\wedge \omega_1\wedge \cdots \wedge
\omega_n=(-1)^n\det[\xi, u_{,1},\cdots,u_{,n}]. $$
Since $u\in W^{2,2} \cap W^{1,\infty} (\Omega,
\bbbr^{n+1})$, it follows from (\ref{0normal}) that $\mathbf{n}\in
W^{1,2}(\Omega, \bbbr^{n+1})$. 

Since $u$ is isometric immersion, $\langle u_{,i}, u_{,j} \rangle =\delta_{ij}$ for
all $1\leq i,j\leq n$. Since  $u\in W^{2,2}(\Omega, \bbbr^{n+1})$, we can
differentiate using the product rule to obtain,
\begin{equation}\label{01}
\langle u_{,ik}, u_{,j }\rangle + \langle u_{,i}, u_{,jk} \rangle =0 \quad
\textrm{a.e.}
\end{equation}
Permutation of indices $i,j,k$ yields,
\begin{equation}\label{02}
\langle u_{,ij}, u_{,k}\rangle + \langle u_{,i}, u_{,kj} \rangle =0 \quad
\textrm{a.e.}
\end{equation}
\begin{equation}\label{03}
\langle u_{,ki}, u_{,j }\rangle + \langle u_{,k}, u_{,ji} \rangle =0 \quad
\textrm{a.e.}
\end{equation}
Using the fact that $u_{,ij}=u_{,ji}$ for all $i,j$, we add (\ref{01}) and
(\ref{02}), then subtract (\ref{03}) to obtain,
\begin{equation}\label{04}
\langle u_{,i}, u_{,jk}\rangle =0 \quad \textrm{a.e.} \quad \textrm{for all}
\quad 1\leq i,j,k\leq n.
\end{equation}
Since for a.e. points in the domain, $\mathbf{n}, u_{,1},,,u_{,j}$ form an orthonormal basis
for $\bbbr^{n+1}$, we can write,
$$u_{,jk}=\sum_{i=1}^n \langle u_{,jk}, u_{,i}\rangle u_{,i} +\langle u_{,jk},
\mathbf{n}\rangle \mathbf{n}.$$
(\ref{04}) then gives,
\begin{equation}\label{05}
u_{,jk}=\langle u_{,jk}, \mathbf{n}\rangle \mathbf{n} \quad \textrm{a.e.} \quad
\textrm{for all} \quad 1\leq j,k\leq n.
\end{equation}
Note that $A_{jk}:=\langle u_{,jk}, \mathbf{n}\rangle$ is the element in row $j$
and column $k$ of the second fundamental form $A$, which is a symmetric $n\times n$
matrix.
In particular, (\ref{05}) holds for each component of $u_{,jk}$ and
$\mathbf{n}$, i.e.,
$$u^\ell_{,jk}=A_{jk}\mathbf{n}^\ell \quad \textrm{for all} \quad 1\leq \ell
\leq n+1, \quad 1\leq j,k \leq n.$$
Thus, the Hessian of $u^\ell$ satisfies,
\begin{equation}\label{06}
\nabla^2 u^\ell =\mathbf{n}^\ell A, \quad 1\leq \ell \leq n+1.
\end{equation}

\begin{lemma}\label{0l1}
The second fundamental form $A\in M^{n\times n}$ has the following properties,
\begin{equation}\label{0id0}
\frac{\partial A_{ij}}{\partial x_k}=\frac{\partial A_{ik}}{\partial x_j} \quad
\textrm{in distributional sense for all} \quad 1\leq i,j,k \leq n,
\end{equation}
and
\begin{equation}\label{0id00}
A_{ij}A_{kl}-A_{il}A_{kj}=0 \quad \textrm{for all} \quad 1\leq i,j,k,l \leq n.
\end{equation}
\end{lemma}
{\em Proof.}
For a smooth immersion $v:\Omega\rightarrow \bbbr^{n+1}$, not necessarily
isometric, let $g_{ij}=\langle v_{,i}, v_{,j}\rangle$ be the first fundamental
forms, then by differentiating $g_{ij}$ twice,
$$g_{ij,kl}=\langle v_{,ikl}, v_{,j}\rangle +\langle v_{,ik},
v_{,jl}\rangle+\langle v_{,il}, v_{,jk}\rangle+\langle v_{,i},
v_{,jkl}\rangle.$$
The summation over the proper permutations of $i,j,k,l$ yields
\begin{equation}\label{0id1}
g_{ij,kl}+g_{kl,ij}-g_{il,kj}-g_{kj,il}=-2\langle v_{,ij},v_{,kl} \rangle +
2\langle v_{,il},v_{,kj} \rangle.
\end{equation}

Given any other smooth immersion $w:\Omega\rightarrow \bbbr^{n+1}$, the
following identity is also obvious,
\begin{equation}\label{0id2}
\langle v_{,ij}, w\rangle_{,k}-\langle v_{,ik}, w\rangle_{,j}=\langle v_{,ij},
w_{,k}\rangle-\langle v_{,ik}, w_{,j}\rangle.
\end{equation}

Now we let a sequence of smooth immersions $u_m\rightarrow u$ in
$W^{2,2}(\Omega, \bbbr^{n+1})$ and $\mathbf{n}_m\rightarrow \mathbf{n}$ in
$W^{1,2}(\Omega, \bbbr^{n+1})$. Writing the left hand sides of (\ref{0id1}) and
(\ref{0id2}) as distributional derivatives and passing to the limit we get,
\begin{equation}\label{0id3}
0=-2\langle u_{,ij},u_{,kl} \rangle + 2\langle u_{,il},u_{,kj} \rangle.
\end{equation}
because $\langle u_{,i}, u_{,j}\rangle =\delta_{ij}$ for all $i,j$. In addition,
since $\mathbf{n}$ is a unit vector, $\langle\mathbf{n}_{,k}, \mathbf{n} \rangle
=0$. Then by (\ref{05}), $\langle u_{,ij}, \mathbf{n}_{,k}\rangle=0$ for all
$i,j,k$, thus,
\begin{equation}\label{0id4}
\langle u_{,ij}, \mathbf{n}\rangle_{,k}-\langle u_{,ik},
\mathbf{n}\rangle_{,j}=0
\end{equation}
The two identities in the lemma follow easily from $A_{ij}=\langle u_{,ij},
\mathbf{n} \rangle$, (\ref{0id3}), and (\ref{0id4}). The proof is complete.
\hfill $\Box$

\begin{corollary}\label{0c1}
The second fundamental form $A$ satisfies $\rm{rank}\, A\leq 1$ and $A$ is symmetric a.e. in $\Omega$. Moreover, the Hessian of 
each component of $u$ satisfies ${\rm rank}\, \nabla^2 u^\ell \leq 1$ for all $1\leq \ell \leq n+1$ a.e. on $\Omega$.  
\end{corollary}
{\em Proof.} By identity (\ref{0id00}), all $2\times 2$ minors of $A$ vanish, hence the rank of $A$ is
less than or equal to $1$. By (\ref{06}), ${\rm rank}\,\nabla^2 u^\ell \leq {\rm rank}\, A\leq 1$ and $A$ is symmetric a.e. 
since $\nabla^2 u^\ell$ is symmetric a.e. The proof is complete.
\hspace*{\fill} $\Box$

\section{Developability and Regularity}\label{developsection}
Our first main result- Theorem \ref{0deve-regu}- follows from the following proposition:
\begin{proposition}\label{0induction}
Let $u\in I^{2,2}(\Omega, \bbbr^{n+1})$, where $\Omega$ is a bounded Lipschitz domain in
$\bbbr^n$. Let $A$ be the second fundamental form of $u$. Let $P_k$ be a $k$-dimensional plane of $\bbbr^n$, $k\leq n$. 
Suppose on $P_k\cap \Omega$ we have the following properties, 
\begin{enumerate}
 \item There exists a sequence of smooth functions $u^\epsilon$ defined in the domain $\Omega$ such that
\begin{equation*}
\int_{P_k\cap \Omega}
|u^\epsilon-u|^2+|\nabla u^\epsilon-\nabla u|^2+ |\nabla^2 u^\epsilon-\nabla^2
u|^2 d\mathcal{H}^{k}\rightarrow 0.
\end{equation*} 
Here $\nabla u^\epsilon$, $\nabla u$, $\nabla^2 u^\epsilon$ and $\nabla^2 u$ denote the first and second full gradients with respect to the domain $\Omega$. 
 \item The full gradient $\nabla u$ satisfies $\nabla u^T \nabla u={\rm I}$ $\mathcal{H}^k$-a.e. on $P_k\cap \Omega$. 
 \item $\nabla^2 u^\ell =\mathbf{n}^\ell A$ for each $1\leq \ell
\leq n+1$ $\mathcal{H}^k$-a.e. on $P_k\cap \Omega$.
 \item ${\rm rank}\, A \leq 1$ and $A$ is symmetric $\mathcal{H}^k$-a.e. on $P_k\cap \Omega$. 
\end{enumerate}
Then $u\in C_{\rm loc}^{1,1/2}(P_k, \bbbr^{n+1})$. Moreover, for
every $x\in P_k\cap \Omega$, either $\nabla u$ is constant in a neighborhood in $P_k\cap \Omega$ of $x$, or
there exists a unique $(k-1)$-dimensional hyperplane $P^x_{k-1}\ni x$ of $P_k$ such that $\nabla u$ is constant on the connected 
component of $x$ in $P_{k-1}^x\cap \Omega$. 
\end{proposition}

The proof of this proposition is based on induction on lower dimensional slices. Before we prove Proposition \ref{0induction},
 we will show that it implies Theorem \ref{0deve-regu}.

{\em Proof of Theorem \ref{0deve-regu}.}
We simply take $k=n$ in Proposition \ref{0induction}, in which case $P_n\cap \Omega=\Omega$. Since $u\in W^{2,2}(\Omega, \bbbr^{n+1})$, 
the convolution of $u$ with the standard mollifier $u^\epsilon$ apparently satisfies assumption (1). By the fact that
$u\in I^{2,2}(\Omega, \bbbr^{n+1})$, $\nabla u^T \nabla u={\rm I}$ a.e. in $\Omega$, which is property (2). Property (3) follows 
from equation (\ref{06}) and property (4) follows from Corollary \ref{0c1} . Therefore, all the assumptions of Proposition \ref{0induction} 
are satisfied, and hence the conclusion of Theorem \ref{0deve-regu} follows from the conclusion of Proposition \ref{0induction}. The proof is complete. 
\hspace*{\fill} $\Box$

\begin{corollary}
Let $u\in I^{2,2}(\Omega, \bbbr^{n+1})$, where $\Omega$ is a Lipschitz domain in
$\bbbr^n$. Then for every $k$ dimensional slice $P_k\cap \Omega$, $\nabla u$ is constant either on $k$ dimensional neighborhoods 
of $P_k\cap \Omega$, or constant on $k-1$-dimensional slice of $P_k\cap \Omega$. 
\end{corollary}
{\em Proof.}
Since assumptions (1)-(4) of Proposition \ref{0induction} are satisfied a.e. in $\Omega$. By Fubini Theorem, assumptions 
(1)-(4) also holds in a.e. $k$-dimensional slice. Thus the conclusion of Proposition \ref{0induction} holds for a.e. $k$-dimensional slices. 
Since $\nabla u$ is continuous, by a simple approximation argument, it holds on \textit{every} $k$-dimensional slices. 
The proof is complete. 
\hspace*{\fill} $\Box$

Assumptions (2) (3) and (4) regard the properties of isometric immersions, while (1) can be formulated for any general Sobolev function. 
This latter assumption is necessary for allowing the use of the chain rule which involves the full 
gradient even in lower dimensional slices. To be precise, we prove the following lemma which will play an important role 
everywhere in the proof of Proposition \ref{0induction}. 
\begin{lemma}[Chain Rule]\label{0chain}
Let $\Psi \in W^{1,2}(\Omega, \bbbr^N)$, $N\geq 1$. Let $\Sigma\subset \Omega$ be a $k$-dimensional flat domain. 
Suppose that there exist a
sequence of smooth functions $\Psi^\epsilon \in C^\infty (\Omega, \bbbr^N)$ such
that
\begin{equation}\label{0good}
\int_{\Sigma} |\Psi^\epsilon-\Psi|^2+|\nabla \Psi^\epsilon-\nabla
\Psi|^2\,d\mathcal{H}^k\rightarrow 0,
\end{equation}
where $\nabla \Psi$ denotes the full gradient with respect to the domain
$\Omega$. Let $\mathbf{v}$ be any directional vector tangent to $\Sigma$, then the chain
rule,
$$\frac{d}{dt}\bigl|_{t=0}\Psi(\cdot+t\mathbf{v})=\nabla \Psi \mathbf{v}$$
holds in the weak sense over the domain $\Sigma$. In particular,
$$\Psi\in W^{1,2}(\Sigma, \bbbr^N). $$
\end{lemma}
{\em Proof.}
Let $\phi\in C_0^\infty (\Sigma)$, then,
$$\int_{\Sigma}
\frac{d}{dt}\bigl|_{t=0}\Psi^\epsilon(x+t\mathbf{v})\phi(x)d\mathcal{H}^k=-\int_
{\Sigma}\Psi^\epsilon(x)\frac{d}{dt}\bigl|_{t=0}\phi(x+t\mathbf{v})d\mathcal{H}
^k.$$
Since $\Psi^\epsilon$ is smooth in $\Omega$, we have,
$$\int_{\Sigma}
\frac{d}{dt}\bigl|_{t=0}\Psi^\epsilon(x+t\mathbf{v})\phi(x)d\mathcal{H}^k=\int_{
\Sigma}\nabla \Psi^\epsilon(x)\mathbf{v}\phi(x)d\mathcal{H}^k.$$
By (\ref{0good}) we pass to the limit to conclude that,
$$\int_{\Sigma}\nabla
\Psi(x)\mathbf{v}\phi(x)d\mathcal{H}^k=-\int_{\Sigma}\Psi(x)\frac{d}{dt}\bigl|_{
t=0}\phi(x+t\mathbf{v})d\mathcal{H}^k.$$
Thus the chain rule as stated in the Lemma hold in the weak sense over the
domain $\Sigma$. The proof is complete.
\hspace*{\fill} $\Box$

\begin{remark}
Note that the above lemma involves the \textit{full} gradient of $\Psi$. The
assumption $\Psi\in W^{1,2}(\Sigma, \bbbr^N)$ by itself is not enough to conclude the
chain rule.
\end{remark}

\subsection{Base case- $2$-dimensional slices.}
Suppose for a $2$-dimensional plane $P_2$ all the assumptions (1)-(4) in Proposition \ref{0induction} are satisfied. 
Without loss of generality, we can assume $P_2$ is parallel to the space spanned by $\mathbf{e}_1$ and $\mathbf{e}_2$. 
Indeed, it is easy to see that assumption (1)-(4) in Proposition \ref{0induction} are invariant under rotating the coordinate 
system. We denote $P_2$ by $P_{\mathbf{e}_1\mathbf{e}_{2}}$ to remind ourselves of this fact. 

Let $f=\nabla u^\ell\in W^{1,2}_{\rm loc} (\Omega, \bbbr^n)$ for some arbitrary $1\leq\ell \leq n+1$. Define,
$$ g:=(f^{1},f^2)|_{P_{\mathbf{e}_1\mathbf{e}_{2}}\cap
\Omega} \in W^{1,2}(P_{\mathbf{e}_1\mathbf{e}_{2}}\cap
\Omega,
\bbbr^{2}).$$
\begin{lemma}\label{0rigid}
Let $f^\epsilon:\Omega\to \bbbr^n$ be a smooth sequence converging strongly to $f$ in $W^{1,2}(\Omega, \bbbr^n)$. 
Let $C$ be a line segment in $P_{\mathbf{e}_1\mathbf{e}_{2}}\cap
\Omega$ such that
\begin{equation}\label{0convergent}
\int_C|f^\epsilon-f|^2+|\nabla f^\epsilon-\nabla
f|^2d\mathcal{H}^{1}\rightarrow 0,
\end{equation}
$\mathrm{rank}\,\nabla f\leq 1$ and $\nabla f$ is symmetric for
$\mathcal{H}^{1}$-a.e. points on $C$. Then if $g$ is constant on $C$, so is
$f$.
\end{lemma}
{\em Proof.}
Let $\mathbf{v}$ be the unit directional vector of $C$. Since $\mathbf{v}$ is a linear combination of $\mathbf{e}_1$ and
$\mathbf{e}_{2}$,
$$\mathbf{v}=(\mathbf{v}_1,\mathbf{v}_{2}, 0,\cdots,0).$$
Let $\tilde{\mathbf{v}}=(\mathbf{v}_1,\mathbf{v}_{2})$, then the first two components of $f$ satisfy $\nabla f^1\cdot 
\mathbf{v}=\nabla g^1 \cdot \tilde{\mathbf{v}}$ a.e. on $C$ and $\nabla f^2\cdot \mathbf{v}=\nabla g^2 \cdot \tilde{\mathbf{v}}$ a.e. on $C$. 

Since $f$ satisfies the
assumption of Lemma \ref{0chain}, the chain rule,
$$\frac{d}{dt}\bigl|_{t=0}f(\cdot+t\mathbf{v})=(\nabla f ) \mathbf{v}$$
holds in the weak sense on $C$. In particular, it holds for it first two component $f^1$ and $f^2$ and, of course, $g$. 

As $g$ is constant on $C$,
\begin{equation}\label{010}
0=\frac{d}{dt}\bigl|_{t=0}g(\cdot+t\tilde{\mathbf{v}})
\end{equation}
in the weak sense. Hence,
$$(\nabla g) \tilde{\mathbf{v}}=0 \quad \textrm{a.e. on } C.$$
This implies, 
$$\nabla f^1\cdot \mathbf{v}=0 \quad \textrm{and} \quad \nabla f^2\cdot \mathbf{v}=0 \quad \textrm{a.e. on } C.$$

For $z\in C$ such that $\nabla f^1(z)\cdot \mathbf{v}=0$ and $\nabla f^2(z)\cdot \mathbf{v}=0$, $\mathrm{rank}\,\nabla
f(z)\leq 1$ and $\nabla f(z)$ is symmetric, we have two cases: 1)$\nabla
f^1(z)\neq 0$ or $\nabla
f^2(z)\neq 0$, 2)$\nabla f^1(z)=\nabla f^2(z) =0$ In the first case, we can assume with loss of
generality that $\nabla f^1(z)\neq 0$. Therefore, $\mathrm{rank}\,\nabla
f(z)= 1$ and
$$\nabla f^i(z)=a^i_z \nabla f^1(z)\quad \textrm{for all } i \geq 1.$$
It then follows that
$$
\nabla f^{i}(z) \cdot \mathbf{v}=a^i_z\bigl(\nabla f^1(z)\cdot
\mathbf{v}\bigr)=0 \quad \textrm{for all } i\geq 1.
$$

In the second case, by symmetry,
$$
f^i_{,j}(z)=f^j_{,i}(z)=0, \quad \textrm{for } j=1,2, \textrm{ and } i=1,\cdots,n.
$$ As $\mathbf{v}=(\mathbf{v}_1,\mathbf{v}_{2}, 0,\cdots,0)$,
$$\nabla f^{i}(z)\cdot \mathbf{v}=0 \quad \textrm{for all } i=1,\cdots,n $$
Therefore, in either cases, we have proved
$$\nabla f^i\cdot \mathbf{v}=0 \quad \textrm{a.e. on } C \quad \textrm{for all } i=1,\cdots,n.$$
Therefore, $f$ is constant on $C$ by
the chain rule in (\ref{010}). The proof is complete.
\hspace*{\fill} $\Box$

\begin{corollary}\label{0rigid1}
If $g$ is constant on a $2$-dimensional region $U$ in $P_{\mathbf{e}_1
\mathbf{e}_{2}}\cap \Omega$, $f$ is constant on $U$ as well.
\end{corollary}
{\em Proof.}
Observe that if $U$ is a $2$-dimensional region of $P_{\mathbf{e}_1
\mathbf{e}_{2}}\cap \Omega$, which has strictly positive $2$-dimensional Hausdorff
measure, then the assumptions (1) and (4) of Proposition \ref{0induction} imply,
\begin{equation*}
\int_{U}|f^\epsilon-f|^2+|\nabla f^\epsilon-\nabla f|^2
d\mathcal{H}^{2}\rightarrow 0,
\end{equation*}
$\mathrm{rank}\,\nabla f\leq 1$ and $\nabla f$ is symmetric for
$\mathcal{H}^{2}$ a.e. points on $U$. Thus the same argument for line segments in Lemma \ref{0rigid} gives for
any directional vector $\mathbf{v}$ of $U$, $\nabla f^i\cdot \mathbf{v}=0$ a.e. on $U$ for all $i=1,\cdots,n$, hence
the chain rule implies $f$ is constant on $U$. The proof is complete.
\hspace*{\fill} $\Box$

\begin{lemma}\label{0continuity}
Suppose assumptions (1)-(4) of Proposition \ref{0induction} are
satisfied on the two dimensional region $P_{\mathbf{e}_1\mathbf{e}_{2}}\cap \Omega$. 
Let $f=\nabla u^\ell\in W^{1,2}_{\rm loc} (\Omega, \bbbr^n)$ for some
arbitrary $1\leq\ell \leq n+1$. Then the restriction $f\in C_{\rm
  loc}^{0,1/2}(P_{\mathbf{e}_1\mathbf{e}_{2}}\cap \Omega, \bbbr^n)$. 
Moreover, for
every point $x \in P_{\mathbf{e}_1\mathbf{e}_{2}}\cap \Omega$, either there exists a neighborhood in $ P_{\mathbf{e}_1\mathbf{e}_{2}}\cap \Omega$ of $x$, or a unique line
segment in $P_{\mathbf{e}_1\mathbf{e}_{2}}\cap \Omega$ passing through $x$ and joining $\partial \Omega$ at both ends, on which
$f$ is constant.
\end{lemma}

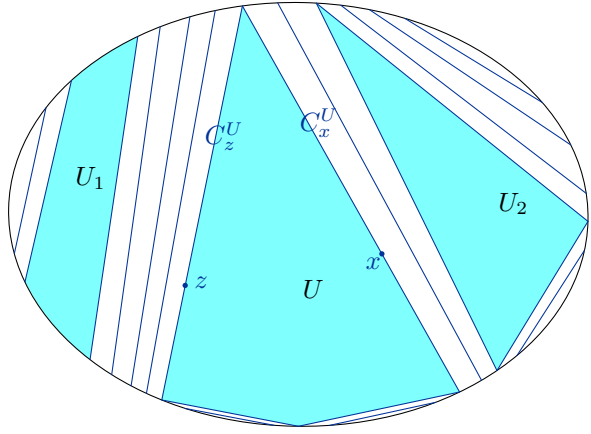
\begin{figure}[ht]
 \centering
 \begin{tikzpicture}[line cap=round,line join=round,>=triangle
45,x=1.0cm,y=1.0cm]
\clip(-4.3,0) rectangle (4.18,6.26);
\fill[line width=0pt,color=qqffff,fill=qqffff,fill opacity=0.5] (-3.11,4.869) --
(-3.729,2.186) -- (-3.532,1.846) -- (-3.315,1.57) -- (-3.092,1.343) --
(-2.864,1.152) -- (-2.229,5.453) -- (-2.522,5.298) -- (-2.819,5.105) -- cycle;
\fill[line width=0pt,color=qqffff,fill=qqffff,fill opacity=0.5] (-0.838,5.858)
-- (-1.91,0.616) -- (-0.094,0.266) -- (2.049,0.722) -- cycle;
\fill[line width=0pt,color=qqffff,fill=qqffff,fill opacity=0.5] (0.144,5.896) --
(2.546,1.007) -- (3.756,2.992) -- cycle;
\draw [rotate around={-1.2:(-0.095,3.085)},color=qqqqqq] (-0.095,3.085) ellipse
(3.852cm and 2.819cm);
\draw [color=qqttzz] (-3.449,4.504)-- (-3.87,2.563);
\draw [color=qqttzz] (-1.931,5.581)-- (-2.576,0.952);
\draw [color=qqttzz] (-1.536,5.714)-- (-2.318,0.804);
\draw [color=qqttzz] (-0.37,5.9)-- (2.341,0.877);
\draw [color=qqttzz] (1.321,5.693)-- (3.15,4.573);
\draw [color=qqttzz] (2.813,1.208)-- (3.709,2.607);
\draw [color=qqttzz] (0.321,0.278)-- (1.729,0.584);
\draw [color=qqttzz] (-2.114,0.704)-- (-1.2,5.796);
\draw [color=qqttzz](-0.2,4.62) node[anchor=north west] {$C_x^{U}$};
\draw [color=qqqqqq](2.44,3.5) node[anchor=north west] {$U_2$};
\draw [color=qqttzz] (-3.11,4.869)-- (-3.729,2.186);
\draw [color=qqttzz] (-2.229,5.453)-- (-2.864,1.152);
\draw [color=qqttzz] (-0.838,5.858)-- (-1.91,0.616);
\draw [color=qqttzz] (-0.838,5.858)-- (2.049,0.722);
\draw [color=qqttzz] (-1.91,0.616)-- (-0.094,0.266);
\draw [color=qqttzz] (-0.094,0.266)-- (2.049,0.722);
\draw [color=qqttzz] (0.144,5.896)-- (2.546,1.007);
\draw [color=qqttzz] (0.144,5.896)-- (3.756,2.992);
\draw [color=qqttzz] (-1.593,0.502)-- (-0.417,0.279);
\draw [color=qqqqqq](-3.18,3.84) node[anchor=north west] {$U_1$};
\draw [color=qqqqqq](-0.16,2.34) node[anchor=north west] {$U$};
\draw [color=qqttzz] (0.852,5.808)-- (3.527,4.009);
\draw [color=qqttzz] (0.47,5.868)-- (3.733,3.364);
\draw [color=qqttzz] (3.756,2.992)-- (2.546,1.007);
\draw [color=qqttzz](-1.46,4.44) node[anchor=north west] {$C_z^{U}$};
\draw [color=qqttzz](-1.6,2.38) node[anchor=north west] {$z$};
\draw [color=qqttzz](0.68,2.64) node[anchor=north west] {$ x $};
\begin{scriptsize}
\fill [color=qqttzz] (1.016,2.56) circle (1.0pt);
\fill [color=qqttzz] (-1.599,2.14) circle (1.0pt);
\end{scriptsize}
\end{tikzpicture}
 \caption{Inverse image of $g$ in $P_{\mathbf{e}_1\mathbf{e}_{2}}\cap \Omega$}
 \label{fig1}
\end{figure}

{\em Proof.} The proof is divided into seven steps. 

\textbf{Step 0.} Preliminary set up: by assumption (4) of Proposition \ref{0induction}, $\nabla f$ satisfies $\textrm{rank}\, \nabla f\leq 1$ and $\nabla f=\nabla^2 u^\ell$ is 
symmetric a.e. on $P_{\mathbf{e}_1\mathbf{e}_{2}}\cap \Omega$. Therefore, 
$g:=(f^{1},f^2)|_{P_{\mathbf{e}_1\mathbf{e}_{2}}\cap
\Omega} \in W^{1,2}(P_{\mathbf{e}_1\mathbf{e}_{2}}\cap
\Omega,
\bbbr^{2})$
also satisfies $\textrm{rank}\, \nabla g \leq 1$ and $\nabla g$ is symmetric a.e. on $P_{\mathbf{e}_1\mathbf{e}_{2}}\cap \Omega$. 
We employ \cite[Prop.1]{Pak}, which is cited above as Proposition \ref{0pakzad2}. The function $g$ satisfies the assumption of this 
proposition on the domain $P_{\mathbf{e}_1\mathbf{e}_{2}}\cap \Omega$ and hence the conclusions holds true for $g$. Suppose
$g$ is constant on some maximal connected neighborhood $U \subset P_{\mathbf{e}_1\mathbf{e}_{2}}\cap \Omega$, by continuity of $g$, it is also constant on
its closure $\overline{U}\cap \Omega$. 

\textbf{Step 1.} We claim that the boundary of $U$ only consists of line segments joining the boundary and none of these line segments intersect inside $\Omega$. Indeed, if $x\in \partial U\cap \Omega$, then $x$ is not
contained in a constancy neighborhood of $g$, therefore by Proposition \ref{0pakzad2},
there exists a unique line segment $C_x^U\subset P_{\mathbf{e}_1\mathbf{e}_{2}}\cap \Omega$ passing through $x$ and joining 
$\partial \Omega$ at both ends on which $g$ is constant,
which implies $\partial U\cap \Omega\subset \bigcup_{x\in \partial U\cap
\Omega}C^U_x$. Moreover, for $x,z\in \partial U\cap \Omega$, $C^U_x=C^U_z$ if $z\in C^U_x$
and $C^U_x\cap C^U_z\cap \Omega=\emptyset$ if $z\notin C^U_x$  (Figure \ref{fig1}). This follows from the fact that if $g$ 
is constant on two such intersecting segments, it must be constant on
their convex hull inside $\Omega$ too.  On the other hand, suppose $g$ is constant on
some line segment $C_x^U$ passing through $x\in \partial U\cap \Omega$ and joining $\partial \Omega$ at both end, since $g$ is constant on
$\overline{U}$ and $C^U_x$, which intersect at
$x$, it must be constant on the convex hull of $\overline{U}$ and $C^U_x$ inside
$\Omega$. But $U$ is maximal, hence $ \bigcup_{x\in
\partial U\cap \Omega}C^U_x\subset \partial U\cap \Omega$. Therefore, 
$$\partial U\cap \Omega= \bigcup_{x\in \partial U\cap \Omega}C^U_x.$$

\begin{figure}[ht]
 \centering
 \begin{tikzpicture}[line cap=round,line join=round,>=triangle
45,x=1.0cm,y=1.0cm]
\clip(-4.3,-0.02) rectangle (4.2,6.3);
\fill[line width=0pt,color=qqffff,fill=qqffff,fill opacity=0.5] (-2.344,5.378)
-- (-3.061,1.335) -- (-2.929,1.223) -- (-2.745,1.084) -- (-2.546,0.952) --
(-2.313,0.82) -- (-2.142,0.735) -- (-1.868,5.588) -- (-2.045,5.518) -- cycle;
\fill[line width=0pt,color=qqffff,fill=qqffff,fill opacity=0.5] (-0.623,0.315)
-- (-0.343,5.884) -- (-0.461,5.878) -- (-1.299,0.436) -- cycle;
\fill[line width=0pt,color=qqffff,fill=qqffff,fill opacity=0.5] (0.382,5.861) --
(-0.009,0.283) -- (0.41,0.302) -- cycle;
\fill[line width=0pt,color=qqffff,fill=qqffff,fill opacity=0.5] (0.913,0.371) --
(0.629,5.83) -- (0.567,5.838) -- (0.704,0.337) -- cycle;
\fill[line width=0pt,color=qqffff,fill=qqffff,fill opacity=0.5] (0.847,5.792) --
(1.439,0.501) -- (1.306,0.463) -- cycle;
\fill[line width=0pt,color=qqffff,fill=qqffff,fill opacity=0.5] (1.711,5.54) --
(2.53,1.015) -- (3.691,2.583) -- cycle;
\draw [rotate around={-1.2:(-0.095,3.085)},color=qqqqqq] (-0.095,3.085) ellipse
(3.84cm and 2.802cm);
\draw [color=qqttzz] (-3.423,4.515)-- (-3.846,2.525);
\draw [color=qqttzz] (-3.132,4.829)-- (-3.599,1.974);
\draw [color=qqttzz] (-2.787,5.109)-- (-3.323,1.599);
\draw [color=qqttzz] (-2.344,5.378)-- (-3.061,1.335);
\draw [color=qqttzz] (-1.868,5.588)-- (-2.142,0.735);
\draw [color=qqttzz] (-1.444,5.721)-- (-1.874,0.619);
\draw [color=qqttzz] (-0.936,5.827)-- (-1.611,0.526);
\draw [color=qqttzz] (-0.461,5.878)-- (-1.299,0.436);
\draw [color=qqttzz] (-0.343,5.884)-- (-0.623,0.315);
\draw [color=qqttzz] (-0.182,5.887)-- (-0.395,0.295);
\draw [color=qqttzz] (0.105,5.881)-- (-0.223,0.286);
\draw [color=qqttzz] (0.567,5.838)-- (0.704,0.337);
\draw [color=qqttzz] (0.629,5.83)-- (0.913,0.371);
\draw [color=qqttzz] (0.847,5.792)-- (1.306,0.463);
\draw [color=qqttzz] (1.139,5.726)-- (1.757,0.612);
\draw [color=qqttzz] (1.422,5.644)-- (2.142,0.786);
\draw [color=qqttzz] (1.711,5.54)-- (2.53,1.015);
\draw [color=qqttzz] (2.336,5.23)-- (3.729,2.8);
\draw [color=qqttzz] (-1.299,0.436)-- (-0.623,0.315);
\draw [color=qqttzz] (0.382,5.861)-- (-0.009,0.283);
\draw [color=qqttzz] (-0.009,0.283)-- (0.41,0.302);
\draw [color=qqttzz] (0.41,0.302)-- (0.382,5.861);
\draw [color=qqttzz] (0.847,5.792)-- (1.439,0.501);
\draw [color=qqttzz] (1.711,5.54)-- (3.691,2.583);
\draw [color=qqttzz] (2.53,1.015)-- (3.691,2.583);
\draw [color=qqttzz] (2.806,1.221)-- (3.586,2.252);
\draw [color=qqttzz] (2.922,4.788)-- (3.733,3.262);
\draw [color=qqqqqq] (1.22,3.16) circle (1.908cm);
\draw [color=qqqqqq] (1.22,3.16)-- (2.639,4.435);
\draw [color=qqqqqq](-2.68,3.34) node[anchor=north west] {$U_1$};
\draw [color=qqqqqq](-0.3,3.42) node[anchor=north west] {$U_m$};
\draw [color=qqqqqq](0.42,3.3) node[anchor=north west] {$\cdots$};
\draw [color=qqqqqq](1.18,3.42) node[anchor=north west] {$x_0$};
\draw [color=qqqqqq](2.16,4.2) node[anchor=north west] {$\delta$};
\draw [color=qqttzz](-0.38,5.9) node[anchor=north west] {$C_{x_1}^{U_m}$};
\draw [color=qqttzz](0.94,5.9) node[anchor=north west] {$C_{x_2}^{U_m}$};
\draw [color=qqttzz](0.2,1.44) node[anchor=north west] {$x_1$};
\draw [color=qqttzz](0.76,1.42) node[anchor=north west] {$x_2$};
\draw [->,>=stealth',color=qqqqqq] (0.225,5.225) -- (0.583,5.212);
\draw [->,>=stealth',color=qqqqqq] (1.04,5.36) -- (0.651,5.398);
\begin{scriptsize}
\fill [color=qqqqqq] (1.22,3.16) circle (1.0pt);
\fill [color=qqttzz] (0.686,1.074) circle (1.0pt);
\fill [color=qqttzz] (0.876,1.085) circle (1.0pt);
\end{scriptsize}
\end{tikzpicture}
 \caption{$B^2(x_0, \delta)$ intersects $\partial U_m$ at two line segments}.
 \label{fig3}
\end{figure}
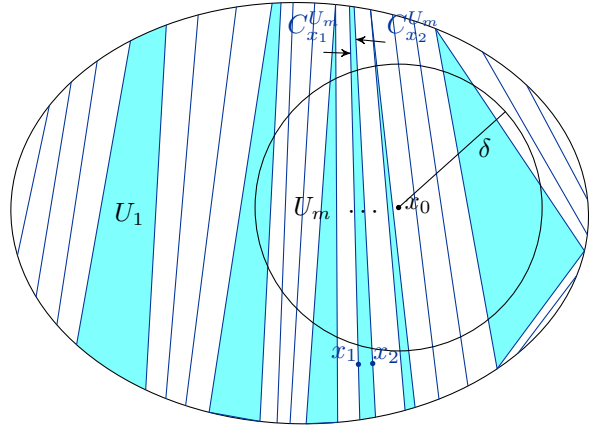

\textbf{Step 2.} We claim that we can choose
small enough $\delta>0$ so that for any region $U$ on which $g$ is constant, the
$2$-dimensional ball $B^2(x_0, \delta)\subset P_{\mathbf{e}_1
\mathbf{e}_{2}}\cap \Omega$ intersects $\partial U$ at
no more than \textit{two} line segments belonging to $\partial U$. Indeed, let $x_0\in P_{\mathbf{e}_1\mathbf{e}_{2}}\cap \Omega$ be such that $g$ is not constant in a neighborhood of $x_0$. We use the fact that for any maximal constant region $U$,  line 
segments in $\partial U$ do not intersect inside $\Omega$. If $x_0$ is at a positive distance of all constancy regions, the conclusion is trivial. 
The same is true if it lies on the boundary of one of the constancy regions and yet is positively distant from all others. Suppose therefore that 
there is a sequence of maximal constancy regions $U_m$ converging to $x_0$ 
in distance, in which case there are two line segments $C^{U_m}_{x_1}$ and $C^{U_m}_{x_2}$ in $\partial U_m$ whose angle (if they 
are nonparallel) or distance (if they are parallel) converges to zero, since both of these sequences of segments must converge to 
the same constancy segment passing through $x_0$. Then since all the other line segments in $\partial U_m$ must be arbitrarily 
close to $\partial \Omega$, we can again choose $\delta$ small enough so that $B^2(x_0, \delta)$ is away from
$\partial \Omega$ and hence it does not intersect a third line segment in
$\partial U_m$ (Figure \ref{fig3}).


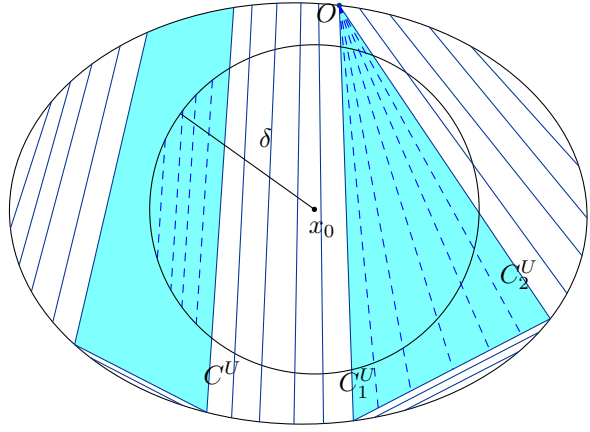
\begin{figure}[ht]
\centering
\begin{tikzpicture}[line cap=round,line join=round,>=triangle
45,x=1.0cm,y=1.0cm]
\clip(-4.3,-0.04) rectangle (4.06,6.28);
\fill[line width=0pt,color=qqffff,fill=qqffff,fill opacity=0.5] (0.351,5.833) -- (0.538,0.307) -- (3.156,1.665) -- cycle;
\fill[line width=0pt,color=qqffff,fill=qqffff,fill opacity=0.5] (-2.201,5.474) -- (-3.168,1.321) -- (-1.411,0.419) -- (-1.053,5.805) -- (-1.274,5.765) -- (-1.483,5.717) -- (-1.647,5.673) -- (-1.859,5.606) -- (-2.008,5.552) -- cycle;
\draw [rotate around={-1.2:(-0.195,3.065)},color=qqqqqq] (-0.195,3.065) ellipse (3.84cm and 2.802cm);
\draw [color=qqttzz] (-2.201,5.474)-- (-3.168,1.321);
\draw [color=qqttzz] (-1.053,5.805)-- (-1.411,0.419);
\draw [color=qqttzz] (0.351,5.833)-- (0.538,0.307);
\draw [color=qqttzz] (0.351,5.833)-- (3.156,1.665);
\draw [color=qqttzz] (0.538,0.307)-- (3.156,1.665);
\draw [color=qqqqqq] (0.02,3.12) circle (2.19cm);
\draw [dash pattern=on 3pt off 3pt,color=qqqqff] (0.351,5.833)-- (0.868,0.478);
\draw [dash pattern=on 3pt off 3pt,color=qqqqff] (0.351,5.833)-- (2.352,1.248);
\draw [dash pattern=on 3pt off 3pt,color=qqqqff] (0.351,5.833)-- (1.303,0.704);
\draw [dash pattern=on 3pt off 3pt,color=qqqqff] (0.351,5.833)-- (2.767,1.463);
\draw [dash pattern=on 3pt off 3pt,color=qqqqff] (0.351,5.833)-- (1.862,0.994);
\draw [dash pattern=on 3pt off 3pt,color=qqqqff] (-1.336,4.84)-- (-1.551,1.594);
\draw [dash pattern=on 3pt off 3pt,color=qqqqff] (-1.554,4.643)-- (-1.742,1.819);
\draw [dash pattern=on 3pt off 3pt,color=qqqqff] (-1.729,4.438)-- (-1.888,2.045);
\draw [dash pattern=on 3pt off 3pt,color=qqqqff] (-1.963,4.051)-- (-2.068,2.46);
\draw [color=qqttzz] (-0.809,5.837)-- (-1.081,0.347);
\draw [color=qqttzz] (-0.478,5.862)-- (-0.711,0.294);
\draw [color=qqttzz] (-0.159,5.867)-- (-0.252,0.264);
\draw [color=qqttzz] (0.076,5.857)-- (0.141,0.27);
\draw [color=qqttzz] (0.618,5.796)-- (3.377,2.003);
\draw [color=qqttzz] (1.022,5.711)-- (3.577,2.504);
\draw [color=qqttzz] (1.551,5.543)-- (3.644,3.005);
\draw [color=qqttzz] (2.102,5.288)-- (3.534,3.695);
\draw [color=qqttzz] (0.818,0.352)-- (3.004,1.485);
\draw [color=qqttzz] (1.226,0.448)-- (2.734,1.225);
\draw [color=qqttzz] (-2.959,1.147)-- (-1.734,0.513);
\draw [color=qqttzz] (-3.168,1.321)-- (-1.411,0.419);
\draw [color=qqttzz] (-3.476,4.552)-- (-4.021,2.866);
\draw [color=qqttzz] (-3.222,4.818)-- (-3.905,2.38);
\draw [color=qqttzz] (-2.937,5.053)-- (-3.742,2.027);
\draw [color=qqttzz] (-2.563,5.293)-- (-3.458,1.621);
\draw [color=qqqqqq](-1.58,1.22) node[anchor=north west] {$C^U$};
\draw [color=qqqqqq](0.22,1.12) node[anchor=north west] {$C_1^U$};
\draw [color=qqqqqq](2.36,2.56) node[anchor=north west] {$C_2^U$};
\draw [color=qqqqqq](-0.08,5.98) node[anchor=north west] {$O$};
\draw [color=qqqqqq](-0.18,3.08) node[anchor=north west] {$x_0$};
\draw [color=qqqqqq] (0.02,3.12)-- (-1.761,4.394);
\draw [color=qqqqqq](-0.84,4.32) node[anchor=north west] {$\delta$};
\begin{scriptsize}
\fill [color=qqttzz] (0.351,5.833) circle (1.0pt);
\fill [color=qqqqqq] (0.02,3.12) circle (1.0pt);
\end{scriptsize}
\end{tikzpicture}
\caption{Constructin of foliations.}
\label{fig4}
\end{figure}

\textbf{Step 3.} We construct a foliation of the ball $B^2(x_0, \delta)\subset P_{\mathbf{e}_1
\mathbf{e}_{2}}\cap \Omega$. Firstly, for any $x\in B^2(x_0, \delta)$, we will construct a line segment $C_x$ in $B^{2}(x_0, \delta)$ passing through $x$ and joining $\partial B^2(x_0, \delta)$ at both ends on
which $g$ is constant and $C_x\cap C_z\cap B^2(x_0, \delta)=\emptyset$ if
$z\notin C_x$. The construction is as follows: for those $x$ not contained in a constant region of $g$, this
line segment is given automatically by Proposition \ref{0pakzad2}. If $x$ is contained in a constant maximal region $U$ of $g$, then it
is constant on every line segment in $U$ that passes through it so we have to
choose the appropriate one: 1) If $B^2(x_0, \delta)$ intersect only one
line segment $C^U$ in $\Omega$ that belongs to $\partial U$, then we define
$C_x$ to be the line segments inside $B^{2}(x_0, \delta)$ passing through $x$ and parallel to
$C^U$; 2) If $B^2(x_0, \delta)$ intersects two line segments $C^U_1,
C^U_2$ in $\Omega$ that belongs to $\partial U$, let $L_1$ and $L_2$ be the two
lines that contain $C^U_1$ and $C^U_2$. If $L_1$ and $L_2$ are not
parallel, let
$O:=L_1\cap L_2$ and let $C_x$ be the segment given by the intersection with $B^2(x_0, \delta)$ of the line passing through $O$ and $x$.
If $L_1$ and $L_2$ are parallel, then we let $C_x$ be the line segment inside $B^2(x_0, \delta)$ passing through $x$ and parallel to $L_1$. (Figure \ref{fig4}). In this way, we have constructed a family of line segments $\{C_x\}_{x\in B^2(x_0,
\delta)}$ in $B^2(x_0, \delta)$ on which $g$ is constant and $C_x\cap
C_z\cap B^2(x_0, \delta)=\emptyset$ if $z\notin C_x$.

\begin{figure}[ht]
\centering
\begin{tikzpicture}[line cap=round,line join=round,>=triangle
45,x=1.0cm,y=1.0cm]
\clip(-4.3,-0.18) rectangle (2.68,6.3);
\fill[color=qqffff,fill=qqffff,fill opacity=0.5] (-3.227,4.293) -- (-3.369,3.969) -- (-3.444,3.726) -- (-3.491,3.503) -- (-3.516,3.318) -- (-3.527,3.053) -- (-3.518,2.842) -- (-3.491,2.628) -- (-3.453,2.442) -- (-3.369,2.164) -- (-3.296,1.983) -- (-3.141,1.687) -- (-2.957,1.422) -- (-2.717,1.158) -- (-2.458,0.938) -- (-2.161,5.387) -- (-2.447,5.204) -- (-2.72,4.973) -- (-2.938,4.735) -- (-3.098,4.515) -- cycle;
\fill[color=qqffff,fill=qqffff,fill opacity=0.1] (-0.184,5.63) -- (-0.071,5.597) -- (-0.258,5.649) -- cycle;
\fill[color=qqffff,fill=qqffff,fill opacity=0.5] (-0.308,5.66) -- (-0.258,5.649) -- (-0.071,5.597) -- (0.052,5.554) -- (0.184,5.501) -- (1.75,2.661) -- (1.678,2.336) -- (1.578,2.048) -- (1.403,1.701) -- (1.222,1.438) -- (0.997,1.183) -- (0.686,0.919) -- (0.455,0.769) -- (0.221,0.649) -- (-0.137,0.517) -- cycle;
\draw [color=qqqqqq] (-0.873,3.067) circle (2.654cm);
\draw [dash pattern=on 3pt off 3pt,color=qqqqff] (-2.447,5.204)-- (-2.717,1.158);
\draw [dash pattern=on 3pt off 3pt,color=qqqqff] (-2.72,4.973)-- (-2.957,1.422);
\draw [dash pattern=on 3pt off 3pt,color=qqqqff] (-2.938,4.735)-- (-3.141,1.687);
\draw [dash pattern=on 3pt off 3pt,color=qqqqff] (-3.227,4.293)-- (-3.369,2.164);
\draw [color=qqqqqq](-0.96,3.14) node[anchor=north west] {$x_0$};
\draw [color=qqttzz] (0.35,5.422)-- (1.77,3.306);
\draw [color=qqttzz] (0.184,5.501)-- (1.75,2.661);
\draw [dash pattern=on 3pt off 3pt,color=qqqqff] (0.052,5.554)-- (1.578,2.048);
\draw [dash pattern=on 3pt off 3pt,color=qqqqff] (-0.071,5.597)-- (1.222,1.438);
\draw [dash pattern=on 3pt off 3pt,color=qqqqff] (-0.184,5.63)-- (0.686,0.919);
\draw [dash pattern=on 3pt off 3pt,color=qqqqff] (-0.258,5.649)-- (0.221,0.649);
\draw [color=qqttzz] (-0.308,5.66)-- (-0.137,0.517);
\draw [color=qqttzz] (-0.671,5.713)-- (-0.612,0.425);
\draw [color=qqttzz] (-0.987,5.718)-- (-1.076,0.42);
\draw [color=qqttzz] (-1.407,5.667)-- (-1.625,0.521);
\draw [color=qqttzz] (-1.833,5.541)-- (-2.074,0.7);
\draw [color=qqttzz] (-2.161,5.387)-- (-2.458,0.938);
\draw [color=qqqqqq](-2.3,3.96) node[anchor=north west] {$\gamma([a,b])$};
\draw [shift={(-0.894,24.372)},color=qqqqqq]  plot[domain=4.646:4.663,variable=\t]({1*21.303*cos(\t r)+0*21.303*sin(\t r)},{0*21.303*cos(\t r)+1*21.303*sin(\t r)});
\draw [shift={(0.942,61.222)},color=qqqqqq]  plot[domain=4.663:4.67,variable=\t]({1*58.198*cos(\t r)+0*58.198*sin(\t r)},{0*58.198*cos(\t r)+1*58.198*sin(\t r)});
\draw [shift={(-0.709,22.185)},color=qqqqqq]  plot[domain=4.67:4.695,variable=\t]({1*19.127*cos(\t r)+0*19.127*sin(\t r)},{0*19.127*cos(\t r)+1*19.127*sin(\t r)});
\draw [shift={(-0.797,16.965)},color=qqqqqq]  plot[domain=4.695:4.724,variable=\t]({1*13.906*cos(\t r)+0*13.906*sin(\t r)},{0*13.906*cos(\t r)+1*13.906*sin(\t r)});
\draw [shift={(-0.853,22.008)},color=qqqqqq]  plot[domain=4.724:4.746,variable=\t]({1*18.949*cos(\t r)+0*18.949*sin(\t r)},{0*18.949*cos(\t r)+1*18.949*sin(\t r)});
\draw [shift={(-0.335,6.443)},color=qqqqqq]  plot[domain=4.746:5.359,variable=\t]({1*3.376*cos(\t r)+0*3.376*sin(\t r)},{0*3.376*cos(\t r)+1*3.376*sin(\t r)});
\draw [color=qqqqqq] (-3.524,3.197)-- (-2.313,3.116);
\begin{scriptsize}
\fill [color=qqqqqq] (-0.873,3.067) circle (1.0pt);
\end{scriptsize}
\end{tikzpicture}
\caption{Foliations of $B^2(x_0, \delta)$.}
\label{fig5}
\end{figure}
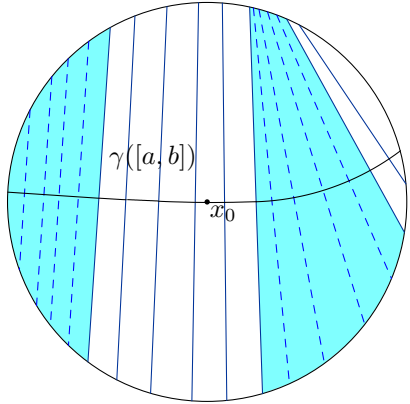

Secondly, for every $x\in
B^2(x_0, \delta)$, let $\mathbf{N}(x)$ be the vector field orthogonal to $C_x$. By making $\delta$ smaller we can make sure that 
none of the $C_x$'s intersect inside $B^k(x_0, 2\delta)$, and therefore we can choose an orientation such that 
$\mathbf{N}$ is a Lipschitz vector
field inside the ball of radius $\delta$. The ODE,
$$ \gamma'(t)= \mathbf{N}(\gamma(t)) \quad \gamma(0)=x_0 ,$$
then has a unique solution $\gamma:(a,b)\rightarrow B^2(x_0, \delta)$ for
some interval $(a,b)\subset \bbbr$ containing $0$. Moreover, if necessary by making $\delta$ smaller, 
$\cup\{C_{\gamma(t)}\}_{t\in (a,b)}=B^2(x_0, \delta)$. Therefore, $\{C_{\gamma(t)}\}_{t\in (a,b)}$ is a foliation of
$B^2(x_0, \delta)$ (Figure \ref{fig5}).

\textbf{Step 4.} We now want to show the assumptions of Lemma \ref{0rigid} are satisfied along
$C_{\gamma(t)}$ for a.e. $t\in (a,b)$. We define the function $h: B^2(x_0, \delta)\rightarrow B^2(x_0, \delta)$ as
$$h(x)=\gamma(t) \quad \textrm{if }\, x\in C_{\gamma(t)}.$$
Since none of the $C_{\gamma(t)}$ intersect inside $B^2(x_0, \delta)$, $h$
is well defined and $h$ is constant along each $C_{\gamma(t)}$, i.e. $h^{-1}(\gamma(t))=C_{\gamma(t)}$. Since $\gamma$
is Lipschitz, $h$ is Lipschitz as well. Moreover, since
$|\gamma''(t)|$ is uniformly bounded, we have the one dimensional
Jacobian $J_h > C>0$. (For the definition of general 
$k$-dimensional Jacobian see \cite[p. 88]{eg}.)

Let $E_0$ be the set of all $x\in
B^2(x_0, \delta)$ such that $\textrm{rank}\,\nabla f(x)>1$ or $\nabla f(x)$
is not symmetric. By assumption (4) of Proposition \ref{0induction} on $f$, $|E_0|=0$. 
As $h$ is Lipschitz, we can apply the general co-area formula
\cite[p. 112]{eg} to $h$ to obtain,
\begin{equation*}
\begin{array}{ll}  \ds 0 = \int_{E_0} J_h(x) dx  & \ds = \int_\gamma \mathcal{H}^1\bigl(E_0 \cap
h^{-1}(w)\bigr)d\mathcal{H}^1(w) 
\ds \\ & \ds =\int_a^b \mathcal{H}^1\bigl(E_0 \cap
h^{-1}(\gamma(t))\bigr)|\gamma'(t)|\,dt \\ & \ds
=\int_a^b \mathcal{H}^1\bigl(E_0 \cap
C_{\gamma(t)}\bigr)|\gamma'(t)|\,dt.
\end{array}
\end{equation*}
Therefore, for a.e. $t\in (a,b)$, $\mathcal{H}^1\bigl(E_0\cap
C_{\gamma(t)}\bigr)=0$ since $|\gamma'|=1$. Moreover, by change of
variable formula, if $f^\epsilon$ is a smooth approximation sequence,
\begin{equation*}
\begin{array}{ll}
\ds \int_{B^2(x_0, \delta)}\bigl(|f^\epsilon-f|^2+|\nabla f^\epsilon-\nabla
f|^2\bigr) J_h \\ & \ds \hspace{-3cm} =\int_\gamma \int_{h^{-1}(w)}\bigl(|f^\epsilon-f|^2+|\nabla f^\epsilon-\nabla
f|^2 \bigr) d\mathcal{H}^1 d\mathcal{H}^1(w)\\ & \ds \hspace{-3cm}
=\int_a^b \int_{h^{-1}(\gamma(t))}\bigl(|f^\epsilon-f|^2+|\nabla
f^\epsilon-\nabla f|^2 \bigr)d\mathcal{H}^1 |\gamma'(t)|\,dt\\ & \ds \hspace{-3cm}
=\int_a^b \int_{C_{\gamma(t)}}\bigl(|f^\epsilon-f|^2+|\nabla
f^\epsilon-\nabla f|^2 \bigr)d\mathcal{H}^1 \,dt.
\end{array}
\end{equation*}
Since $J_h$ is bounded, together with assumption (1) in Proposition \ref{0induction}, we then have for a.e.
$t\in (a,b)$,
$$\int_{C_{\gamma(t)}}\bigl(|f^\epsilon-f|^2+|\nabla f^\epsilon-\nabla f|^2
\bigr) d\mathcal{H}^1
\rightarrow 0. $$
Therefore, the assumptions of
Lemma \ref{0rigid} are satisfied along $C_{\gamma(t)}$ for a.e. $t\in (a,b)$. 

\textbf{Step 5.} We now prove the Lemma for the ball $B^2(x_0, \delta)$. Step 4 and Lemma \ref{0rigid} imply that $f$ is constant on $C_{\gamma(t)}$ for a.e. $t\in (a,b)$. By choosing an initial value for $\gamma$ arbitrary 
close to $x_0$ and applying the general co-area formula in a similar manner we can make sure that 
$f$ is of class $W^{1,2}$ on $\gamma$. Hence we conclude that $f$ is
$C^{0,1/2}$ on $\gamma$ by the Sobolev embedding theorem. Let $F$ be the set of
$t\in (a,b)$ such that $f$ is not constant along $C_{\gamma(t)}$, then $\mathcal{H}^1(F)=0$. We modify $f$ 
to be constant along $C_{\gamma(t)}$
for each $t\in F$. Note that, 
$$\mathcal{H}^2(\bigcup \{C_{\gamma(t)}:t\in F\})\leq 2\delta\sup{J_h ^{-1}} \mathcal{H}^1
(\{\gamma(t):t\in F\})= 2\delta\sup{ J_h ^{-1}} \mathcal{H}^1(F)=0 . $$
Hence $f$ is $C^{0,1/2}$ up to modification of a set of measure zero in
$B^2(x_0, \delta)$. Moreover, $f$ is constant on $C_{\gamma(t)}$ for all $t$, which foliates $B^2(x_0, \delta)$. In addition, by Corollary \ref{0rigid1}, 
$f$ is constant on every $2$-dimensional region in $B^2(x_0, \delta)$ on which $g$ is
constant. Therefore, $f$ is either constant on a line segment joining $\partial B^2(x_0, \delta)$ at both ends, or
constant on a $2$-dimensional region in $B^2(x_0, \delta)$. This proves
Lemma \ref{0continuity} for the ball $B^2(x_0, \delta)$. 

\textbf{Step 6.} Finally, we prove the lemma for the entire domain
$P_{\mathbf{e}_1\mathbf{e}_{2}}\cap \Omega$. Suppose there is some $x\in P_{\mathbf{e}_1\mathbf{e}_{2}}\cap \Omega$ that is not contained in a constant region of $f$.
Then by what we have proved, $f$ is constant on a line segment passing through $x$ and joining the boundary of $B^2(x,\delta_x)\subset P_{\mathbf{e}_1\mathbf{e}_{2}}\cap 
\Omega$ for some $\delta_x>0$. Let $\overline{y_1y_2}$ be the largest line segment containing this segment on which $f$ is constant. Suppose 
$y_1\in P_{\mathbf{e}_1\mathbf{e}_{2}}\cap \Omega$, then from what we have proved, $f$ is either constant on $2$-dimensional regions or line 
segments passing through $y_1$ and joining the boundary of $B^2(y_1, \delta_{y_1})\subset P_{\mathbf{e}_1\mathbf{e}_{2}}\cap \Omega$ for some 
$\delta_{y_1}>0$. Firstly, $y_1$ cannot be contained in a constant region of $f$, otherwise we can prolong the segment $[y_1,y_2]$. Thus, 
there must be a line segment $\overline{z_1z_2}$ passing through $y_1$ and joining the boundary of $B^2(y_1, \delta_{y_1})$ at 
both end on which $f$ is constant. Secondly, $\overline{z_1z_2}$ cannot have the same 
direction as $\overline{y_1y_2}$, otherwise, we can again prolong the segment $\overline{y_1y_2}$. Then we consider the region $\Delta$ bounded by $\overline{y_2z_1}$, 
$\overline{z_1z_2}$ and $\overline{z_2y_2}$. Since $g$ is constant on $\overline{y_1y_2}$ and  $\overline{z_1z_2}$, by Proposition \ref{0pakzad2}, $g$ must be constant on 
$\Delta$ because no line segment can join the boundary of $P_{\mathbf{e}_1\mathbf{e}_{2}}\cap \Omega$ passing through a point inside $\Delta$ without intersecting either 
$\overline{y_1y_2}$ or $\overline{z_1z_2}$ (Figure \ref{fig20}). Hence by Corollary \ref{0rigid2}, $f$ is constant on $\Delta$ as well, contradiction 
to our assumption $x$ is not contained in a constant region of $f$. The proof is complete. \hspace*{\fill} $\Box$
\begin{figure}[ht]
 \centering
 \begin{tikzpicture}[line cap=round,line join=round,>=triangle 45,x=1.0cm,y=1.0cm]
\clip(-3.7,-0.5) rectangle (4.62,5.62);
\fill[line width=0pt,color=qqffff,fill=qqffff,fill opacity=0.5] (-0.61,4.12) -- (-0.21,0.33) -- (1.37,3.8) -- cycle;
\draw [rotate around={-1.2:(0.53,2.59)}] (0.53,2.59) ellipse (3.85cm and 2.82cm);
\draw(0.08,2.12) circle (1.5cm);
\draw (0.38,3.96)-- (-0.21,0.33);
\draw(0.38,3.96) circle (1cm);
\draw (-0.61,4.12)-- (1.37,3.8);
\draw (-0.61,4.12)-- (-0.21,0.33);
\draw (-0.21,0.33)-- (1.37,3.8);
\draw (0.08,2.36) node[anchor=north west] {$x$};
\draw (-0.22,0.56) node[anchor=north west] {$y_2$};
\draw (0.38,4.2) node[anchor=north west] {$y_1$};
\draw (-0.6,4.36) node[anchor=north west] {$z_1$};
\draw (1.36,4.04) node[anchor=north west] {$z_2$};
\draw [color=qqttzz] (-2.66,4.2)-- (3.72,0.98);
\draw [color=qqttzz] (-2.89,3.92)-- (2.67,0.22);
\draw [color=qqttzz] (-2.31,4.52)-- (4.27,1.88);
\draw [color=qqttzz] (-3.11,3.55)-- (1.58,-0.14);
\draw [color=qqttzz] (-3.24,3.21)-- (0.36,-0.23);
\draw [color=qqttzz] (-1.9,4.8)-- (4.35,2.89);
\draw (2.14,0.34) node[anchor=north west] {$P_{\mathbf{e}_1\mathbf{e}_2}\cap\Omega$};
\draw (0.08,2.12)-- (-1.25,1.42);
\draw (0.38,3.96)-- (1.22,4.5);
\draw (-0.92,1.78) node[anchor=north west] {\footnotesize $\delta_x$};
\draw (0.84,4.5) node[anchor=north west] {\footnotesize $\delta_{y_1}$};
\draw (0.28,2.78) node[anchor=north west] {$\Delta$};
\begin{scriptsize}
\fill [color=black] (0.08,2.12) circle (1.0pt);
\fill [color=black] (0.38,3.96) circle (1.0pt);
\fill [color=black] (-0.21,0.33) circle (1.0pt);
\fill [color=black] (-0.61,4.12) circle (1.0pt);
\fill [color=black] (1.37,3.8) circle (1.0pt);
\end{scriptsize}
\end{tikzpicture}
 \caption{}
 \label{fig20}
\end{figure}

Now we are ready to prove Proposition \ref{0induction} for the domain $P_{\mathbf{e}_1\mathbf{e}_{2}}\cap \Omega$. Since we take $f=\nabla u^\ell$ 
for arbitrary $1\leq \ell \leq n+1$, Lemma \ref{0continuity} gives all
$\nabla u^\ell$ are continuous on $P_{\mathbf{e}_1\mathbf{e}_{2}}\cap \Omega$ and constant either on $2$-dimensional neighborhoods or line segments 
in $P_{\mathbf{e}_1\mathbf{e}_{2}}\cap \Omega$ joining $\partial \Omega$ at both ends. 
Therefore, what is left is to prove that they are constant on the \textit{same}
neighborhoods or line segments in $P_{\mathbf{e}_1\mathbf{e}_{2}}\cap \Omega$. 

Recall from equation (\ref{0normal}) that $\mathbf{n}$ is the wedge product of entries of $\nabla u$, hence is continuous. Let 
$$\Delta_\ell=\{x\in P_{\mathbf{e}_1\mathbf{e}_{2}}\cap \Omega: \mathbf{n}^\ell(x)\neq 0\}.$$
Apparently each $\Delta_\ell$ is open by continuity. Moreover, since $|\mathbf{n}|=1$ everywhere, 
$$ \bigcup_{1\leq\ell \leq n+1} \Delta_\ell= P_{\mathbf{e}_1\mathbf{e}_{2}}\cap \Omega.$$

Let $x_0\in P_{\mathbf{e}_1\mathbf{e}_{2}}\cap \Omega$, then $x_0\in \Delta_\ell$ for some $\ell$. Without loss of generality, we assume $x_0\in \Delta_1$. 
Then as the same argument as in the proof of Lemma \ref{0continuity}, there exist $B^2(x_0, \delta)\subset \Delta_1$ for some $\delta>0$, 
on which we can construct a foliation $\{C_{\gamma(t)}\}_{t\in (a,b)}$, i.e. $\cup\{C_{\gamma(t)}\}_{t\in (a,b)}=B^2(x_0, \delta)$ and 
$C_\gamma(t)\cap C_\gamma(t')\cap B^2(x_0, \delta)=\emptyset$ for $t'\neq t$. Moreover, $\nabla^2 u^1$ is constant on $C_\gamma(t)$ for every $t\in (a,b)$. 
Assumption (1) and (3) in Proposition \ref{0induction}, together with the same argument using co-area and change of variable 
formulas as in the proof of Lemma \ref{0continuity} yield for a.e. $t\in (a,b)$
$$
\int_{C_{\gamma(t)}}|\nabla u^\epsilon-\nabla u|^2+|\nabla^2 u^\epsilon-\nabla^2 u
|^2d\mathcal{H}^{1}\rightarrow 0
$$
and $\nabla^2 u^\ell =(\mathbf{n}^\ell/\mathbf{n}^1) \nabla^2 u^1, 2\leq \ell \leq n+1$, $\mathcal{H}^1$ a.e on $C_{\gamma(t)}$. 

Let $\mathbf{v}$ be the
directional vector of one such $C_{\gamma(t)}$, then the chain rule in Lemma \ref{0chain} and the fact that $\nabla u^1$ is constant on $C_{\gamma(t)}$
imply
$$0=\frac{d}{dt}\bigl|_{t=0}\nabla u^1(\cdot +t\mathbf{v})=(\nabla^2 u^1)\mathbf{v}
$$
in the weak sense in $C_{\gamma(t)}$. Therefore, 
$$\bigl(\nabla^2 u^\ell \bigr) \mathbf{v} =\frac{\mathbf{n}^\ell}{\mathbf{n}^1} (\nabla^2 u^1)\mathbf{v}=0,
\quad 2\leq \ell \leq n+1 \textrm{ a.e. on } C_{\gamma(t)}.$$
Hence again by the
chain rule in Lemma
\ref{0chain}, $\nabla u^\ell, 2\leq \ell \leq n+1$, is constant on $C_{\gamma(t)}$. 
Therefore, each $\nabla u^\ell$ is constant on $C_{\gamma(t)}$ for a.e. $t\in (a,b)$. 
Furthermore, since for each $1\leq \ell \leq n+1 $, $\nabla u^\ell$ is continuous on $P_{\mathbf{e}_1\mathbf{e}_{2}}\cap \Omega$, we conclude
that $\nabla u^\ell$ for all $1\leq \ell \leq n+1 $ are constant on all
$C_{\gamma(t)}$ that foliates $B^2(x_0,\delta)$. 
On the other hand, each $2$-dimensional region $U$ of $B^2(x_0,\delta)$ automatically
satisfies all the assumptions (1) and (3) in Proposition \ref{0induction}, hence the same argument for each $C_{\gamma(t)}$ gives that $\nabla u^\ell$ for all $2\leq \ell \leq n+1 $ is constant on the same region
on which $\nabla u^1$ is constant. This proves $\nabla u$ is either constant on $2$-dimensional regions or constant on line segments in 
$B^2(x_0, \delta)$ joining the boundary. The proof of Proposition \ref{0induction} for the domain $P_{\mathbf{e}_1\mathbf{e}_{2}}\cap \Omega$ 
follows from exactly the same argument as the last step of the proof of Lemma \ref{0continuity}. The proof for the base case is complete. 
\hspace*{\fill} $\Box$

\subsection{Inductive step- $k$-dimensional slices}

In this subsection, we will prove that Proposition \ref{0induction} holds true for $k$ if it holds true for $k-1$ when $2< k \le n$. 
This, combined with the base case $k=2$ established in the previous step, completes the proof of Proposition  \ref{0induction}. 

\subsubsection{Developability} 
Based on the induction hypothesis for $k-1$, we first prove a weaker result in $k$-dimensional slices of $\Omega$ than Proposition 
\ref{0induction}. That is, we prove that $u$ is developable on all $k$-dimensional slices satisfying assumptions 
(1)-(4) of Proposition \ref{0induction} in the following sense: 

\begin{proposition}\label{0developability}
Suppose Proposition \ref{0induction} is true for any $(k-1)$-dimensional slice of $\Omega$ on which assumptions (1)-(4) 
are satisfied. Let $P_k$ be any $k$-dimensional plane such that assumptions (1)-(4) for $u$ holds on $P_k\cap \Omega$, then 
for every $x\in \Omega$, either
$u$ is {\rm affine} in a neighborhood in $P_k\cap \Omega$ of $x$, or there exists a unique
$(k-1)$-dimensional hyperplane $P^x_{k-1}\ni x$ of $P_k$ such that $u$ is {\rm affine} on
the connected component of $x$ in $P^x_{k-1}\cap \Omega$.
\end{proposition}

{\em Proof.} We first need to define a terminology that is the higher dimensional version 
of ``line segments joining the boundary of some domain at both ends''. 

\begin{definition}\label{0kplane}
By a \textit{$k$-plane $P$ in $\Sigma$}, we mean a connected component of a
$k$-dimensional plane $P\cap \Sigma$, where $\Sigma$ is any $N$-dimensional
region with $N\geq k\geq 1$.
\end{definition}
\begin{remark}
We emphasize here that such $k$-plane $P$ in $\Sigma$ refers to \textit{not} the
entire plane, but just the part inside a region. On the other hand, it refers to
the \textit{entire} connected part inside this region. 
\end{remark}
 
Let $\mathbf{v}$ be any unit directional vector of $P_k$,
let $\mathbf{v}_1,\cdots,\mathbf{v}_{k-1}$ be a set of \textit{linearly independent unit} vectors of $P_k$
perpendicular to $\mathbf{v}$. We parametrize the family of
$(k-1)$-dimensional planes parallel to the space spanned by these vectors as
follows:
\begin{equation*}
\ds P_{\mathbf{v}_1\cdots\mathbf{v}_{k-1}}^y=\{z: 
z=y+ \sum_{i=1}^{k-1} s_i\mathbf{v}_i, \, s_i \in
\bbbr\}, \quad y\in \mathrm{span}\langle\mathbf{v}\rangle.
\end{equation*}

\begin{lemma}\label{0l3}
Given the direction $\mathbf{v}$, for a.e. $y\in
\mathrm{span}\langle\mathbf{v}\rangle$, $u$ is
$C_{\rm loc}^{1,1/2}$ and is an isometry on $P_{\mathbf{v}_1\cdots\mathbf{v}_{k-1}}^y\cap \Omega$. Moreover for every 
$x\in P_{\mathbf{v}_1\cdots\mathbf{v}_{k-1}}^y\cap \Omega$, $u$ is
either affine on a $(k-1)$-dimensional region in
$P_{\mathbf{v}_1\cdots\mathbf{v}_{k-1}}^y\cap \Omega$ containing $x$, or affine on a
$(k-2)$-plane in $P_{\mathbf{v}_1\cdots\mathbf{v}_{k-1}}^y\cap \Omega$ passing through $x$. 
\end{lemma}
{\em Proof.}
Since $u$ satisfies assumptions (1)-(4) on $P_k\cap \Omega$, by Fubini Theorem, for a.e. 
$y\in
\mathrm{span}\langle\mathbf{v}\rangle$, assumptions (1)-(4) are also satisfied on $P_{\mathbf{v}_1\cdots\mathbf{v}_{k-1}}^y\cap \Omega$. 
Hence by our induction hypothesis on $(k-1)$ slices of $\Omega$, $\nabla u$ is $C_{\rm loc}^{0,1/2}$ on $P_{\mathbf{v}_1\cdots\mathbf{v}_{k-1}}^y\cap \Omega$. 
By assumption (2) $\nabla u^T\nabla u={\rm I}$ a.e., and hence everywhere 
in $P_{\mathbf{v}_1\cdots\mathbf{v}_{k-1}}^y\cap \Omega$ 
by continuity. Therefore, by assumption (1) and the chain rule in Lemma \ref{0chain}, $u$ is an isometry on $P_{\mathbf{v}_1\cdots\mathbf{v}_{k-1}}^y\cap \Omega$. 

Moreover by our induction hypothesis, for every $x\in P_{\mathbf{v}_1\cdots\mathbf{v}_{k-1}}^y\cap \Omega$, $\nabla u$ is 
either constant on a $(k-1)$-dimensional region in
$P_{\mathbf{v}_1\cdots\mathbf{v}_{k-1}}^y\cap \Omega$ containing $x$, or constant on an
$(k-2)$-plane in $P_{\mathbf{v}_1\cdots\mathbf{v}_{k-1}}^y\cap \Omega$ passing through $x$. Hence by the 
the chain rule in Lemma \ref{0chain}, $u$ is either affine on $(k-1)$ dimensional regions in $P_{\mathbf{v}_1\cdots\mathbf{v}_{k-1}}^y\cap \Omega$, or 
affine on $(k-2)$-planes in $P_{\mathbf{v}_1\cdots\mathbf{v}_{k-1}}^y\cap \Omega$. 
The proof is complete.
\hspace*{\fill} $\Box$

Now we want to show that a substantial part of Lemma \ref{0l3} is true for \textit{every} rather
than a.e. $(k-1)$-dimensional planes in $\Omega$.
\begin{lemma}\label{0l4}
Given direction $\mathbf{v}$, for all $y\in
\mathrm{span}\langle\mathbf{v}\rangle$ and for all $x\in P_{\mathbf{v}_1\cdots\mathbf{v}_{k-1}}^y\cap \Omega$, $u$ is
either an {\rm affine isometry} on a $(k-1)$-dimensional region in
$P_{\mathbf{v}_1\cdots\mathbf{v}_{k-1}}^y\cap \Omega$ containing $x$, or an {\rm affine isometry}
on a
$(k-2)$-plane in $P_{\mathbf{v}_1\cdots\mathbf{v}_{k-1}}^y\cap \Omega$ passing through $x$.
\end{lemma}
\begin{remark}
We obtain from the proof of Lemma \ref{0l3} that $u$ is $C^1$ on a.e.
planes. However, Lemma \ref{0l3} does \textit{not} imply $u$ is $C^1$ on
\textit{every} plane because even though $\nabla u$ is continuous on a.e.
planes, we cannot conclude from here that $\nabla u$ is continuous in $\Omega$,
so we cannot pass to the limit to conclude as Lemma \ref{0l3}.
\end{remark}
{\em Proof.}
Given $y\in
\mathrm{span}\langle\mathbf{v}\rangle$,
Lemma \ref{0l3} guarantees a sequence $y_m\in
\mathrm{span}\langle\mathbf{v}\rangle$,
$y_m\rightarrow y$ such that Lemma \ref{0l3} is true on
$P_{\mathbf{v}_1\cdots\mathbf{v}_{k-1}}^{y_m}\cap \Omega$ for every $m$.

Let $x\in P_{\mathbf{v}_1\cdots\mathbf{v}_{k-1}}^y\cap \Omega$, we divide the
proof into the following two cases:
\begin{enumerate}
	\item There is a sequence of $(k-2)$-planes $P_m$ in
$P_{\mathbf{v}_1\cdots\mathbf{v}_{k-1}}^{y_m} \cap \Omega$ on which $u$ is an
affine isometry and $P_m$ converges to $x$ in distance.
        \item There does not exist such a sequence of $(k-2)$-planes.
\end{enumerate}
Suppose we are in case (1), then the limit of $P_m$ must also be a $(k-2)$-plane
$P$ in $P_{\mathbf{v}_1\cdots\mathbf{v}_{k-1}}^y\cap \Omega$ passing through $x$. Also since $u$ is
Lipschitz continuous, $u$ must also be an
affine isometry on $P$, which proves the Lemma in this case (Figure \ref{fig6}). 
\begin{figure}[ht]
\centering
\begin{tikzpicture}[line cap=round,line join=round,>=triangle
45,x=1.0cm,y=1.0cm]
\clip(-4.28,1.3) rectangle (6.18,6.26);
\clip(-4.28,1.3) rectangle (6.18,6.26);
\fill[line width=0pt,color=qqffff,fill=qqffff,fill opacity=0.5] (-1.5,5.4) --
(-0.207,5.4) -- (-2.85,3) -- (-3.9,3) -- cycle;
\fill[line width=0pt,color=qqffff,fill=qqffff,fill opacity=0.5] (0.806,5.4) --
(-0.897,3) -- (1.292,3) -- cycle;
\fill[line width=0pt,color=qqffff,fill=qqffff,fill opacity=0.5] (4.188,5.4) --
(5.188,4.888) -- (5.7,5.4) -- cycle;
\fill[line width=0pt,color=ccffcc,fill=ccffcc,fill opacity=0.5] (-0.178,4.8) --
(-3.649,2.4) -- (-1.551,2.4) -- (0.211,4.8) -- cycle;
\fill[line width=0pt,color=ccffcc,fill=ccffcc,fill opacity=0.5] (2.248,4.8) --
(2.736,2.4) -- (3.3,2.4) -- (5.7,4.8) -- cycle;
\fill[line width=0pt,color=qqffff,fill=qqffff,fill opacity=0.5] (1.964,5.4) --
(2.666,3) -- (3.3,3) -- (4.223,3.923) -- (2.742,5.4) -- cycle;
\draw [color=qqqqqq] (-1.5,5.4)-- (-3.9,3);
\draw [color=qqqqqq] (-3.9,3)-- (3.3,3);
\draw [color=qqqqqq] (3.3,3)-- (5.7,5.4);
\draw [color=qqqqqq] (5.7,5.4)-- (-1.5,5.4);
\draw [color=qqqqqq] (-3.9,2.4)-- (3.3,2.4);
\draw [color=qqqqqq] (3.3,2.4)-- (5.7,4.8);
\draw [color=qqqqqq] (-3.9,1.8)-- (3.3,1.8);
\draw [color=qqqqqq] (3.3,1.8)-- (5.7,4.2);
\draw [color=qqqqqq] (-0.207,5.4)-- (-2.85,3);
\draw [color=qqqqqq] (0.246,5.4)-- (-1.747,3);
\draw [line width=1.2pt,color=qqqqqq] (0.806,5.4)-- (-0.897,3);
\draw [color=qqqqqq] (0.806,5.4)-- (1.292,3);
\draw [color=qqqqqq] (1.439,5.4)-- (2.015,3);
\draw [color=qqqqqq] (1.964,5.4)-- (2.666,3);
\draw [color=qqqqqq] (2.742,5.4)-- (4.223,3.923);
\draw [color=qqqqqq] (3.447,5.4)-- (4.611,4.311);
\draw [color=qqqqqq] (4.188,5.4)-- (5.188,4.888);
\draw [line width=0.2pt,dash pattern=on 3pt off 3pt,color=yqyqyq] (-1.5,4.2)--
(-3.3,2.4);
\draw [line width=0.2pt,dash pattern=on 3pt off 3pt,color=yqyqyq] (-1.5,4.2)--
(5.1,4.2);
\draw [line width=0.2pt,dash pattern=on 3pt off 3pt,color=qqqqqq] (-1.5,4.8)--
(-3.3,3);
\draw [line width=0.2pt,dash pattern=on 3pt off 3pt,color=qqqqqq] (-1.5,4.8)--
(5.1,4.8);
\draw [line width=0.2pt,dash pattern=on 3pt off 3pt,color=qqqqqq] (-0.827,4.8)--
(-2.599,3.701);
\draw [line width=0.2pt,dash pattern=on 3pt off 3pt,color=qqqqqq] (-0.178,4.8)--
(-3.649,2.4);
\draw [line width=0.2pt,dash pattern=on 3pt off 3pt,color=qqqqqq] (0.211,4.8)--
(-1.551,2.4);
\draw [dash pattern=on 3pt off 3pt,color=qqqqqq] (0.692,4.8)-- (-0.249,2.4);
\draw [line width=0.2pt,dash pattern=on 3pt off 3pt,color=qqqqqq] (1.007,4.8)--
(1.09,2.4);
\draw [line width=0.2pt,dash pattern=on 3pt off 3pt,color=qqqqqq] (1.581,4.8)--
(1.849,2.4);
\draw [line width=0.2pt,dash pattern=on 3pt off 3pt,color=qqqqqq] (2.248,4.8)--
(2.736,2.4);
\draw [color=qqqqqq] (-2.781,3)-- (-3.649,2.4);
\draw [color=qqqqqq] (-1.111,3)-- (-1.551,2.4);
\draw [line width=1.2pt,color=qqqqqq] (-0.013,3)-- (-0.249,2.4);
\draw [color=qqqqqq] (1.069,3)-- (1.09,2.4);
\draw [color=qqqqqq] (1.782,3)-- (1.849,2.4);
\draw [color=qqqqqq] (2.614,3)-- (2.736,2.4);
\draw [color=qqqqqq] (5.1,4.8)-- (5.7,4.8);
\draw [color=qqqqqq] (5.1,4.2)-- (5.7,4.2);
\draw [color=qqqqqq] (-3.3,2.4)-- (-3.9,1.8);
\draw [color=qqqqqq] (-3.3,3)-- (-3.9,2.4);
\draw [dash pattern=on 3pt off 3pt,color=qqqqqq] (0.72,4.2)-- (0.366,1.8);
\draw [line width=1.2pt,color=qqqqqq] (0.455,2.4)-- (0.366,1.8);
\draw [color=qqqqqq](0.7,3.9) node[anchor=north west] {$x $};
\draw [color=qqqqqq](-0.42,4.12) node[anchor=north west] {$P_m $};
\draw [->,>=stealth',color=qqqqqq] (0.23,5.77) -- (0.381,4.8);
\draw [->,>=stealth',line width=0.2pt,dash pattern=on 3pt off 3pt,color=qqqqqq] (0.23,5.77)
-- (0.687,4.787);
\draw [color=qqqqqq](-1.86,6.28) node[anchor=north west] {\scriptsize $u$ is
affine on $P_m$ for each $m$, $P_m\rightarrow P\Rightarrow u$ is affine on
$P$.};
\draw [color=qqqqqq](-3.84,4.86) node[anchor=north west] {\scriptsize$
m\rightarrow \infty $};
\draw [color=qqqqqq](-3.74,5.66) node[anchor=north west] {\scriptsize
$P_{\mathbf{v}_1\cdots\mathbf{v}_{k-1}}^{y_m} \cap \Omega$};
\draw [color=qqqqqq](-3.4,2.44) node[anchor=north west] {\scriptsize $
P_{\mathbf{v}_1\cdots\mathbf{v}_{k-1}}^{y} \cap \Omega$};
\draw [color=qqqqqq] (0.23,5.77)-- (0.405,5.393);
\draw [->,>=stealth',dash pattern=on 3pt off 3pt,color=qqqqqq] (-3.18,5.02) --
(-3.16,3.92);
\draw [->,>=stealth',color=qqqqqq] (-2.16,5.32) -- (-1.38,5.04);
\draw [->,>=stealth',line width=0.2pt,dash pattern=on 3pt off 3pt,color=qqqqqq]
(-2.16,5.32) -- (-1.98,4.08);
\draw [color=qqqqqq] (-2.16,5.32)-- (-2.088,4.822);
\draw [color=qqqqqq](0.38,2.48) node[anchor=north west] {$P$};
\begin{scriptsize}
\fill [color=qqqqqq] (0.659,3.788) circle (1.0pt);
\end{scriptsize}
\end{tikzpicture}
\caption{Case (1).}
\label{fig6}
\end{figure}

Suppose now we are in case (2). If we cannot find such a sequence of
$(k-2)$-planes, then we can find $x_m\in
P_{\mathbf{v}_1\cdots\mathbf{v}_{k-1}}^{y_m} \cap \Omega$, $x_m\rightarrow x$
with the property that there is $\epsilon>0$ such that $u$ is an affine isometry on
$B^{k-1}(x_m, \epsilon)\subset P_{\mathbf{v}_1\cdots\mathbf{v}_{k-1}}^{y_m} \cap
\Omega$. Otherwise, there will again be a sequence of $(k-2)$-planes (i.e. the
boundaries of the maximal affine regions containing $x_m$) converging to $x$ in
distance, contradiction to the fact that we are in case (2). Continuity of $u$
then must force $u$ to be an affine isometry on $B^{k-1}(x, \epsilon)\subset
P_{\mathbf{v}_1\cdots\mathbf{v}_{n-1}}^y\cap \Omega$, which again proves the
lemma in this case (Figure \ref{fig7}). The proof is complete.
\hspace*{\fill} $\Box$
\begin{figure}[ht]
\centering
\begin{tikzpicture}[line cap=round,line join=round,>=triangle
45,x=1.0cm,y=1.0cm]
\clip(-4,1.06) rectangle (6.26,5.94);
\fill[line width=0pt,color=cqcqcq,fill=cqcqcq,fill opacity=0.5] (1.949,3.8) --
(-1.881,1.4) -- (1.269,1.4) -- cycle;
\fill[line width=0pt,color=ccffcc,fill=ccffcc,fill opacity=0.5] (2.918,4.4) --
(4.073,2.613) -- (5.86,4.4) -- cycle;
\fill[line width=0pt,color=ccffcc,fill=ccffcc,fill opacity=0.5] (-1.34,4.4) --
(-3.74,2) -- (-2.384,2) -- (0.44,4.4) -- cycle;
\fill[line width=0pt,color=ccffcc,fill=ccffcc,fill opacity=0.5] (0.965,4.4) --
(-1.701,2) -- (0.889,2) -- (1.409,4.4) -- cycle;
\fill[line width=0pt,color=qqffff,fill=qqffff,fill opacity=0.5] (1.002,5) --
(-1.376,2.6) -- (1.972,2.6) -- cycle;
\fill[line width=0pt,color=qqffff,fill=qqffff,fill opacity=0.5] (1.592,5) --
(2.541,2.6) -- (3.46,2.6) -- (4.071,3.211) -- (2.798,5) -- cycle;
\draw [color=qqqqqq] (-1.34,5)-- (-3.74,2.6);
\draw [color=qqqqqq] (-3.74,2.6)-- (3.46,2.6);
\draw [color=qqqqqq] (3.46,2.6)-- (5.86,5);
\draw [color=qqqqqq] (5.86,5)-- (-1.34,5);
\draw [line width=0.2pt,dash pattern=on 3pt off 3pt,color=qqqqqq] (-1.34,4.4)--
(-3.74,2);
\draw [color=qqqqqq] (-3.74,2)-- (3.46,2);
\draw [color=qqqqqq] (3.46,2)-- (5.86,4.4);
\draw [line width=0.2pt,dash pattern=on 3pt off 3pt,color=qqqqqq] (5.86,4.4)--
(-1.34,4.4);
\draw [line width=0.2pt,dash pattern=on 3pt off 3pt,color=yqyqyq] (-1.34,3.8)--
(-3.74,1.4);
\draw [color=qqqqqq] (-3.74,1.4)-- (3.46,1.4);
\draw [color=qqqqqq] (3.46,1.4)-- (5.86,3.8);
\draw [line width=0.2pt,dash pattern=on 3pt off 3pt,color=yqyqyq] (5.86,3.8)--
(-1.34,3.8);
\draw [color=qqqqqq] (-0.727,5)-- (-3.237,2.6);
\draw [color=qqqqqq] (0.301,5)-- (-2.514,2.6);
\draw [color=qqqqqq] (1.002,5)-- (-1.376,2.6);
\draw [color=qqqqqq] (1.002,5)-- (1.972,2.6);
\draw [color=qqqqqq] (1.592,5)-- (2.541,2.6);
\draw [color=qqqqqq] (2.798,5)-- (4.071,3.211);
\draw [color=qqqqqq] (3.395,5)-- (4.497,3.637);
\draw [color=qqqqqq] (4.209,5)-- (4.914,4.054);
\draw [color=qqqqqq] (4.878,5)-- (5.312,4.452);
\draw [color=qqqqqq] (-3.14,2.6)-- (-3.74,2);
\draw [color=qqqqqq] (-3.14,2)-- (-3.74,1.4);
\draw [color=qqqqqq] (5.26,3.8)-- (5.86,3.8);
\draw [color=qqqqqq] (5.26,4.4)-- (5.86,4.4);
\draw [line width=0.2pt,dash pattern=on 3pt off 3pt,color=qqqqqq] (0.44,4.4)--
(-2.384,2);
\draw [line width=0.2pt,dash pattern=on 3pt off 3pt,color=qqqqqq] (0.965,4.4)--
(-1.701,2);
\draw [line width=0.2pt,dash pattern=on 3pt off 3pt,color=qqqqqq] (1.409,4.4)--
(0.889,2);
\draw [line width=0.2pt,dash pattern=on 3pt off 3pt,color=qqqqqq] (1.779,4.4)--
(1.845,2);
\draw [line width=0.2pt,dash pattern=on 3pt off 3pt,color=qqqqqq] (2.321,4.4)--
(3.697,2.237);
\draw [line width=0.2pt,dash pattern=on 3pt off 3pt,color=qqqqqq] (2.918,4.4)--
(4.073,2.613);
\draw [line width=0.2pt,dash pattern=on 3pt off 3pt,color=qqqqqq] (1.949,3.8)--
(-1.881,1.4);
\draw [line width=0.2pt,dash pattern=on 3pt off 3pt,color=qqqqqq] (1.949,3.8)--
(1.269,1.4);
\draw [color=qqqqqq] (-1.678,2.6)-- (-2.384,2);
\draw [color=qqqqqq] (-1.035,2.6)-- (-1.701,2);
\draw [color=qqqqqq] (1.019,2.6)-- (0.889,2);
\draw [color=qqqqqq] (2.541,2.6)-- (1.828,2.6);
\draw [color=qqqqqq] (1.828,2.6)-- (1.845,2);
\draw [color=qqqqqq] (3.838,2.978)-- (4.073,2.613);
\draw [color=qqqqqq] (-0.924,2)-- (-1.881,1.4);
\draw [color=qqqqqq] (1.439,2)-- (1.269,1.4);
\draw [color=qqqqqq](-2.52,5.82) node[anchor=north west] {\scriptsize  $u$ is
affine on $B(x_m,\epsilon)$ for each $m$, $x\rightarrow x\Rightarrow u$ is
affine on $B(x,\epsilon)$. };
\draw [line width=0.2pt,dash pattern=on 3pt off 3pt,color=qqqqqq] (2.095,4.4)--
(2.662,2);
\draw [color=qqqqqq] (2.52,2.6)-- (2.662,2);
\draw [rotate around={-1.42:(0.685,3.765)},color=qqqqqq] (0.685,3.765) ellipse
(0.668cm and 0.283cm);
\draw [color=qqqqqq] (0.685,3.765)-- (1.173,3.948);
\draw [color=qqqqqq](-3.66,5.4) node[anchor=north west] {\scriptsize $
P_{\mathbf{v}_1\cdots\mathbf{v}_{k-1}}^{y_m} \cap \Omega $};
\draw [color=qqqqqq](-3.4,2.02) node[anchor=north west] {\scriptsize $
P_{\mathbf{v}_1\cdots\mathbf{v}_{k-1}}^{y} \cap \Omega$};
\draw [->,>=stealth',dash pattern=on 3pt off 3pt,color=qqqqqq] (-2.78,4.62) -- (-2.75,3.7);
\draw [rotate around={-1.893:(0.335,3.03)},line width=0.2pt,dash pattern=on 3pt
off 3pt,color=qqqqqq] (0.335,3.03) ellipse (0.67cm and 0.286cm);
\draw [rotate around={-1.909:(0.83,2.44)},line width=0.2pt,dash pattern=on 3pt
off 3pt,color=qqqqqq] (0.83,2.44) ellipse (0.665cm and 0.286cm);
\draw [color=qqqqqq](-3.46,4.64) node[anchor=north west] {\scriptsize $
m\rightarrow \infty $};
\draw [line width=0.2pt,dash pattern=on 3pt off 3pt,color=qqqqqq] (0.335,3.03)--
(-0.003,3.286);
\draw [line width=0.2pt,dash pattern=on 3pt off 3pt,color=qqqqqq] (0.83,2.44)--
(0.268,2.303);
\draw [color=qqqqqq](0.56,3.84) node[anchor=north west] {\scriptsize $ x_m$};
\draw [color=qqqqqq](0.56,2.5) node[anchor=north west] {\scriptsize $ x $};
\draw [color=qqqqqq](0.6,4.16) node[anchor=north west] {\scriptsize $ \epsilon
$};
\draw [color=qqqqqq](0,3.2) node[anchor=north west] {\scriptsize $ \epsilon $};
\draw [color=qqqqqq](0.42,2.74) node[anchor=north west] {\scriptsize $ \epsilon
$};
\draw [->,>=stealth',color=qqqqqq] (0.26,5.24) -- (0.276,3.997);
\draw [->,>=stealth',line width=0.2pt,dash pattern=on 3pt off 3pt,color=qqqqqq] (0.26,5.24)
-- (0.543,3.297);
\draw [color=qqqqqq] (0.26,5.24)-- (0.456,3.891);
\draw [->,>=stealth',color=qqqqqq] (-1.9,5.05) -- (-1.31,4.59);
\draw [->,>=stealth',line width=0.2pt,dash pattern=on 3pt off 3pt,color=qqqqqq] (-1.9,5.05)
-- (-1.5,3.9);
\draw [color=qqqqqq] (-1.9,5.05)-- (-1.749,4.617);
\begin{scriptsize}
\fill [color=qqqqqq] (0.685,3.765) circle (1.0pt);
\fill [color=qqqqqq] (0.335,3.03) circle (1.0pt);
\fill [color=qqqqqq] (0.83,2.44) circle (1.0pt);
\end{scriptsize}
\end{tikzpicture}
\caption{Case (2)}
\label{fig7}
\end{figure}

\begin{lemma}\label{0l6}
Suppose $u$ is an affine isometry on two line segments $C_1$ and $C_2$ in $P_k\cap \Omega$ intersecting at 
a point $x$ in the interior of both $C_1$ and $C_2$.
Let $H$ be the convex hull of the line segments $C_1$ and $C_2$, then $u$ is an
affine isometry on $H\cap \Omega$.
\end{lemma}
{\em Proof.}
We parametrize $C_1$ and $C_2$ by $\{x+t\mathbf{v}_1, t\in [-a,b]\}$ and
$\{x+s\mathbf{v}_2, s\in [-c,d]\}$, respectively, with both $\mathbf{v}_1$ and
$\mathbf{v}_2$ unit vectors. We can assume $\mathbf{v}_1$ and $\mathbf{v}_2$ are linearly
independent, otherwise, the conclusion of the lemma is obvious. Since $u$ is
affine on both $C_1$ and $C_2$, $u(C_1)$ and $u(C_2)$ are both line segments in
$\bbbr^{n+1}$. We can again parametrize the lines that contains the line
segments $u(C_1)$ and $u(C_2)$ by $u(x)+t\tilde{\mathbf{v}}_1$ and
$u(x)+s\tilde{\mathbf{v}}_2$, where both $\tilde{\mathbf{v}}_1$ and
$\tilde{\mathbf{v}}_2$ are unit vectors due to the isometry assumption. 

Let $y\in H\cap \Omega$, we can of course
assume that $y$ is neither in $C_1$ nor $C_2$, otherwise, there is nothing to
prove. In this way, we can find a line $L_3$ passing through $y$ and intersecting
$C_1$ at only one point, denoted $x_{13}$; and $C_2$ at only one point, denoted
$x_{23}$, where the segment $\overline{x_{13}x_{23}}$ lies inside $\Omega$. Since $x_{13}\in C_1$, $x_{13}=x+t_0\mathbf{v}_1$ for some $t_0\in
[-a,b]$. Similarly $x_{23}=x+s_0\mathbf{v}_2$ for some $s_0\in [-c,d]$.
Then since
\begin{equation}\label{080}
 y=wx_{13}+(1-w)x_{23} \quad \textrm{for some } w\in [0,1],
\end{equation}
it follows
$$y=x+wt_0\mathbf{v}_1+(1-w)s_0\mathbf{v}_2.$$
To prove that $u$ is an affine isometry on $H$, we need to prove
\begin{equation}\label{08}
\texttt{}u(y)=u(x)+wt_0\tilde{\mathbf{v}}_1+(1-w)s_0\tilde{\mathbf{v}_2}.
\end{equation}

We first claim that the angle between line segments $u(C_1)$ and $u(C_2)$ is the
same as the angle between $C_1$ and $C_2$. Since $x$ is in the interior of $C_1$
and $C_2$, we can construct a parallelogram $ABCD$ centered at $x$, with $A,C\in
C_1$ and $B,D\in C_2$. Since $u$ is an affine isometry on $C_1$ and $C_2$, $|u(A)-u(x)|=|A-x|, |u(B)-u(x)|=|B-x|, |u(C)-u(x)|=|C-x|$ and
$|u(D)-u(x)|=|D-x|$. On the other hand, $|u(A)-u(B)|\leq |A-B|$ and
$|u(B)-u(C)|\leq |B-C|$ since $u$ is 1-Lipschitz (Figure \ref{fig9}).
\begin{figure}[ht]
\centering
\begin{tikzpicture}[line cap=round,line join=round,>=triangle
45,x=1.0cm,y=1.0cm]
\clip(-4.3,1.02) rectangle (4.72,6.3);
\draw [shift={(-2.178,3.542)}] (0,0) -- (122.327:0.6) arc (122.327:230.136:0.6) -- cycle;
\draw [shift={(-2.178,3.542)}] (0,0) -- (50.136:0.6) arc (50.136:122.327:0.6) -- cycle;
\draw [shift={(2.458,3.542)}] (0,0) -- (57.673:0.6) arc (57.673:129.864:0.6) -- cycle;
\draw [shift={(2.458,3.542)}] (0,0) -- (-50.136:0.6) arc (-50.136:57.673:0.6) -- cycle;
\draw [color=qqqqqq] (-2.86,4.62)-- (-0.86,5.12);
\draw [color=qqqqqq] (-2.86,4.62)-- (-3.501,1.957);
\draw [color=qqqqqq] (-0.86,5.12)-- (-1.524,2.51);
\draw [color=qqqqqq] (-3.501,1.957)-- (-1.524,2.51);
\draw [color=qqqqqq] (-3.519,5.661)-- (-1.014,1.704);
\draw [color=qqqqqq] (-0.382,5.693)-- (-3.947,1.423);
\draw [color=qqqqqq] (-3.579,5.086)-- (-0.334,5.412);
\draw [color=qqqqqq] (3.799,5.661)-- (1.294,1.704);
\draw [color=qqqqqq] (1.14,5.12)-- (3.14,4.62);
\draw [color=qqqqqq] (1.14,5.12)-- (1.804,2.51);
\draw [color=qqqqqq] (0.662,5.693)-- (4.227,1.423);
\draw [color=qqqqqq] (1.804,2.51)-- (3.781,1.957);
\draw [color=qqqqqq] (3.14,4.62)-- (3.781,1.957);
\draw [color=qqqqqq] (0.614,5.412)-- (3.859,5.086);
\draw [shift={(0.038,3.336)},color=qqqqqq] 
plot[domain=0.843:2.171,variable=\t]({1*1.034*cos(\t r)+0*1.034*sin(\t
r)},{0*1.034*cos(\t r)+1*1.034*sin(\t r)});
\draw [color=qqqqqq](-0.86,5.36) node[anchor=north west] {\scriptsize $ A$};
\draw [color=qqqqqq](-3.46,4.96) node[anchor=north west] {\scriptsize $ B$};
\draw [color=qqqqqq](-3.56,2.2) node[anchor=north west] {\scriptsize $ C$};
\draw [color=qqqqqq](-1.4,2.86) node[anchor=north west] {\scriptsize $ D$};
\draw [color=qqqqqq](-2.36,3.52) node[anchor=north west] {\scriptsize $ x $};
\draw [color=qqqqqq](-2.18,5.46) node[anchor=north west] {\scriptsize $  y $};
\draw [color=qqqqqq](-3.32,5.98) node[anchor=north west] {\scriptsize $ C_2$};
\draw [color=qqqqqq](-1.06,6.22) node[anchor=north west] {\scriptsize $ C_1$};
\draw [color=qqqqqq](-3.58,5.26) node[anchor=north west] {\scriptsize $
x_{23}$};
\draw [color=qqqqqq](-0.7,5.54) node[anchor=north west] {\scriptsize $ x_{13}$};
\draw [color=qqqqqq](-2.34,4.08) node[anchor=north west] {\scriptsize $
\alpha_1$};
\draw [color=qqqqqq](-2.7,3.88) node[anchor=north west] {\scriptsize $
\beta_1$};
\draw [color=qqqqqq](0.38,5.3) node[anchor=north west] {\scriptsize $ u(A) $};
\draw [color=qqqqqq](3.16,4.98) node[anchor=north west] {\scriptsize $ u(B)$};
\draw [color=qqqqqq](0.9,3) node[anchor=north west] {\scriptsize $ u(D)$};
\draw [color=qqqqqq](3.76,2.36) node[anchor=north west] {\scriptsize $ u(C)$};
\draw [color=qqqqqq](2.14,3.46) node[anchor=north west] {\scriptsize $ u(x)$};
\draw [color=qqqqqq](2.22,4.16) node[anchor=north west] {\scriptsize $ \alpha_2
$};
\draw [color=qqqqqq](2.58,3.94) node[anchor=north west] {\scriptsize $ \beta_2
$};
\draw [color=qqqqqq](1.42,2.08) node[anchor=north west] {\scriptsize $ u(C_2)
$};
\draw [color=qqqqqq](3.2,1.76) node[anchor=north west] {\scriptsize $ u(C_1)$};
\draw [color=qqqqqq](0.84,5.82) node[anchor=north west] {\scriptsize $
u(x_{13})$};
\draw [color=qqqqqq](2.76,5.66) node[anchor=north west] {\scriptsize $
u(x_{23})$};
\draw [color=qqqqqq](2.06,5.5) node[anchor=north west] {\scriptsize $ u(y)$};
\draw [color=qqqqqq] (0.726,4.107)-- (0.595,4.336);
\draw [color=qqqqqq] (0.726,4.107)-- (0.514,4.091);
\draw [color=qqqqqq](-0.06,4.42) node[anchor=north west] {\footnotesize $ u$};
\draw [color=qqqqqq](-1.08,3.98) node[anchor=north west] {\scriptsize $
\alpha_2\leq \alpha_1, \beta_2\leq \beta_1$};
\draw [color=qqqqqq](-1.16,3.46) node[anchor=north west] {\scriptsize $
\Rightarrow \alpha_2=\alpha_1, \beta_2=\beta_1$};
\begin{scriptsize}
\fill [color=qqqqqq] (-3.501,1.957) circle (1.0pt);
\fill [color=qqqqqq] (-2.86,4.62) circle (1.0pt);
\fill [color=qqqqqq] (-0.86,5.12) circle (1.0pt);
\fill [color=qqqqqq] (-1.524,2.51) circle (1.0pt);
\fill [color=qqqqqq] (-3.519,5.661) circle (1.0pt);
\fill [color=qqqqqq] (-1.014,1.704) circle (1.0pt);
\fill [color=qqqqqq] (-0.382,5.693) circle (1.0pt);
\fill [color=qqqqqq] (-3.947,1.423) circle (1.0pt);
\fill [color=qqqqqq] (-2.177,5.227) circle (1.0pt);
\fill [color=qqqqqq] (-2.178,3.542) circle (1.0pt);
\fill [color=qqqqqq] (1.14,5.12) circle (1.0pt);
\fill [color=qqqqqq] (3.14,4.62) circle (1.0pt);
\fill [color=qqqqqq] (3.781,1.957) circle (1.0pt);
\fill [color=qqqqqq] (1.804,2.51) circle (1.0pt);
\fill [color=qqqqqq] (3.799,5.661) circle (1.0pt);
\fill [color=qqqqqq] (1.294,1.704) circle (1.0pt);
\fill [color=qqqqqq] (0.662,5.693) circle (1.0pt);
\fill [color=qqqqqq] (4.227,1.423) circle (1.0pt);
\fill [color=qqqqqq] (2.457,5.227) circle (1.0pt);
\fill [color=qqqqqq] (2.458,3.542) circle (1.0pt);
\fill [color=qqqqqq] (-3.18,5.126) circle (1.0pt);
\fill [color=qqqqqq] (-0.642,5.381) circle (1.0pt);
\fill [color=qqqqqq] (0.922,5.381) circle (1.0pt);
\fill [color=qqqqqq] (3.46,5.126) circle (1.0pt);
\end{scriptsize}
\end{tikzpicture}
\caption{}
\label{fig9}
\end{figure}

This implies the angle $\alpha_2$ between
the line segments $\overline{u(x)u(A)}$ and $\overline{u(x)u(B)}$ must be
smaller than or equal to the angle $\alpha_1$ between $\overline{xA}$ and
$\overline{xB}$,
and the angle $\beta_2$ between the line segments $\overline{u(x)u(B)}$ and
$\overline{u(x)u(C)}$ must be smaller than or equal to the angle $\beta_1$
between
$\overline{xB}$ and $\overline{xC}$. Hence $\alpha_2=\alpha_1$ and
$\beta_2=\beta_1$. This proves our claim.

Since by assumption, $u$ is an affine isometry on $\overline {x_{13}x}$ and $\overline{x_{23}x}$, we have
\begin{equation*}
u(x_{13})-u(x)=t_0\tilde{\mathbf{v}}_1\quad\textrm{ and}\quad
u(x_{23})-u(x)=s_0\tilde{\mathbf{v}}_2
\end{equation*}
for the same $t_0, s_0$ and unit vector $\tilde{\mathbf{v}}_1$, $\tilde{\mathbf{v}}_2$ as defined before. 
In particular, $|u(x_{13})-u(x)|=|x_{13}-x|$ and $|u(x_{23})-u(x)|=|x_{23}-x|$. Moreover,
since the angle between line segments $u(C_1)$ and $u(C_2)$ is the same as the
angle between $C_1$ and $C_2$, we have $|x_{13}-x_{23}|=|u(x_{13})-u(x_{23})|$.

On the other hand, $u(\overline{x_{13}x_{23}})$ is a 1-Lipschitz curve, hence
the length the the curve $u(\overline{x_{13}x_{23}})$, denoted by
$|u(\overline{x_{13}x_{23}})|$, satisfies
$|u(\overline{x_{13}x_{23}})|\leq |x_{13}-x_{23}|$. Altogether we have
$$|u(x_{13})-u(x_{23})|\leq |u(\overline{x_{13}x_{23}})|\leq |x_{13}-x_{23}|
=|u(x_{13})-u(x_{23})|. $$
This implies
$$|u(\overline{x_{13}x_{23}})|= |u(x_{13})-u(x_{23})|.$$
Hence the curve $u(\overline{x_{13}x_{23}})$ must coincide with line segment
$\overline{u(x_{13})u(x_{23})}$. Therefore, $u$ also maps the line segment
$\overline{x_{13}x_{23}}$ onto a line segment $\overline{u(x_{13})u(x_{23})}$,
which means $u$ is affine on $\overline{x_{13}x_{23}}$.

Finally, since $u$ is 1-Lipschitz, $|u(x_{13})-u(y)|\leq |x_{13}-y|$ and $|u(x_{23})-u(y)|\leq |x_{23}-y|$. 
However, since $u$ is affine on $\overline{x_{13}x_{23}}$, 
\begin{equation*}
|u(x_{13})-u(x_{23})|=|u(x_{13})-u(y)|+|u(y)-u(x_{23})| 
\leq |x_{13}-y|+|y-x_{23}|=|x_{13}-x_{23}|.
\end{equation*}
But we already showed that $|x_{13}-x_{23}|=|u(x_{13})-u(x_{23})|$. Hence $|u(x_{13})-u(y)|=|x_{13}-y|$ and $|u(x_{23})-u(y)|=|x_{23}-y|$. Therefore, 
$$u(y)=wu(x_{13})+(1-w)u(x_{23}) $$
for the same $w$ as (\ref{080}), which yields (\ref{08}). The proof is complete.
\hspace*{\fill} $\Box$

\begin{corollary}\label{0c3}
Given a $\ell$-dimensional ($\ell\leq k$) region $U$ in $P_k\cap \Omega$, and a line segment $C$ in 
$P_k\cap \Omega$ for which there exists $x\in C\cap U$ that lies in the interior of both $U$ and $C$,
if $u$ is an affine isometry on both $U$ and $C$, then $u$ is an affine isometry on the convex hull $H$
of $U$ and $C$ inside $\Omega$
\end{corollary}
{\em Proof.}
Let $y\in H\cap \Omega$. We need to show that $u(y)=u(x)+t\tilde{\mathbf{v}}$ for some
$\tilde{\mathbf{v}}$ given by a linear combination of directional vectors in
$u(U)$ and $u(C)$ and $|t\tilde{\mathbf{v}}|=|y-x|$. Let
$P_y$ be a $2$-dimensional plane that contains $y$ and $C$. Then $P_y$
intersects $U$ at some line segment $C_y$. Since $u$ is an affine isometry on both $C$ and $C_y$, by
Lemma \ref{0l6}, $u$ is an affine isometry on the convex hull of $C$ and $C_y$ (Figure \ref{fig8}).
\begin{figure}[ht]
 \centering
 \begin{tikzpicture}[line cap=round,line join=round,>=triangle 45,x=1.0cm,y=1.0cm]
\clip(-0.18,-0.9) rectangle (7.82,3.56);
\fill[line width=0pt,color=qqffff,fill=qqffff,fill opacity=0.5] (6.17,2.89) -- (0.82,0.82) -- (6.88,1.52) -- cycle;
\fill[line width=0pt,color=qqffff,fill=qqffff,fill opacity=0.5] (0.72,0.03) -- (0.64,-0.6) -- (2.74,0.12) -- cycle;
\fill[line width=0pt,color=qqffff,fill=qqffff,fill opacity=0.2] (0.82,0.82) -- (0.72,0.03) -- (2.74,0.12) -- (6.88,1.52) -- cycle;
\draw [line width=0.2pt,dash pattern=on 3pt off 3pt] (2,2)-- (7.6,2.24);
\draw [line width=0.2pt,dash pattern=on 3pt off 3pt] (2,2)-- (0,0);
\draw (0,0)-- (5.6,0.24);
\draw (7.6,2.24)-- (5.6,0.24);
\draw [line width=0.2pt,dash pattern=on 3pt off 3pt] (1.5,2.42)-- (0.2,-0.64);
\draw (1.5,2.42)-- (7.52,3.02);
\draw (7.52,3.02)-- (6.22,-0.04);
\draw (0.2,-0.64)-- (6.22,-0.04);
\draw (0.82,0.82)-- (6.88,1.52);
\draw [line width=0.2pt,dash pattern=on 3pt off 3pt] (6.17,2.89)-- (0.64,-0.6);
\draw (0.82,0.82)-- (0,0);
\draw (1.5,2.42)-- (0.82,0.82);
\draw (0.48,0.02)-- (0.2,-0.64);
\draw (1.7,0.07)-- (0.64,-0.6);
\draw (3.35,1.11)-- (6.17,2.89);
\draw (3.26,1.16) node[anchor=north west] {$x$};
\draw (5.88,2.26) node[anchor=north west] {$y$};
\draw (3.68,1.80) node[anchor=north west] {$C$};
\draw (4.1,1.42) node[anchor=north west] {$C_y$};
\draw (6.17,2.89)-- (6.88,1.52);
\draw [line width=0.2pt,dash pattern=on 3pt off 3pt] (0.82,0.82)-- (0.64,-0.6);
\draw (0.72,0.03)-- (0.64,-0.6);
\draw (0.82,0.82)-- (6.17,2.89);
\draw [line width=0.2pt,dash pattern=on 3pt off 3pt] (0.64,-0.6)-- (6.88,1.52);
\draw (0.64,-0.6)-- (2.74,0.12);
\draw (5.2,0.8) node[anchor=north west] {$U$};
\draw (6.72,3.00) node[anchor=north west] {$P_y$};
\draw (7.18,2.22)-- (7.6,2.24);
\draw (0.68,3.28) node[anchor=north west] {\scriptsize $u$ is affine on the convex hull of $C$ and $C_y$};
\draw [->,>=stealth',] (3.08,2.84) -- (3.36,1.5);
\begin{scriptsize}
\fill [color=black] (3.35,1.11) circle (1.0pt);
\fill [color=black] (5.88,2.02) circle (1.0pt);
\end{scriptsize}
\end{tikzpicture}
 \caption{}
 \label{fig8}
\end{figure}
Since this
convex hull
contains both $y$ and $x$, this implies $u(y)=u(x)+t\tilde{\mathbf{v}}$ for some vector $\tilde{\mathbf{v}}$, 
$|t\tilde{\mathbf{v}}|=|y-x|$, and $\tilde{\mathbf{v}}$ 
is a linear combination of directional vectors of $u(C)$
and $u(C_y)$. Our claim then follows because $C_y\subset U$ and $u$ is an affine isometry on $U$, so any 
vectors of $u(C_y)$ is a linear combination of 
vectors in $u(U)$. The proof is complete.
\hspace*{\fill} $\Box$

By obvious induction we then have
\begin{corollary}\label{0c4}
Suppose $U_1$ and $U_2$ are $k_1$ and $k_2$-dimensional regions ($k_1,k_2\leq k$) in
$P_k\cap \Omega$ with nonempty intersections. Moreover, there exists a point $x\in U_1\cap
U_2$ belonging to the interior of both $U_1$ and $U_2$. If $u$ is an affine isometry on
both $U_1$ and $U_2$, then $u$ is an affine isometry on the convex hull of $U_1$ and $U_2$ inside $\Omega$.
\end{corollary}

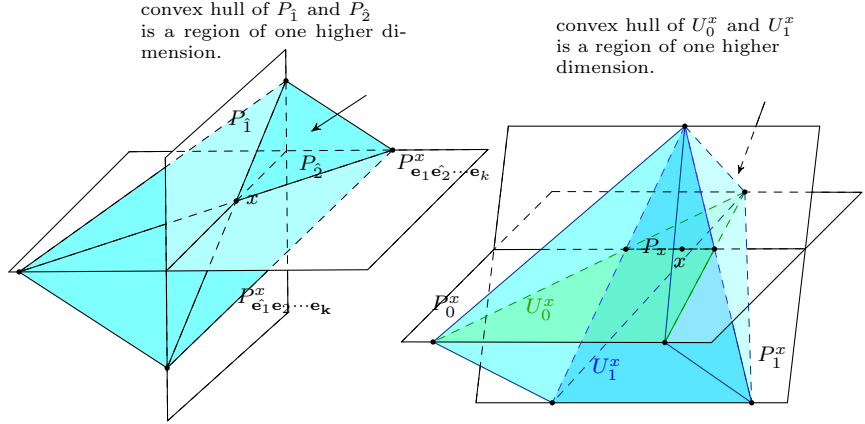
\begin{figure}[ht]
\centering
\begin{tikzpicture}[line cap=round,line join=round,>=triangle
45,x=1.0cm,y=1.0cm]
\clip(-4.18,0.2) rectangle (7.58,6.22);
\fill[line width=0pt,color=qqffff,fill=qqffff,fill opacity=0.5] (-0.306,5.033)
-- (-0.968,3.433) -- (1.112,4.112) -- cycle;
\fill[line width=0pt,color=qqffff,fill=qqffff,fill opacity=0.5] (-1.909,3.885)
-- (-3.854,2.491) -- (-1.886,1.213) -- (-0.543,2.512) -- (-1.897,2.503) --
cycle;
\fill[line width=0pt,color=qqffff,fill=qqffff,fill opacity=0.3] (-0.306,5.033)
-- (-1.909,3.885) -- (-1.897,2.503) -- (-0.543,2.512) -- (1.112,4.112) --
(-0.968,3.433) -- cycle;
\fill[line width=0pt,color=ccffcc,fill=ccffcc,fill opacity=0.4] (5.794,3.554) --
(4.213,2.793) -- (5.391,2.796) -- cycle;
\fill[line width=0pt,color=qqzzff,fill=qqzzff,fill opacity=0.4] (4.997,4.43) --
(4.213,2.793) -- (5.391,2.796) -- cycle;
\fill[line width=0pt,color=qqzzff,fill=qqzzff,fill opacity=0.4] (3.62,1.554) --
(3.235,0.75) -- (5.885,0.75) -- (5.626,1.826) -- (5.35,1.55) -- cycle;
\fill[line width=0pt,color=qqzzff,fill=qqzzff,fill opacity=0.3] (4.213,2.793) --
(3.62,1.554) -- (5.35,1.55) -- (5.626,1.826) -- (5.391,2.796) -- cycle;
\fill[line width=0pt,color=zzffqq,fill=zzffqq,fill opacity=0.3] (4.213,2.793) --
(1.649,1.559) -- (4.73,1.552) -- (5.391,2.796) -- cycle;
\fill[line width=0pt,color=qqffff,fill=qqffff,fill opacity=0.4] (4.997,4.43) --
(1.649,1.559) -- (3.235,0.75) -- (5.885,0.75) -- (5.626,1.826) -- (5.35,1.55) --
(4.73,1.552) -- (5.391,2.796) -- cycle;
\fill[line width=0pt,color=qqffff,fill=qqffff,fill opacity=0.2] (4.997,4.43) --
(5.794,3.554) -- (5.885,0.75) -- (5.626,1.826) -- (5.35,1.55) -- (4.73,1.552) --
(5.391,2.796) -- cycle;
\draw [color=qqqqqq] (-2.39,4.09)-- (-3.99,2.49);
\draw [color=qqqqqq] (-3.99,2.49)-- (0.78,2.52);
\draw [line width=0.2pt,dash pattern=on 3pt off 3pt,color=qqttzz] (5.794,3.554)-- (3.235,0.75);
\draw [line width=0.2pt,dash pattern=on 3pt off 3pt,color=qqqqqq] (-2.39,4.09)-- (2.38,4.12);
\draw [line width=0.2pt,dash pattern=on 3pt off 3pt,color=qqqqqq] (2.38,4.12)-- (0.78,2.52);
\draw [line width=0.2pt,dash pattern=on 3pt off 3pt,color=qqqqqq] (-1.91,4.01)-- (-1.88,0.51);
\draw [color=qqqqqq] (-1.91,4.01)-- (-0.31,5.45);
\draw [line width=0.2pt,dash pattern=on 3pt off 3pt,color=qqqqqq] (-0.31,5.45)-- (-0.28,1.95);
\draw [line width=0.2pt,dash pattern=on 3pt off 3pt,color=qqqqqq] (-1.88,0.51)-- (-0.28,1.95);
\draw [line width=0.2pt,dash pattern=on 3pt off 3pt,color=qqqqqq] (-0.298,4.103)-- (-1.897,2.503);
\draw [line width=0.2pt,dash pattern=on 3pt off 3pt,color=qqqqqq] (-0.306,5.033)-- (-1.886,1.213);
\draw [color=qqqqqq] (-0.306,5.033)-- (1.112,4.112);
\draw [color=qqqqqq] (-0.31,5.45)-- (-0.306,5.033);
\draw [color=qqqqqq] (1.112,4.112)-- (2.38,4.12);
\draw [color=qqqqqq] (2.38,4.12)-- (0.78,2.52);
\draw [color=qqqqqq] (-1.91,4.01)-- (-1.88,0.51);
\draw [color=qqqqqq] (-0.285,2.513)-- (-0.28,1.95);
\draw [color=qqqqqq] (-1.88,0.51)-- (-0.28,1.95);
\draw [color=qqqqqq] (-0.306,5.033)-- (-0.968,3.433);
\draw [color=qqqqqq] (-0.968,3.433)-- (-1.897,2.503);
\draw [line width=0.2pt,dash pattern=on 3pt off 3pt,color=qqqqqq] (1.112,4.112)-- (-1.886,1.213);
\draw [color=qqqqqq] (-0.543,2.512)-- (-1.886,1.213);
\draw [color=qqqqqq] (-1.351,2.507)-- (-1.886,1.213);
\draw [line width=0.2pt,dash pattern=on 3pt off 3pt,color=qqqqqq] (3.23,3.56)-- (1.23,1.56);
\draw [color=qqqqqq] (1.23,1.56)-- (5.35,1.55);
\draw [line width=0.2pt,dash pattern=on 3pt off 3pt,color=qqqqqq] (3.23,3.56)-- (7.35,3.55);
\draw [color=qqqqqq] (7.35,3.55)-- (5.35,1.55);
\draw [color=qqqqqq] (-2.39,4.09)-- (-1.817,4.094);
\draw [line width=0.2pt,dash pattern=on 3pt off 3pt,color=qqqqqq] (1.112,4.112)-- (-3.854,2.491);
\draw [color=qqqqqq] (1.112,4.112)-- (-0.968,3.433);
\draw [color=qqqqqq] (-3.854,2.491)-- (-1.886,1.213);
\draw [line width=0.2pt,dash pattern=on 3pt off 3pt,color=qqqqqq] (-0.306,5.033)-- (-3.854,2.491);
\draw [color=qqqqqq] (-1.909,3.885)-- (-3.854,2.491);
\draw [color=qqqqqq] (2.65,4.43)-- (6.79,4.43);
\draw [color=qqqqqq] (6.79,4.43)-- (6.36,0.75);
\draw [line width=0.2pt,dash pattern=on 3pt off 3pt,color=qqqqqq] (2.65,4.43)-- (2.22,0.75);
\draw [color=qqqqqq] (2.22,0.75)-- (6.36,0.75);
\draw [line width=0.2pt,dash pattern=on 3pt off 3pt,color=qqqqqq] (2.458,2.788)-- (6.599,2.799);
\draw [line width=0.2pt,dash pattern=on 3pt off 3pt,color=qqzzqq] (5.794,3.554)-- (1.649,1.559);
\draw [dash pattern=on 3pt off 3pt,color=qqzzqq] (5.794,3.554)-- (4.73,1.552);
\draw [line width=0.2pt,dash pattern=on 3pt off 3pt,color=ttttff] (4.997,4.43)-- (3.235,0.75);
\draw [line width=0.2pt,dash pattern=on 3pt off 3pt,color=ttttff] (4.997,4.43)-- (5.885,0.75);
\draw [line width=0.2pt,dash pattern=on 3pt off 3pt,color=qqttzz] (4.997,4.43)-- (5.794,3.554);
\draw [line width=0.2pt,dash pattern=on 3pt off 3pt,color=qqttzz] (5.794,3.554)-- (5.885,0.75);
\draw [color=qqttzz] (4.997,4.43)-- (4.73,1.552);
\draw [color=qqttzz] (4.73,1.552)-- (5.885,0.75);
\draw [color=qqttzz] (1.649,1.559)-- (3.235,0.75);
\draw [color=qqttzz] (4.997,4.43)-- (1.649,1.559);
\draw [color=qqqqqq] (6.687,3.552)-- (7.35,3.55);
\draw [color=qqqqqq] (2.458,2.788)-- (1.23,1.56);
\draw [color=qqqqqq] (2.65,4.43)-- (2.458,2.788);
\draw [color=ttttff] (4.997,4.43)-- (5.391,2.796);
\draw [color=qqzzqq] (5.391,2.796)-- (4.73,1.552);
\draw [color=qqqqqq] (2.277,1.239)-- (2.22,0.75);
\draw [color=qqqqqq] (2.458,2.788)-- (3.084,2.79);
\draw [color=qqqqqq] (5.819,2.797)-- (6.599,2.799);
\draw [color=qqqqqq](-0.96,3.64) node[anchor=north west] {\footnotesize $ x$};
\draw [color=qqqqqq](1.04,4.18) node[anchor=north west] {\footnotesize $
P^x_{\mathbf{e}_1 \hat{\mathbf{e}_2}\cdots \mathbf{e}_k} $};
\draw [color=qqqqqq](-1.08,2.4) node[anchor=north west] {\footnotesize $ P^x_{
\hat{\mathbf{e}_1} \mathbf{e}_2\cdots \mathbf{e_k}}$};
\draw [color=qqqqqq](-0.26,4.16) node[anchor=north west] {\footnotesize $
P_{\hat{2}}$};
\draw [color=qqqqqq](-1.2,4.76) node[anchor=north west] {\footnotesize $
P_{\hat{1}} $};
\draw [color=ttzzqq](2.76,2.26) node[anchor=north west] {\footnotesize $
U^x_0$};
\draw [color=qqqqqq](1.5,2.3) node[anchor=north west] {\footnotesize $ P^x_0 $};
\draw [color=qqqqff](3.64,1.46) node[anchor=north west] {\footnotesize $
U^x_1$};
\draw [color=qqqqqq](5.84,1.62) node[anchor=north west] {\footnotesize $ P^x_1
$};
\draw [color=qqqqqq](4.3,3.06) node[anchor=north west] {\footnotesize $ P_x $};
\draw [color=qqqqqq](4.72,2.8) node[anchor=north west] {\footnotesize $ x$};
\draw [color=qqqqqq](-2.36,6.2) node[anchor=north west] {\parbox{3.533
cm}{\scriptsize  convex hull of  $P_{\hat{1}}$ and $P_{\hat{2}}$ \\ is a region
of one higher dimension.}};
\draw [color=qqqqqq](3.16,6) node[anchor=north west] {\parbox{3.272
cm}{\scriptsize  convex hull of  $U^x_0$ and  $U^x_1 $ \\ is a region of one
higher \\ dimension.}};
\draw [->,>=stealth',color=qqqqqq] (0.76,4.84) -- (0.03,4.36);
\draw [->,>=stealth',line width=0.2pt,dash pattern=on 3pt off 3pt,color=qqqqqq] (6.06,4.75) -- (5.71,3.8);
\draw [color=qqqqqq] (6.06,4.75)-- (5.938,4.42);
\draw [color=qqqqqq] (-1.902,3.128)-- (-3.854,2.491);
\draw [color=qqttzz] (5.626,1.826)-- (5.885,0.75);
\draw [dash pattern=on 3pt off 3pt,color=qqqqqq] (4.213,2.793)-- (5.391,2.796);
\begin{scriptsize}
\fill [color=qqqqqq] (-0.306,5.033) circle (1.0pt);
\fill [color=qqqqqq] (-1.886,1.213) circle (1.0pt);
\fill [color=qqqqqq] (-0.968,3.433) circle (1.0pt);
\fill [color=qqqqqq] (1.112,4.112) circle (1.0pt);
\fill [color=qqqqqq] (-3.854,2.491) circle (1.0pt);
\fill [color=qqqqqq] (5.794,3.554) circle (1.0pt);
\fill [color=qqqqqq] (1.649,1.559) circle (1.0pt);
\fill [color=qqqqqq] (4.73,1.552) circle (1.0pt);
\fill [color=qqqqqq] (4.213,2.793) circle (1.0pt);
\fill [color=qqqqqq] (5.391,2.796) circle (1.0pt);
\fill [color=qqqqqq] (4.997,4.43) circle (1.0pt);
\fill [color=qqqqqq] (3.235,0.75) circle (1.0pt);
\fill [color=qqqqqq] (5.885,0.75) circle (1.0pt);
\fill [color=qqqqqq] (4.96,2.795)  circle (1.0pt);
\end{scriptsize}
\end{tikzpicture}
\caption{Case 1 (left) and Case 2 (right).}
\label{fig10}
\end{figure}

Now we are ready to prove Proposition \ref{0developability}. Given $x\in P_k\cap
\Omega$, we first claim that there is a $(k-1)$-dimensional hyperplane $P_0^x$ in $P_k$
and a $(k-1)$-dimensional
neighborhood $U_0^x\subset P_0^x\cap \Omega$ containing $x$ on which $u$ is an affine isometry.
Otherwise, for all $(k-1)$-dimensional hyperplanes in $P_k\cap
\Omega$ that pass through $x$, $x$ is
not contained in any $(k-1)$-dimensional neighborhood on which $u$ is an affine isometry. In particular, let $\mathbf{v}_1, \cdots,\mathbf{v}_k$ 
be linearly independent vectors of $P_k$ and let
$P^x_{\mathbf{v}_1\cdots\hat{\mathbf{v}_i}\cdots\mathbf{v}_k}$, $i=1,...,k$
be the $(k-1)$-dimensional hyperplanes in $\Omega$ passing through $x$ and parallel to the space spanned by 
$\mathbf{v}_1,\cdots,\mathbf{v}_{i-1}, \mathbf{v}_{i+1},\cdots,\mathbf{v}_k$. Since $x$ is
not contained in any $(k-1)$-dimensional neighborhood in $P^x_{\mathbf{v}_1\cdots\hat{\mathbf{v}_i}\cdots\mathbf{v}_k}\cap \Omega$ on which $u$ is an affine isometry, by Lemma
\ref{0l4} there exists $(k-2)$-planes
$P_{\hat{i}}\ni x$ in
$P^x_{\mathbf{v}_1\cdots\hat{\mathbf{v}_i}\cdots\mathbf{v}_k}\cap\Omega$ such that
$u$ is an affine isometry on $P_{\hat{i}}$. By Corollary \ref{0c4}, $u$ is an
affine isometry on the convex hull of $P_{\hat{i}}$ for all $1\leq i\leq k$ (Figure
\ref{fig10} Case 1). Let $\mathbf{v}_{\hat{i}}$ be a directional vector of
$P_{\hat{i}}$. Since $P_{\hat{i}}\subset P^x_{\mathbf{v}_1\cdots
\hat{\mathbf{v}_i} \cdots \mathbf{e}_k}$, which is orthogonal to $\mathbf{v}_i$, at least $k-1$ out of
these $k$ vectors are linearly independent. This convex hull has $k-1$ linearly
independent directional
vectors, hence it must contain a $(k-1)$-dimensional neighborhood of $x$, contradiction to
our
assumption, which proves our claim.

Therefore, we have proved that $x$ must be contained in a $(k-1)$-dimensional neighborhood
$U_0^x\subset P_0^x \cap \Omega$ for some $(k-1)$-dimensional hyperplane
$P_0^x$ and $u$ is an affine isometry on $U_0^x$. If $U_0^x$ is the entire connected component containing $x$ in $P_0^x\cap\Omega$, 
then the conclusion of the proposition is achieved.
Otherwise, we can find a maximal $(k-2)$-plane $P_x$ in
$U_0^x$, which is \textit{not} a $(k-2)$-plane in $P_0^x \cap \Omega$, i.e., it is away from
$\partial \Omega$, on which $u$ is an affine isometry. Let $P^x_1$ be any other  
$(k-1)$-dimensional hyperplane containing the region $P_x$. We have $P_x=U_0^x\cap P_1^x$ and since
the maximal affine region $P_x\subset P^x_1\cap\Omega$
is not a $(k-2)$-plane in $P^x_1\cap\Omega$, by Lemma \ref{0l4}, $x$ must be contained in
a $(k-1)$-dimensional neighborhood $U_1^x\subset P_1^x\cap\Omega$ on which $u$
is an affine isometry (Figure \ref{fig10} Case 2). By Corollary \ref{0c4}, $u$ is an affine isometry on
the convex
hull of $U_0^x$
and $U_1^x$, whose interior is a $k$-dimensional region, which also
achieves the
conclusion of Proposition \ref{0developability}. The proof is complete.
\hspace*{\fill} $\Box$

\subsubsection{Regularity and the conclusion of the inductive step}\label{sregu}

In out last step, we will essentially show that Proposition \ref{0developability} combined with assumptions (1)-(4) 
of Proposition \ref{0induction} for a $k$-dimensional slice $P_k$, implies the conclusion 
of the latter proposition. This will hence conclude the inductive step.  The key point is to show that if $u$ is affine on a $(k-1)$-plane, then its full gradient must be 
constant on the same region. The arguments are very similar to what we used in the proofs of Lemmas \ref{0rigid}-\ref{0continuity}. We will first  prove the following lemma.

\begin{lemma}\label{0rigid2}
Suppose on a $k$-plane $P$ ($1\leq k \leq n$) in $\Omega$ we have the following:
\begin{enumerate}
        \item There is a sequence of smooth functions $u^\epsilon \in C^\infty
(\Omega, \bbbr^{n+1})$ such that
\begin{equation*}
\int_{P} |u^\epsilon-u|^2+|\nabla u^\epsilon-\nabla u|^2+ |\nabla^2
u^\epsilon-\nabla^2 u|^2 d\mathcal{H}^k\rightarrow 0.
\end{equation*}
        \item $\mathrm{rank}\, \nabla^2 u^\ell \leq 1$ and $\nabla^2 u^\ell$ is
symmetric a.e. on $P$ for all $1\leq \ell\leq n+1$.
\end{enumerate}
Then if $u$ is affine on $P$, $\nabla u$ is constant on $P$.
\end{lemma}
{\em Proof.}
Let $\mathbf{v}$ be any unit directional vector in $P$. By assumption (1) and
the chain rule in Lemma \ref{0chain}, $u$ is affine on $P$ implies
$$\nabla u(x) \mathbf{v} =\textrm{constant} \quad \textrm{for a.e. } x\in P. $$
Take the directional derivative one more time, together with assumption (1) we
obtain,
\begin{equation}\label{0aff}
(\mathbf{v})^T \nabla^2 u^\ell \mathbf{v}=0 \quad \textrm{for a.e. } x\in P
\end{equation}
for all $ 1\leq \ell \leq n+1 $.
However, to show that $\nabla u$ is constant on $P$, we need a conclusion stronger
than (\ref{0aff}), i.e.,
\begin{equation}\label{0const}
\nabla^2 u^\ell \mathbf{v}=0 \quad \textrm{for a.e. } x\in P
\end{equation}
for all $1\leq \ell \leq n+1 $. Indeed, by assumption (2), we can write $\nabla ^2 u^\ell$ as
$$\nabla ^2 u^\ell (x)=\lambda(x)\mathbf{b}(x) \otimes \mathbf{b}(x) \quad
\textrm{a.e.}$$
for some scalar function $\lambda$ and $\mathbf{b}\in \mathbb{S}^{n-1}$.
Then (\ref{0aff}) implies,
$$(\mathbf{v})^T\lambda(x)\mathbf{b}(x) \otimes \mathbf{b}(x)\mathbf{v}=
\lambda(x) \langle \mathbf{v}, \mathbf{b}(x)\rangle ^2 =0 \quad \textrm{a.e.}$$
This then implies
$$\lambda(x)\langle \mathbf{v}, \mathbf{b}(x)\rangle=0 \quad \textrm{a.e.} $$
Therefore,
$$ \nabla^2 u^\ell \mathbf{v}=\lambda(x)\langle \mathbf{v}, \mathbf{b}(x)\rangle
\mathbf{b}(x) =0 \quad \textrm{a.e.}$$
which is exactly (\ref{0const}). The proof of the lemma is complete.
\hspace*{\fill} $\Box$
 
Let $P_k$ be any $k$-dimensional plane such that assumptions (1)-(4) in Proposition \ref{0induction} hold on $P_k\cap \Omega$. 

By means of Lemma \ref{0continuity} and making using of Proposition \ref{0developability} and Corollary \ref{0c4}, similar as before we get
$$\partial
U\cap \Omega= \bigcup_{x\in \partial U\cap \Omega}P^U_x$$
where $P^U_x$ is some $(n-1)$-plane
in $\Omega$ containing $x$ with the property that for $x,z\in \partial U\cap \Omega$, $P^U_x=P^U_z$ if $z\in P^U_x$
and $P^U_x\cap P^U_z\cap \Omega=\emptyset$ if $z\notin P^U_x$. 

Similarly as in the proof of Lemma \ref{0continuity} (Figure \ref{fig20}), it suffices to show that the conclusions hold 
locally true. If $x_0\in P_k\cap \Omega$ is a point lying in an affine neighborhood for $u$ in $P_k$, then Lemma \ref{0rigid2} and the assumptions of Proposition 
\ref{0induction} immediately 
imply that $\nabla u$ must be constant in the same neighborhood, which is the desired conclusion. Otherwise, we may and do choose a small $\delta>0$ so that for any region $U$ on which $u$ is affine, the
$k$-dimensional ball $B^k(x_0, \delta)\subset P_k\cap \Omega$ intersects $\partial U$ at
no more than \textit{two} $(k-1)$-planes belonging to $\partial U$. 

We now focus on $B^k(x_0, \delta)\subset P_k\cap \Omega$. For any $x\in B^k(x_0, \delta)$, as in Lemma \ref{0continuity}, we
construct a $(k-1)$-plane $P_x$ in $B^{k}(x_0, \delta)$ passing through $x$ on
which $u$ is affine and $P_x\cap P_z\cap B^k(x_0, \delta)=\emptyset$ if
$z\notin P_x$, see Figure \ref{fig4}. We then construct a foliation of $B^k(x_0, \delta)$, see Figure \ref{fig5} and obtain that the assumptions of Lemma \ref{0continuity} are satisfied along $P_{\gamma(t)}$ for a.e. $t\in (a,b)$ by the same argument as Step 4 and Step 5 of Lemma \ref{0continuity}. It then follows that $\nabla u$ is constant on $P_{\gamma(t)}$ for a.e. $t\in (a,b)$. 

By choosing an initial value for $\gamma$ arbitrary 
close to $x_0$ and applying the co-area formula in a similar manner we can make sure that 
$\nabla u$ is of class $W^{1,2}$ on $\gamma$. Hence we conclude that $\nabla u$ is
$C^{0,1/2}$ on $\gamma$ by the Sobolev embedding theorem. Let $F$ be the set of
$t\in (a,b)$ such that $\nabla u$ is not constant along $P_{\gamma(t)}$, then $\mathcal{H}^1(F)=0$. We modify $\nabla u$ to be constant along $P_{\gamma(t)}$
for each $t\in F$. Note that, 
$$\mathcal{H}^k(\bigcup  \{P_{\gamma(t)}:t\in F\})\leq c(2\delta)^{k-1} \mathcal{H}^1
(\{\gamma(t):t\in F\})= c(2\delta)^{k-1}  \mathcal{H}^1(F)=0  $$
for some constant $c$. Hence $\nabla u$ is $C^{0,1/2}$ up to modification of a set of measure zero in
$B^k(x_0, \delta)$. Moreover, $\nabla u$ is constant on $P_{\gamma(t)}$ for all $t$, which foliates $B^k(x_0,\delta)$. 
Thus $\nabla u$ is constant on any region on which $u$ is affine. Therefore, $\nabla u$ is constant either
on a $(k-1)$-plane or $k$-dimensional region in $B^k(x_0, \delta)$. This implies that the conclusions of 
Proposition \ref{0induction} under the induction hypothesis are true and hence the inductive step is established. As a conclusion the proofs of Proposition \ref{0induction} 
and Theorem \ref{0deve-regu} are complete. 
\hspace*{\fill} $\Box$

\section{Density: Proof of Theorem \ref{0density}}\label{sdensity}
In this section we show that isometric immersions smooth up to the boundary are strongly dense in
$I^{2,2}(\Omega, \bbbr^{n+1})$ if $\Omega \subset \bbbr^n$ is a convex $C^1$ domain.
Note that it is sufficient to prove that $I^{2,2} \cap C^\infty (\Omega, \R^{n+1})$ is strongly dense in 
$I^{2,2}  (\Omega, \R^{n+1})$. Having this result at hand, and since $\Omega$ is assumed convex, the approximating sequence 
can be easily rescaled to be smooth up to 
the boundary. 
 
\subsection{Foliations of the domain.}\label{sfoliation}
We have argued in the proof of Theorem \ref{0deve-regu} in section \ref{sregu} that for every
maximal region $U\subset \Omega$ on which $u$ is affine, $\partial U\cap \Omega=
\bigcup_{x\in \partial U\cap \Omega}P^U_x$, where $P^U_x$ is some $(n-1)$-plane
in $\Omega$ containing $x$ with the property that for $x_1, x_2\in \partial
U\cap \Omega$, $P^U_{x_1}=P^U_{x_2}$ if $x_2\in P^U_{x_1}$ and $P^U_{x_1}\cap
P^U_{x_2}\cap \Omega=\emptyset$ if $x_2\notin P^U_{x_1}$

We say a maximal region on which $u$ is affine is a \textit{body} if its
boundary contains more than \textit{two} different $(n-1)$-planes in $\Omega$.
\begin{lemma}\label{0body}
It is sufficient to prove Theorem \ref{0density} for a function in
$I^{2,2}(\Omega, \bbbr^{n+1})$ with a finite number of bodies.
\end{lemma}
{\em Proof.} The proof is similar to the proof of \cite[Lemma 3.8]{Pak}
and is omitted for brevity. \hspace*{\fill} $\Box$

\begin{figure}[ht]
\centering
\begin{tikzpicture}[line cap=round,line join=round,>=triangle
45,x=1.0cm,y=1.0cm]
\clip(-4.28,-0.12) rectangle (4.1,6.26);
\fill[line width=0pt,color=qqffff,fill=qqffff,fill opacity=0.5] (-0.915,5.805) -- (-1.89,0.576) -- (-0.074,0.226) -- (1.965,0.634) -- cycle;
\fill[line width=0pt,color=qqffff,fill=qqffff,fill opacity=0.5] (0.164,5.856) -- (2.407,0.865) -- (3.742,2.63) -- cycle;
\fill[line width=0pt,color=qqffff,fill=qqffff,fill opacity=0.5] (-2.897,4.991) -- (-2.687,5.143) -- (-2.605,5.196) -- (-2.534,5.239) -- (-2.461,5.281) -- (-2.682,0.995) -- (-2.863,1.127) -- (-2.998,1.238) -- (-3.138,1.365) -- (-3.253,1.483) -- (-3.378,1.627) -- (-3.487,1.77) -- (-3.576,1.903) -- (-3.66,2.049) -- (-3.747,2.229) -- (-3.814,2.403) -- cycle;
\fill[line width=0pt,color=qqffff,fill=qqffff,fill opacity=0.5] (0.872,5.768) -- (1.097,5.719) -- (1.341,5.653) -- (1.544,5.587) -- (1.718,5.523) -- (1.973,5.413) -- (2.195,5.301) -- (2.344,5.216) -- (2.474,5.134) -- (2.599,5.048) -- (2.727,4.952) -- (2.874,4.83) -- (2.985,4.728) -- (3.121,4.587) -- (3.286,4.39) -- (3.409,4.214) -- (3.501,4.058) -- (3.596,3.864) -- (3.667,3.676) -- (3.741,3.392) -- cycle;
\draw [rotate around={-1.2:(-0.075,3.045)},color=qqqqqq] (-0.075,3.045) ellipse (3.852cm and 2.819cm);
\draw [color=qqttzz] (-3.21,4.713)-- (-3.916,2.873);
\draw [color=qqttzz] (-2.005,5.503)-- (-2.516,0.888);
\draw [color=qqttzz] (-1.618,5.643)-- (-2.298,0.764);
\draw [color=qqttzz] (-0.35,5.86)-- (2.229,0.764);
\draw [dash pattern=on 3pt off 3pt,color=qqqqff] (1.341,5.653)-- (3.596,3.864);
\draw [color=qqttzz] (2.71,1.07)-- (3.683,2.389);
\draw [color=qqttzz] (0.341,0.238)-- (1.749,0.544);
\draw [color=qqttzz] (-2.094,0.664)-- (-1.259,5.739);
\draw [color=qqttzz] (-2.897,4.991)-- (-3.814,2.403);
\draw [color=qqttzz] (-2.461,5.281)-- (-2.682,0.995);
\draw [color=qqttzz] (-0.915,5.805)-- (-1.89,0.576);
\draw [color=qqttzz] (-0.915,5.805)-- (1.965,0.634);
\draw [color=qqttzz] (-1.89,0.576)-- (-0.074,0.226);
\draw [color=qqttzz] (-0.074,0.226)-- (1.965,0.634);
\draw [color=qqttzz] (0.164,5.856)-- (2.407,0.865);
\draw [color=qqttzz] (0.164,5.856)-- (3.742,2.63);
\draw [color=qqttzz] (-1.573,0.462)-- (-0.397,0.239);
\draw [color=qqttzz] (0.872,5.768)-- (3.741,3.392);
\draw [color=qqttzz] (0.49,5.828)-- (3.777,2.99);
\draw [color=qqttzz] (3.742,2.63)-- (2.407,0.865);
\draw [dash pattern=on 3pt off 3pt,color=qqqqff] (1.973,5.413)-- (3.286,4.39);
\draw [dash pattern=on 3pt off 3pt,color=qqqqff] (-2.687,5.143)-- (-3.487,1.77);
\draw [dash pattern=on 3pt off 3pt,color=qqqqff] (-2.605,5.196)-- (-3.253,1.483);
\draw [dash pattern=on 3pt off 3pt,color=qqqqff] (-2.534,5.239)-- (-2.998,1.238);
\draw [dash pattern=on 3pt off 3pt,color=qqqqff] (-2.769,5.086)-- (-3.66,2.049);
\draw [color=qqttzz](-2.26,4.02) node[anchor=north west] {$x $};
\draw [color=qqttzz](-1.64,4.04) node[anchor=north west] {$z $};
\draw [color=qqttzz](-2.42,5.08) node[anchor=north west] {$P_x$};
\draw [color=qqttzz](-1.8,4.96) node[anchor=north west] {$P_z $};
\draw [color=qqqqqq](-0.56,2.7) node[anchor=north west] {$B_1 $};
\draw [color=qqqqqq](2.28,3.22) node[anchor=north west] {$B_2$};
\begin{scriptsize}
\fill [color=qqttzz] (-2.188,3.853) circle (1.0pt);
\fill [color=qqttzz] (-1.578,3.799) circle (1.0pt);
\end{scriptsize}
\end{tikzpicture}
\caption{Construction of global foliations in $\widetilde{\Omega}$.}
\label{fig11}
\end{figure}
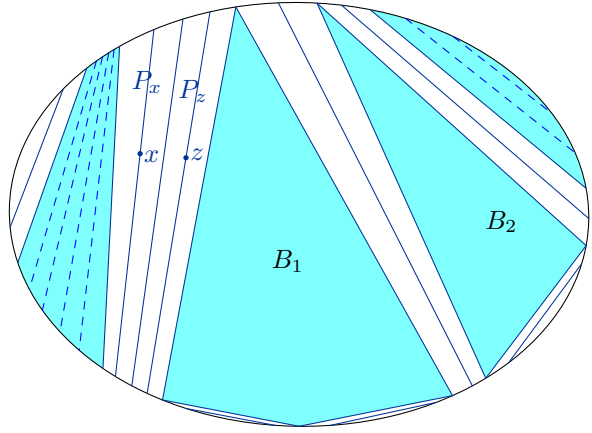

Now we can just assume $u\in I^{2,2}(\Omega, \bbbr^{n+1})$ has finite number of
bodies. Each body is closed and so is therefore their union, whose complement we denote by 
$\widetilde{\Omega}$. Note that now for every $n$-dimensional maximal-affine region $U\subset
\widetilde{\Omega}$, $\partial U\cap \widetilde{\Omega}$ consists of at most two
$(n-1)$-planes.

Similarly as in the proof of Lemma \ref{0continuity}, for every $x\in \widetilde{\Omega}$, 
we will construct an $(n-1)$-plane
$P_x$ in $\widetilde{\Omega}$ passing through it on which $\nabla u $ is
constant and $P_x\cap P_z\cap \widetilde{\Omega}=\emptyset$ if $z\notin P_x$. To apply the same construction in Lemma \ref{0continuity}, 
we makes use of Theorem \ref{0deve-regu} and the fact that $\partial U\cap \widetilde{\Omega}$ consists of at most two
$(n-1)$-planes, see in Figure \ref{fig11}.
 
For every $x\in \widetilde{\Omega}$, we define the normal vector field
$\mathbf{N}(x)$ as the unit vector orthogonal to the family $P_x$ constructed
above. Since none of the $P_x$'s intersect inside $\widetilde{\Omega}$ 
we can choose an orientation such that $\mathbf{N}$ is a Lipschitz vector
fields. The ODE,
\begin{equation}\label{ode}
 \gamma'(t)= \mathbf{N}(\gamma(t)) \quad \gamma(0)=x_0
\end{equation}
has a unique solution $\gamma:(a,b)\rightarrow \widetilde{\Omega}$ for some
interval $(a,b)\subset \bbbr$ containing $0$. Note that $P_x=P_{\gamma(t)}$ if $x\in
P_{\gamma(t)}$, therefore, $\{P_{\gamma(t)}\}_{t\in (a,b)}$ is a local foliation
of $\widetilde{\Omega}$ such that $\nabla u$ is constant on $P_{\gamma(t)}$ for
all $t\in (a,b)$ (Figure \ref{fig12}).
 
\subsection{Leading curves in the domain}\label{sleading}

\begin{figure}[ht]
\centering
\begin{tikzpicture}[line cap=round,line join=round,>=triangle
45,x=1.0cm,y=1.0cm]
\clip(-4.3,-0.22) rectangle (4.32,6.24);
\fill[line width=0pt,color=qqffff,fill=qqffff,fill opacity=0.5] (-0.855,5.745) -- (-1.83,0.516) -- (-0.014,0.166) -- (2.025,0.574) -- cycle;
\fill[line width=0pt,color=qqffff,fill=qqffff,fill opacity=0.5] (0.224,5.796) -- (2.467,0.805) -- (3.802,2.57) -- cycle;
\fill[line width=0pt,color=qqffff,fill=qqffff,fill opacity=0.5] (-2.837,4.931) -- (-2.627,5.083) -- (-2.545,5.136) -- (-2.474,5.179) -- (-2.401,5.221) -- (-2.622,0.935) -- (-2.803,1.067) -- (-2.938,1.178) -- (-3.078,1.305) -- (-3.193,1.423) -- (-3.318,1.567) -- (-3.427,1.71) -- (-3.516,1.843) -- (-3.6,1.989) -- (-3.687,2.169) -- (-3.754,2.343) -- cycle;
\fill[line width=0pt,color=qqffff,fill=qqffff,fill opacity=0.5] (0.932,5.708) -- (1.157,5.659) -- (1.401,5.593) -- (1.604,5.527) -- (1.778,5.463) -- (2.033,5.353) -- (2.255,5.241) -- (2.404,5.156) -- (2.534,5.074) -- (2.659,4.988) -- (2.787,4.892) -- (2.934,4.77) -- (3.045,4.668) -- (3.181,4.527) -- (3.346,4.33) -- (3.469,4.154) -- (3.561,3.998) -- (3.656,3.804) -- (3.727,3.616) -- (3.801,3.332) -- cycle;
\draw [rotate around={-1.2:(-0.015,2.985)},color=qqqqqq] (-0.015,2.985) ellipse (3.852cm and 2.819cm);
\draw [color=qqttzz] (-3.11,4.693)-- (-3.861,2.862);
\draw [color=qqttzz] (-2.114,5.369)-- (-2.456,0.828);
\draw [color=qqttzz] (-1.681,5.543)-- (-2.238,0.704);
\draw [color=qqttzz] (-0.29,5.8)-- (2.289,0.704);
\draw [dash pattern=on 3pt off 3pt,color=qqqqff] (1.401,5.593)-- (3.656,3.804);
\draw [color=qqttzz] (2.77,1.01)-- (3.743,2.329);
\draw [color=qqttzz] (0.28,0.171)-- (1.801,0.481);
\draw [color=qqttzz] (-2.034,0.604)-- (-1.27,5.663);
\draw [color=qqttzz] (-2.837,4.931)-- (-3.754,2.343);
\draw [color=qqttzz] (-2.401,5.221)-- (-2.622,0.935);
\draw [color=qqttzz] (-0.855,5.745)-- (-1.83,0.516);
\draw [color=qqttzz] (-0.855,5.745)-- (2.025,0.574);
\draw [color=qqttzz] (-1.83,0.516)-- (-0.014,0.166);
\draw [color=qqttzz] (-0.014,0.166)-- (2.025,0.574);
\draw [color=qqttzz] (0.224,5.796)-- (2.467,0.805);
\draw [color=qqttzz] (0.224,5.796)-- (3.802,2.57);
\draw [color=qqttzz] (-1.513,0.402)-- (-0.31,0.177);
\draw [color=qqttzz] (0.932,5.708)-- (3.801,3.332);
\draw [color=qqttzz] (0.55,5.768)-- (3.837,2.93);
\draw [color=qqttzz] (3.802,2.57)-- (2.467,0.805);
\draw [dash pattern=on 3pt off 3pt,color=qqqqff] (2.033,5.353)-- (3.346,4.33);
\draw [dash pattern=on 3pt off 3pt,color=qqqqff] (-2.627,5.083)-- (-3.427,1.71);
\draw [dash pattern=on 3pt off 3pt,color=qqqqff] (-2.545,5.136)-- (-3.193,1.423);
\draw [dash pattern=on 3pt off 3pt,color=qqqqff] (-2.474,5.179)-- (-2.938,1.178);
\draw [dash pattern=on 3pt off 3pt,color=qqqqff] (-2.709,5.026)-- (-3.6,1.989);
\draw [color=qqqqqq](-0.46,3.14) node[anchor=north west] {$B_1$};
\draw [color=qqqqqq](1.98,3.62) node[anchor=north west] {$B_2$};
\draw [shift={(-2.35,6.25)},color=qqqqqq]  plot[domain=4.204:4.662,variable=\t]({1*2.817*cos(\t r)+0*2.817*sin(\t r)},{0*2.817*cos(\t r)+1*2.817*sin(\t r)});
\draw [color=qqqqqq] (-1.007,0.357)-- (-1.027,0.275);
\draw [color=qqqqqq] (0.962,0.361)-- (0.981,0.252);
\draw [color=qqqqqq] (3.11,1.656)-- (3.299,1.515);
\draw [shift={(2.704,0.242)},color=qqqqqq]  plot[domain=1.989:2.193,variable=\t]({1*3.782*cos(\t r)+0*3.782*sin(\t r)},{0*3.782*cos(\t r)+1*3.782*sin(\t r)});
\draw [shift={(-3.249,9.062)},color=qqqqqq]  plot[domain=5.536:5.678,variable=\t]({1*7.366*cos(\t r)+0*7.366*sin(\t r)},{0*7.366*cos(\t r)+1*7.366*sin(\t r)});
\draw [shift={(-2.953,-5.5)},color=qqqqqq]  plot[domain=1.386:1.519,variable=\t]({1*8.943*cos(\t r)+0*8.943*sin(\t r)},{0*8.943*cos(\t r)+1*8.943*sin(\t r)});
\draw [color=qqqqqq](-2.46,4.2) node[anchor=north west] {$ \gamma_1$};
\draw [color=qqqqqq](0.72,3.76) node[anchor=north west] {$\gamma_2$};
\draw [color=qqqqqq](2.36,4.58) node[anchor=north west] {$\gamma_3$};
\draw [color=qqqqqq](-1.12,1.16) node[anchor=north west] {$\gamma_4$};
\draw [color=qqqqqq](0.4,1.26) node[anchor=north west] {$\gamma_5$};
\draw [color=qqqqqq](2.56,2.46) node[anchor=north west] {$\gamma_6$};
\end{tikzpicture}
\caption{}
\label{fig12}
\end{figure}

\begin{figure}[ht]
\centering
\begin{tikzpicture}[line cap=round,line join=round,>=triangle
45,x=1.0cm,y=1.0cm]
\clip(-4.54,0.2) rectangle (7.98,6.26);
\fill[line width=0pt,color=qqffff,fill=qqffff,fill opacity=0.5] (-1.887,5.198) -- (-1.308,1.299) -- (1.107,2.119) -- cycle;
\fill[line width=0pt,color=qqffff,fill=qqffff,fill opacity=0.1] (3.04,1.18) -- (3.35,5.05) -- (5.13,5.93) -- (4.78,2.15) -- cycle;
\fill[line width=0pt,color=qqffff,fill=qqffff,fill opacity=0.1] (3.14,1.06) -- (2.42,4.87) -- (4.19,6.06) -- (4.86,2.33) -- cycle;
\fill[line width=0pt,color=qqffff,fill=qqffff,fill opacity=0.1] (3.4,1.06) -- (3.38,4.89) -- (5.08,6) -- (5.07,2.25) -- cycle;
\fill[line width=0pt,color=qqffff,fill=qqffff,fill opacity=0.1] (4.1,1.02) -- (4.47,4.88) -- (6.18,5.99) -- (5.77,2.21) -- cycle;
\fill[line width=0pt,color=qqffff,fill=qqffff,fill opacity=0.1] (5.22,1.01) -- (4.85,4.78) -- (6.55,5.88) -- (6.89,2.2) -- cycle;
\fill[line width=0pt,color=qqffff,fill=qqffff,fill opacity=0.1] (6.01,0.91) -- (5.38,4.57) -- (7.03,5.99) -- (7.62,2.41) -- cycle;
\draw [color=qqttzz] (2.42,4.87)-- (4.19,6.06);
\draw [color=qqttzz] (2.42,4.87)-- (3.14,1.06);
\draw [line width=0.2pt,dash pattern=on 3pt off 3pt,color=qqttzz] (3.14,1.06)-- (4.86,2.33);
\draw [line width=0.2pt,dash pattern=on 3pt off 3pt,color=qqttzz] (4.19,6.06)-- (4.86,2.33);
\draw [color=qqttzz] (3.35,5.05)-- (5.13,5.93);
\draw [color=qqttzz] (3.35,5.05)-- (3.04,1.18);
\draw [line width=0.2pt,dash pattern=on 3pt off 3pt,color=qqttzz] (3.04,1.18)-- (4.78,2.15);
\draw [line width=0.2pt,dash pattern=on 3pt off 3pt,color=qqttzz] (5.13,5.93)-- (4.78,2.15);
\draw [color=qqttzz] (3.38,4.89)-- (5.08,6);
\draw [color=qqttzz] (3.38,4.89)-- (3.4,1.06);
\draw [line width=0.2pt,dash pattern=on 3pt off 3pt,color=qqttzz] (3.4,1.06)-- (5.07,2.25);
\draw [line width=0.2pt,dash pattern=on 3pt off 3pt,color=qqttzz] (5.08,6)-- (5.07,2.25);
\draw [color=qqttzz] (4.47,4.88)-- (6.18,5.99);
\draw [color=qqttzz] (4.47,4.88)-- (4.1,1.02);
\draw [line width=0.2pt,dash pattern=on 3pt off 3pt,color=qqttzz] (4.1,1.02)-- (5.77,2.21);
\draw [line width=0.2pt,dash pattern=on 3pt off 3pt,color=qqttzz] (6.18,5.99)-- (5.77,2.21);
\draw [color=qqttzz] (4.85,4.78)-- (6.55,5.88);
\draw [color=qqttzz] (4.85,4.78)-- (5.22,1.01);
\draw [line width=0.2pt,dash pattern=on 3pt off 3pt,color=qqttzz] (5.22,1.01)-- (6.89,2.2);
\draw [line width=0.2pt,dash pattern=on 3pt off 3pt,color=qqttzz] (6.55,5.88)-- (6.89,2.2);
\draw [color=qqttzz] (5.38,4.57)-- (7.03,5.99);
\draw [color=qqttzz] (5.38,4.57)-- (6.01,0.91);
\draw [color=qqttzz] (6.01,0.91)-- (7.62,2.41);
\draw [color=qqttzz] (7.03,5.99)-- (7.62,2.41);
\draw [line width=0.2pt,dash pattern=on 3pt off 3pt,color=qqttzz] (3.063,1.467)-- (4.818,2.563);
\draw [color=qqttzz] (5.08,6)-- (5.08,5.905);
\draw [color=qqttzz] (3.04,1.18)-- (3.11,1.219);
\draw [color=qqttzz] (6.55,5.88)-- (6.576,5.599);
\draw [color=qqttzz] (5.22,1.01)-- (5.908,1.501);
\draw [color=qqttzz] (6.18,5.99)-- (6.139,5.614);
\draw [color=qqttzz] (4.1,1.02)-- (5.146,1.765);
\draw [color=qqttzz] (5.08,5.905)-- (5.13,5.93);
\draw [color=qqttzz] (3.4,1.06)-- (4.155,1.598);
\draw [color=qqttzz] (3.14,1.06)-- (3.399,1.251);
\draw [color=qqttzz] (4.19,6.06)-- (4.293,5.486);
\draw [color=qqttzz] (3.063,1.467)-- (3.397,1.676);
\draw [line width=0.2pt,dash pattern=on 3pt off 3pt,color=qqqqqq] (3.69,3.76)-- (3.88,3.8);
\draw [line width=0.2pt,dash pattern=on 3pt off 3pt,color=qqqqqq] (3.88,3.8)-- (4.09,3.79);
\draw [color=qqqqqq] (4.09,3.79)-- (4.348,3.749);
\draw [line width=0.2pt,dash pattern=on 3pt off 3pt,color=qqqqqq] (4.348,3.749)-- (4.51,3.76);
\draw [line width=0.2pt,dash pattern=on 3pt off 3pt,color=qqqqqq] (4.51,3.76)-- (4.72,3.77);
\draw [line width=0.2pt,dash pattern=on 3pt off 3pt,color=qqqqqq] (4.72,3.77)-- (4.95,3.74);
\draw [line width=0.2pt,dash pattern=on 3pt off 3pt,color=qqqqqq] (4.95,3.74)-- (5.02,3.7);
\draw [color=qqqqqq] (5.02,3.7)-- (5.19,3.66);
\draw [line width=0.2pt,dash pattern=on 3pt off 3pt,color=qqqqqq] (5.19,3.66)-- (5.4,3.59);
\draw [line width=0.2pt,dash pattern=on 3pt off 3pt,color=qqqqqq] (5.4,3.59)-- (5.72,3.55);
\draw [line width=0.2pt,dash pattern=on 3pt off 3pt,color=qqqqqq] (5.72,3.55)-- (5.96,3.56);
\draw [line width=0.2pt,dash pattern=on 3pt off 3pt,color=qqqqqq] (5.96,3.56)-- (6.35,3.63);
\draw [line width=0.2pt,dash pattern=on 3pt off 3pt,color=qqqqqq] (6.35,3.63)-- (6.58,3.7);
\draw [->,>=stealth',color=qqqqqq] (6.58,3.7) -- (6.47,4.3);
\draw [->,>=stealth',color=qqqqqq] (6.58,3.7) -- (6.19,3.38);
\draw [->,>=stealth',line width=0.2pt,dash pattern=on 3pt off 3pt,color=qqqqqq] (5.96,3.56) -- (5.97,4.19);
\draw [->,>=stealth',line width=0.2pt,dash pattern=on 3pt off 3pt,color=qqqqqq] (5.96,3.56) -- (5.59,3.23);
\draw [->,>=stealth',line width=0.2pt,dash pattern=on 3pt off 3pt,color=qqqqqq] (5.02,3.7) -- (5.21,4.36);
\draw [->,>=stealth',color=qqqqqq] (5.02,3.7) -- (4.63,3.31);
\draw [->,>=stealth',line width=0.2pt,dash pattern=on 3pt off 3pt,color=qqqqqq] (4.51,3.76) -- (4.46,4.37);
\draw [->,>=stealth',line width=0.2pt,dash pattern=on 3pt off 3pt,color=qqqqqq] (4.51,3.76) -- (4.01,3.46);
\draw [->,>=stealth',line width=0.2pt,dash pattern=on 3pt off 3pt,color=qqqqqq] (4.09,3.79) -- (4.14,4.4);
\draw [->,>=stealth',line width=0.2pt,dash pattern=on 3pt off 3pt,color=qqqqqq] (4.09,3.79) -- (3.62,3.42);
\draw [->,>=stealth',line width=0.2pt,dash pattern=on 3pt off 3pt,color=qqqqqq] (3.69,3.76) -- (3.52,4.32);
\draw [->,>=stealth',line width=0.2pt,dash pattern=on 3pt off 3pt,color=qqqqqq] (3.69,3.76) -- (3.24,3.44);
\draw [color=qqttzz] (5.13,5.93)-- (5.078,5.372);
\draw [color=qqttzz] (5.078,5.372)-- (5.078,5.275);
\draw [color=qqttzz] (4.483,5.61)-- (5.078,5.372);
\draw [rotate around={0:(-1.085,3.29)},color=qqqqqq] (-1.085,3.29) ellipse (2.706cm and 1.998cm);
\draw [color=qqttzz] (-3.699,3.806)-- (-3.478,2.358);
\draw [color=qqttzz] (-3.33,2.175)-- (-3.215,4.522);
\draw [color=qqttzz] (-3.121,4.605)-- (-2.764,1.723);
\draw [color=qqttzz] (-2.525,1.598)-- (-2.494,4.995);
\draw [color=qqttzz] (-2.167,5.121)-- (-1.904,1.386);
\draw [color=qqttzz] (-1.887,5.198)-- (-1.308,1.299);
\draw [color=qqttzz] (-1.887,5.198)-- (1.107,2.119);
\draw [color=qqttzz] (-1.308,1.299)-- (1.107,2.119);
\draw [color=qqttzz] (-0.835,1.301)-- (0.825,1.875);
\draw [color=qqttzz] (-1.434,5.271)-- (1.487,2.669);
\draw [color=qqttzz] (-0.574,5.252)-- (1.62,3.24);
\draw [color=qqttzz] (0.297,5.007)-- (1.527,3.81);
\draw [dash pattern=on 3pt off 3pt,color=qqttzz] (-3.89,5.06)-- (-3.699,3.806);
\draw [dash pattern=on 3pt off 3pt,color=qqttzz] (-3.478,2.358)-- (-3.27,0.99);
\draw [dash pattern=on 3pt off 3pt,color=qqttzz] (-3.17,5.45)-- (-3.215,4.522);
\draw [dash pattern=on 3pt off 3pt,color=qqttzz] (-3.33,2.175)-- (-3.39,0.94);
\draw [dash pattern=on 3pt off 3pt,color=qqttzz] (-3.22,5.4)-- (-3.121,4.605);
\draw [dash pattern=on 3pt off 3pt,color=qqttzz] (-2.764,1.723)-- (-2.67,0.97);
\draw [dash pattern=on 3pt off 3pt,color=qqttzz] (-2.49,5.5)-- (-2.494,4.995);
\draw [dash pattern=on 3pt off 3pt,color=qqttzz] (-2.525,1.598)-- (-2.53,0.99);
\draw [dash pattern=on 3pt off 3pt,color=qqttzz] (-2.2,5.59)-- (-2.167,5.121);
\draw [dash pattern=on 3pt off 3pt,color=qqttzz] (-1.904,1.386)-- (-1.87,0.91);
\draw [dash pattern=on 3pt off 3pt,color=qqttzz] (-1.95,5.62)-- (-1.887,5.198);
\draw [dash pattern=on 3pt off 3pt,color=qqttzz] (-1.308,1.299)-- (-1.25,0.91);
\draw [line width=0.2pt,dash pattern=on 3pt off 3pt,color=qqqqqq] (3.69,3.76)-- (3.295,3.665);
\draw [line width=0.2pt,dash pattern=on 3pt off 3pt,color=qqqqqq] (3.295,3.665)-- (2.626,3.482);
\draw [color=qqqqqq] (2.626,3.482)-- (2.376,3.406);
\draw [color=qqqqqq](-4.42,3.28) node[anchor=north west] {\scriptsize $
\gamma(0)$};
\draw [color=qqqqqq](-1.58,3.42) node[anchor=north west] {\scriptsize $
\gamma(\ell) $};
\draw [color=qqqqqq](-2.8,3.32) node[anchor=north west] {\scriptsize $
\gamma(t)$};
\draw [color=qqqqqq](-2.88,1.32) node[anchor=north west] {\scriptsize $
F_\gamma(t) $};
\draw [color=qqqqqq](-1.62,1.2) node[anchor=north west] {\scriptsize $
F_\gamma(\ell)$};
\draw [color=qqqqqq](1.82,3.56) node[anchor=north west] {\scriptsize $
\gamma(0)$};
\draw [color=qqqqqq](4.94,3.86) node[anchor=north west] {\scriptsize $ \gamma(t)
$};
\draw [color=qqqqqq](6.46,3.94) node[anchor=north west] {\scriptsize $
\gamma(\ell) $};
\draw [color=qqqqqq](6.54,4.42) node[anchor=north west] {\scriptsize $
\mathbf{N}_1 (\ell)$};
\draw [color=qqqqqq](6.86,2.9) node[anchor=north west] {\scriptsize $
F_\gamma(\ell) $};
\draw [color=qqqqqq](4.16,2.02) node[anchor=north west] {\scriptsize $
F_\gamma(t)$};
\draw [color=qqqqqq](3.76,0.8) node[anchor=north west] {\scriptsize leading
fronts};
\draw [color=qqqqqq](0.96,4.94) node[anchor=north west] {\scriptsize leading
curve};
\draw [->,>=stealth',color=qqqqqq] (4.29,0.81) -- (4.35,1.54);
\draw [->,>=stealth',color=qqqqqq] (4.29,0.81) -- (6.24,1.36);
\draw [->,>=stealth',color=qqqqqq] (4.29,0.81) -- (3.52,1.341);
\draw [->,>=stealth',color=qqqqqq] (1.97,4.53) -- (2.578,3.467);
\draw [color=qqqqqq](5.92,3.58) node[anchor=north west] {\scriptsize $
\mathbf{N}_2(\ell)$};
\draw [shift={(1.429,3.257)},color=qqqqqq]  plot[domain=1.102:2.01,variable=\t]({1*2.324*cos(\t r)+0*2.324*sin(\t r)},{0*2.324*cos(\t r)+1*2.324*sin(\t r)});
\draw [color=qqqqqq] (2.48,5.33)-- (2.38,5.52);
\draw [color=qqqqqq] (2.48,5.33)-- (2.28,5.27);
\draw [shift={(-3.359,1.575)},color=qqqqqq]  plot[domain=1.522:1.842,variable=\t]({1*1.537*cos(\t r)+0*1.537*sin(\t r)},{0*1.537*cos(\t r)+1*1.537*sin(\t r)});
\draw [shift={(-2.539,-0.081)},color=qqqqqq]  plot[domain=1.562:1.694,variable=\t]({1*3.231*cos(\t r)+0*3.231*sin(\t r)},{0*3.231*cos(\t r)+1*3.231*sin(\t r)});
\draw [shift={(-3.186,5.125)},color=qqqqqq]  plot[domain=4.664:4.836,variable=\t]({1*2.018*cos(\t r)+0*2.018*sin(\t r)},{0*2.018*cos(\t r)+1*2.018*sin(\t r)});
\draw [shift={(-2.457,9.233)},color=qqqqqq]  plot[domain=4.704:4.783,variable=\t]({1*6.083*cos(\t r)+0*6.083*sin(\t r)},{0*6.083*cos(\t r)+1*6.083*sin(\t r)});
\draw [shift={(-2.43,8.848)},color=qqqqqq]  plot[domain=4.783:4.86,variable=\t]({1*5.697*cos(\t r)+0*5.697*sin(\t r)},{0*5.697*cos(\t r)+1*5.697*sin(\t r)});
\begin{scriptsize}
\fill [color=qqqqqq] (3.69,3.76) circle (1.0pt);
\fill [color=qqqqqq] (4.09,3.79) circle (1.0pt);
\fill [color=qqqqqq] (4.51,3.76) circle (1.0pt);
\fill [color=qqqqqq] (5.02,3.7) circle (1.0pt);
\fill [color=qqqqqq] (5.96,3.56) circle (1.0pt);
\fill [color=qqqqqq] (6.58,3.7) circle (1.0pt);
\fill [color=qqqqqq] (2.376,3.406) circle (1.0pt);
\fill [color=qqqqqq] (-2.511,3.15) circle (1.0pt);
\end{scriptsize}
\end{tikzpicture}
\caption{Leading curve and leading fronts.}
\label{fig13}
\end{figure}
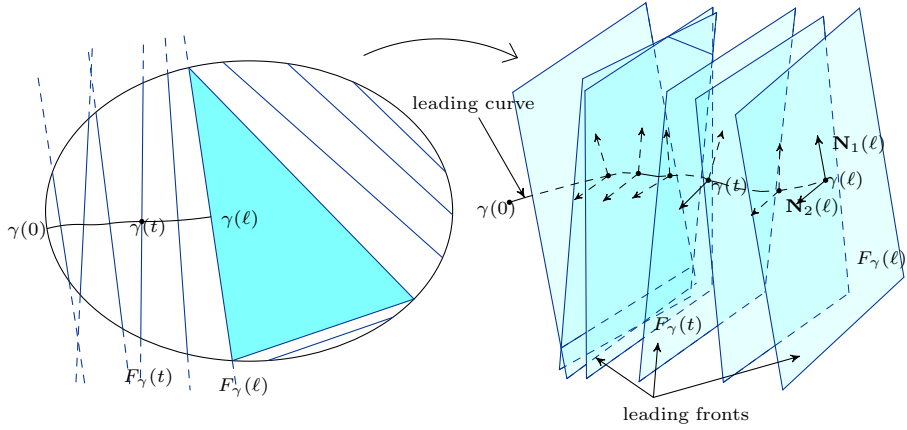

\begin{definition}\label{0leadingcurve}
Let $\{P_x\}_{x\in \widetilde{\Omega}}$ be a family of $(n-1)$-planes in
$\widetilde{\Omega}$ passing through $x$ on which $\nabla u$ is constant, satisfying  
$P_x\cap P_z\cap \widetilde{\Omega}=\emptyset$ if $z\notin P_x$ and $P_x = P_z$ if $z\in P_x$. 
We say that a curve $\gamma\in C^{1,1}([0,\ell],  \widetilde{\Omega})$ parametrized by arclength  
is a \textit{leading curve} if it is orthogonal at any possible point of intersection $z\in \gamma([0,\ell])\cap P_x$ to $P_x= P_z$ for all $x \in \widetilde{\Omega}$ (Figure \ref{fig13}).
\end{definition}
It is easy to see that $\gamma$ constructed in Subsection \ref{sfoliation} when
restricted to the interval $[0,\ell]$ is a leading curve, since by the ODE (\ref{ode}),
$|\gamma'|=1$ and $|\gamma''|$ is bounded as $\mathbf{N}$ is Lipschitz.

\begin{definition}\label{0leadingfront}
The $(n-1)$-dimensional hyperplane $F_{\gamma}(t)$ orthogonal to $\gamma(t)$ at
$t\in [0,\ell]$ is called the \textit{leading front} of $\gamma$ at $t\in [0,\ell]$
(Figure \ref{fig13}).
\end{definition}

\begin{remark}\label{0remark}
It then follows from the definition of the leading curve that $F_\gamma(t)\cap
F_\gamma(\tilde{t})\cap \widetilde{\Omega}=\emptyset$ for all $t,\tilde{t}\in
[0,\ell]$ such that $t\neq \tilde{t}$. Moreover,
$F_\gamma(t)\subset \widetilde{\Omega}$, otherwise,
$F_\gamma(t)\cap B\neq \emptyset$ where $B$ is one of the bodies in
$\Omega\setminus \widetilde{\Omega}$. Since $\nabla u$, being continuous, is
constant on $F_\gamma(t)\cap \widetilde{\Omega}$ and $B$, it must be constant on
their convex hull, which is again a body, contradiction to that a body is a
maximal region. 
\end{remark}

We say that a curve $\gamma$ covers the domain $A\subset \Omega$ if
$$A\subset  \bigcup \{F_\gamma(t): t\in [0,\ell] \} . $$
By $\Omega(\gamma)$ we refer to the \textit{biggest set covered by} $\gamma$ \textit{in} $\Omega$.
We now restrict our attention to the covered domain $\Omega(\gamma)$. It is
obvious that $\Omega(\gamma)$ is convex since it is bounded by $F_\gamma(0)$,
$F_\gamma(\ell)$ and $\partial \Omega$.

From the construction in subsection \ref{sfoliation}, the $(n-1)$-planes
$P_{\gamma(t)}$ in $\widetilde{\Omega}$, $t\in [0,\ell]$ which constitute a
local foliation of $\widetilde{\Omega}$ are global foliations of
$\Omega(\gamma)$. Moreover, $P_{\gamma(t)}=F_{\gamma}(t)\cap
\Omega(\gamma)=F_{\gamma}(t)\cap \Omega$ for all $t\in [0,\ell]$. We relabel
them $P_\gamma(t)$ to be in consistence of notation and we name them:
\begin{definition}\label{0leadingplane}
The component $P_\gamma (t) :=F_{\gamma}(t)\cap \Omega$ is called the leading
$(n-1)$-planes in $\Omega$ of $\gamma$ at $t\in [0,\ell]$.
\end{definition}

Let $\{\mathbf{N}_i(t)\}_{i=1}^{n-1}$, be an orthonormal
basis for the leading front $F_{\gamma}(t)$ (Figure \ref{fig13}) such that
$\mathbf{N}_i$ is Lipschitz for
all $1\leq i \leq n-1$ and
$\det[\gamma'(t),\mathbf{N}_1(\tilde{t}),\cdots,\mathbf{N}_{n-1}(\tilde{t})]=1$.
It is obvious such orthonormal basis exists because we can pick
$\{\mathbf{N}_i(0)\}_{i=1}^{n-1}$ as an orthonormal
basis for $F_{\gamma}(0)$ that form a positive orientation with $\gamma'(0)$
and then move this frame along $\gamma$ in an orientation preserving way (note
that $\gamma$ is not a closed curve so this is possible).
Let
$\Phi:[0,\ell]\times \bbbr^{n-1}\rightarrow \bbbr^n$ be defined as,
\begin{equation}\label{0phi}
\ds \Phi(t,s):=\gamma(t)+ \sum_{i=1}^{n-1} s_i\mathbf{N}_i(t) ,
\end{equation}
where $s=(s_1,\cdots,s_{n-1})$. Then we can represent the leading front at $t\in
[0,\ell]$ as,
\begin{equation}\label{0leadingfront1}
F_\gamma(t)= \{\Phi(t,s), s=(s_1,\cdots,s_{n-1})\in \bbbr^{n-1}\}.
\end{equation}
For each $t\in [0,\ell]$, define the open set,
\begin{equation} \label{0domain1}
\Sigma^\gamma(t)=\{s=(s_1,\cdots,s_{n-1})\in \bbbr^{n-1}: \Phi(t,s)\in \Omega\} .
\end{equation}
It is obvious that $0\in \Sigma^\gamma(t)$, hence it is non-empty. Then we can also parametrize the leading planes as
\begin{equation}\label{0leadingplane2}
P_\gamma (t)=\{\Phi(t,s), s=(s_1,\cdots,s_{n-1})\in \Sigma^\gamma(t)\}.
\end{equation}
Now define,
\begin{equation}\label{0domain3}
\Sigma^\gamma:=\{(t,s), \Phi(t,s)\in \Omega\}.
\end{equation}
Of course we can also write,
$$\Sigma^\gamma=\{(t,s), t\in [0,\ell], s=(s_1,\cdots,s_{n-1})\in
\Sigma^\gamma(t)\}.$$

We will focus on the restriction of $\Phi$ in $\Sigma^\gamma$. However, if no
confusion is caused, we still denote such restriction $\Phi$. It is easy to see that
$\Phi$ maps $\Sigma^\gamma$ into $\Omega(\gamma)$. Indeed, if $x=\Phi(t,s)$ for
some $(t,s) \in \Sigma^\gamma$, by definition of $\Sigma^\gamma$, $\Phi(t,s)\in
\Omega$. On the other hand, $\Phi(t,s)\in F_\gamma(t)$, thus, $x=\Phi(t,s)\in
F_\gamma(t)\cap\Omega \subset \Omega(\gamma)$.
\begin{lemma}\label{0phi2}
$\Phi: \Sigma^\gamma\rightarrow \Omega(\gamma)$ is one-to-one and onto. In
particular,
$$\Omega(\gamma)=\{\Phi(t,s), (t,s)\in \Sigma^\gamma\}=\bigcup\{P_\gamma (t):
t\in [0,\ell]\}.$$
\end{lemma}
{\em Proof.}
We first show one-to-one. Suppose $\Phi(t_1,s_1)=\Phi(t_2,s_2)$ while
$(t_1,s_1)\neq (t_2,s_2)$. Since $s\rightarrow \Phi(t,s)$ is obviously
one-to-one by the definition of $\Phi$, it must be $t_1\neq t_2$. We
have argued in Remark \ref{0remark} that $F_\gamma(t_1)\cap F_\gamma(t_2)\cap
\Omega=\emptyset$. Therefore, $F_\gamma(t_1)\cap F_\gamma(t_2)\cap
\Omega(\gamma)=\emptyset$ since $\Omega(\gamma)\subset \Omega$.
However, $\Phi(t_1,s_1)\in F_\gamma(t_1)$ and $\Phi(t_2,s_2)\in F_\gamma(t_2)$,
contradiction to $\Phi(t_1,s_1)=\Phi(t_2,s_2)$.

We will now show onto. Let $x\in \Omega(\gamma)$, then $x=\Phi(t,s)$ for some
$t\in [0,\ell]$ and $s\in \bbbr^{n-1}$. Since $x\in \Omega(\gamma)$,
$\Phi(t,s)\in \Omega(\gamma)\subset \Omega$, hence $(t,s) \in \Sigma^\gamma$.
The proof is complete.
\hspace*{\fill} $\Box$

Apparently we can rewrite
$\ds \Phi(t,s):=\gamma(t)+ \sum_{i=1}^{n-1} s_i\mathbf{N}_i(t), t\in
[0,\ell], s\in \bbbr^{n-1}$ as
\begin{equation*}
\ds \Phi(t,S\cdot s)=\gamma(t)+S \Big ( \sum_{i=1}^{n-1} s_i\mathbf{N}_i(t) \Big ), t\in [0,\ell], 
s\in \mathbb{S}^{n-2}, S\geq 0.
\end{equation*}
We then rewrite the representation of leading front in (\ref{0leadingfront1}) in
an equivalent way:
\begin{equation}\label{0leadingfront3}
F_\gamma(t)= \{\Phi(t,S \cdot s), S\geq 0, s=(s_1,\cdots,s_{n-1})\in \mathbb{S}^{n-2}\}.
\end{equation}
For each $t\in [0,\ell]$ and $s=(s_1,\cdots,s_{n-1})\in \mathbb{S}^{n-2}$, define
the scalar function,
\begin{equation} \label{0length}
S_s^\gamma(t):=\sup\{S\geq 0: \Phi(t,S\cdot s)\in \Omega\} .
\end{equation}
That is, $S_s^\gamma(t)$ is the distance from $\gamma(t)$ to $\partial \Omega$
in the direction $\ds \sum_{i=1}^{n-1} s_i\mathbf{N}_i(t)$.
From the definition of $\Sigma^\gamma(t)$ and $\Sigma^\gamma$,
\begin{equation}\label{0domain5}
\Sigma^\gamma(t)=\{(S\cdot s): s=(s_1,\cdots,s_{n-1})\in \mathbb{S}^{n-2},
0<S<S_s^\gamma(t)\},
\end{equation} and
\begin{equation}\label{0domain6}
\Sigma^\gamma=\{(t,S\cdot s), t\in [0,\ell], s=(s_1,\cdots,s_{n-1})\in \mathbb{S}^{n-2},
0<S<S_s^\gamma(t)\}.
\end{equation}

Since $|\gamma'(t)|=1$, $\gamma
^{\prime \prime}(t)\cdot \gamma'(t)=0$, we can then write
$
\displaystyle \ds \gamma ^{\prime \prime }(t)= \sum_{i=1}^{n-1} \kappa_i(t)\mathbf{N}_i(t).
$ Similarly we can also write
\begin{equation}\label{N}
 \ds \mathbf{N}_i'=\kappa_{i_0}\gamma'+ \sum_{j=1}^{n-1} \kappa_{i_j}\mathbf{N}_j. 
\end{equation}

It is easy to see that $\kappa_{i_0}=-\kappa_i$, $\kappa_{i_i}=0$ and
$\kappa_{i_j}=-\kappa_{j_i}$. These equations can  be written as the matrix equation
\begin{align}\label{R^TR'}
   \left ( \begin{array}{c}  \gamma' \\ \mathbf{N}_1  \\ \mathbf{N}_2\\ \vdots
\\ \mathbf{N}_{n-1} \end{array} \right )'
=  \begin{pmatrix}
 0 & \kappa_1 & \kappa_2 & \cdots & \kappa_{n-1} \\
 -\kappa_1 & 0 & \kappa_{1_2} & \cdots & \kappa_{1_{n-1}} \\
 -\kappa_2 & -\kappa_{1_2} & 0 & \cdots & \kappa_{2_{n-1}} \\
 \vdots & \vdots & \vdots & \cdots  & \vdots \\
 -\kappa_{n-1} & -\kappa_{1_{n-1}} & -\kappa_{2_{n-1}} & \cdots & 0 \\
\end{pmatrix}
\left ( \begin{array}{c}  \gamma' \\ \mathbf{N}_1  \\ \mathbf{N}_2\\ \vdots  \\
\mathbf{N}_{n-1} \end{array} \right )
\end{align}

Given two \textit{non-parallel} leading fronts $F_\gamma(t)$ and
$F_\gamma(\tilde{t})$, denote their $(n-2)$-plane intersection by
$F(t,\tilde{t})$. Given $s=(s_1,\cdots, s_{n-1})\in \mathbb{S}^{n-2}$, define $L_s(t,
\tilde{t})$ as the distance from $\gamma(t)$ to $F(t,\tilde{t})$ along the
direction $\ds \sum_{i=1}^{n-1} s_i\mathbf{N}_i(t)$ (we set $L_s(t,
\tilde{t})=+\infty$ if it does not hit $F(t,\tilde{t})$ along this
direction) (Figure \ref{fig14}). We then define,
\begin{equation}\label{0distance}
L_s^\gamma(t):=\inf \{L_s(t,\tilde{t}): \tilde{t}\neq t\}.
\end{equation}
Since all $F(t,\tilde{t})$ are outside $\Omega$, $L_s^\gamma(t)\geq
S_s^\gamma(t)$ for all $s\in \mathbb{S}^{n-2}$ and $t\in [0,\ell]$.
\begin{figure}[ht]
\centering
\begin{tikzpicture}[line cap=round,line join=round,>=triangle
45,x=1.0cm,y=1.0cm]
\clip(-4.24,0.2) rectangle (5.54,6.28);
\fill[line width=0pt,color=ffffzz,fill=ffffzz,fill opacity=0.5] (-2.64,4.075) --
(-3.9,1.04) -- (-2.06,1.44) -- (-0.671,4.912) -- (-1.675,4.36) -- cycle;
\fill[line width=0pt,color=zzffzz,fill=zzffzz,fill opacity=0.5] (-2.713,3.534)
-- (-1.675,4.36) -- (-0.671,4.912) -- (-0.97,6.05) -- (-2.89,4.21) -- cycle;
\fill[line width=0pt,color=qqffff,fill=qqffff,fill opacity=0.1] (-2.55,4.75) --
(-3.1,0.64) -- (-1.1,1.84) -- (-0.63,5.95) -- cycle;
\fill[line width=0pt,color=ffffzz,fill=ffffzz,fill opacity=0.5] (-2.64,4.075) --
(-1.675,4.36) -- (-0.671,4.912) -- (-0.3,5.84) -- (-2.14,5.28) -- cycle;
\fill[line width=0pt,color=zzffzz,fill=zzffzz,fill opacity=0.5] (-2.713,3.534)
-- (-1.93,0.55) -- (-0.01,2.39) -- (-0.671,4.912) -- (-1.675,4.36) -- cycle;
\draw [color=qqqqqq] (-2.14,5.28)-- (-0.3,5.84);
\draw [color=qqqqqq] (-2.14,5.28)-- (-3.9,1.04);
\draw [line width=0.2pt,dash pattern=on 3pt off 3pt,color=qqqqqq] (-0.3,5.84)--
(-2.06,1.44);
\draw [line width=0.2pt,dash pattern=on 3pt off 3pt,color=qqqqqq] (-3.9,1.04)--
(-2.06,1.44);
\draw [line width=0.2pt,dash pattern=on 3pt off 3pt,color=qqqqqq] (-0.63,5.95)--
(-2.55,4.75);
\draw [line width=0.2pt,dash pattern=on 3pt off 3pt,color=qqqqqq] (-0.63,5.95)--
(-1.1,1.84);
\draw [color=qqqqqq] (-2.55,4.75)-- (-3.1,0.64);
\draw [line width=0.2pt,dash pattern=on 3pt off 3pt,color=qqqqqq] (-3.1,0.64)--
(-1.1,1.84);
\draw [line width=0.2pt,dash pattern=on 3pt off 3pt,color=qqqqqq] (-0.97,6.05)--
(-2.89,4.21);
\draw [color=qqqqqq] (-2.89,4.21)-- (-1.93,0.55);
\draw [line width=0.2pt,dash pattern=on 3pt off 3pt,color=qqqqqq] (-0.97,6.05)--
(-0.01,2.39);
\draw [color=qqqqqq] (-1.93,0.55)-- (-0.01,2.39);
\draw [line width=0.2pt,dash pattern=on 3pt off 3pt,color=qqqqqq]
(-2.713,3.534)-- (-0.725,5.117);
\draw [line width=0.2pt,dash pattern=on 3pt off 3pt,color=qqqqqq]
(-2.772,3.758)-- (-0.671,4.912);
\draw [line width=0.2pt,dash pattern=on 3pt off 3pt,color=qqqqqq]
(-2.64,4.075)-- (-0.78,4.641);
\draw [color=qqqqqq] (-2.64,4.075)-- (-1.675,4.36);
\draw [color=qqqqqq] (-2.713,3.534)-- (-1.675,4.36);
\draw [color=qqqqqq] (-1.675,4.36)-- (-0.671,4.912);
\draw [color=qqqqqq] (-2.772,3.758)-- (-2.676,3.811);
\draw [color=qqqqqq] (-3.9,1.04)-- (-3.021,1.231);
\draw [color=qqqqqq] (-3.1,0.64)-- (-2.109,1.234);
\draw [color=qqqqqq] (-0.671,4.912)-- (-0.01,2.39);
\draw [color=qqqqqq] (-2.89,4.21)-- (-2.583,4.504);
\draw [color=qqqqqq] (-2.55,4.75)-- (-2.293,4.91);
\draw [color=qqqqqq] (-1.603,5.443)-- (-0.97,6.05);
\draw [color=qqqqqq] (-0.97,6.05)-- (-0.9,5.782);
\draw [color=qqqqqq] (-1.286,5.54)-- (-0.63,5.95);
\draw [color=qqqqqq] (-0.63,5.95)-- (-0.655,5.732);
\draw [color=qqqqqq] (-0.3,5.84)-- (-0.671,4.912);
\draw [shift={(-1.402,5.944)},line width=0.2pt,dash pattern=on 3pt off
3pt,color=qqqqqq]  plot[domain=4.03:5.179,variable=\t]({1*2.916*cos(\t
r)+0*2.916*sin(\t r)},{0*2.916*cos(\t r)+1*2.916*sin(\t r)});
\draw [shift={(-1.402,5.944)},color=qqqqqq] 
plot[domain=4.809:5.179,variable=\t]({1*2.916*cos(\t r)+0*2.916*sin(\t
r)},{0*2.916*cos(\t r)+1*2.916*sin(\t r)});
\draw [rotate around={-74.116:(-1.275,3.185)},color=qqqqqq] (-1.275,3.185)
ellipse (1.074cm and 0.863cm);
\draw [color=qqqqqq](-2.1,1.46) node[anchor=north west] {\scriptsize $
F_\gamma(t) $};
\draw [color=qqqqqq](-3.04,1.84) node[anchor=north west] {\scriptsize $
F_\gamma(\tilde{t}) $};
\draw [color=qqqqqq](-0.3,4.94) node[anchor=north west] {\scriptsize $
F(t,\tilde{t})$};
\draw [->,>=stealth',color=qqqqqq] (-0.42,4.65) -- (-1.038,4.868);
\draw [color=qqqqqq](-1.34,3.68) node[anchor=north west] {\scriptsize $
\gamma(t) $};
\draw [color=qqqqqq](-2.12,3.7) node[anchor=north west] {\scriptsize $
\gamma(\tilde{t})$};
\draw [color=qqqqqq](-1.92,2.84) node[anchor=north west] {\scriptsize $
F_\gamma(t)\cap\Omega$};
\draw [color=qqqqqq] (1.51,5.48)-- (1.54,0.8);
\draw [color=qqqqqq] (1.54,0.8)-- (4.93,0.8);
\draw [color=qqqqqq] (1.51,5.48)-- (4.9,5.48);
\draw [color=qqqqqq] (4.9,5.48)-- (4.93,0.8);
\draw [color=qqqqqq] (1.515,4.716)-- (4.905,4.65);
\draw [rotate around={-89.4:(3.18,2.895)},color=qqqqqq] (3.18,2.895) ellipse
(1.433cm and 1.068cm);
\draw [line width=0.2pt,dash pattern=on 3pt off 3pt,color=qqqqqq]
(1.516,4.502)-- (4.904,4.858);
\draw [line width=0.2pt,dash pattern=on 3pt off 3pt,color=qqqqqq]
(1.513,4.971)-- (4.905,4.65);
\draw [line width=0.2pt,dash pattern=on 3pt off 3pt,color=qqqqqq]
(1.516,4.502)-- (4.902,5.141);
\draw [line width=0.2pt,dash pattern=on 3pt off 3pt,color=qqqqqq]
(1.513,4.971)-- (4.901,5.29);
\draw [->,>=stealth',color=qqqqqq] (4.12,1.26) -- (4.76,1.26);
\draw [->,>=stealth',color=qqqqqq] (4.1,1.26) -- (4.1,1.9);
\draw [->,>=stealth',color=qqqqqq] (2.04,2.78) -- (2.42,3.29);
\draw [color=qqqqqq](3.92,1.34) node[anchor=north west] {\scriptsize $
\mathbf{N}_1(t)$};
\draw [color=qqqqqq](3.18,1.74) node[anchor=north west] {\scriptsize $
\mathbf{N}_2(t) $};
\draw [color=qqqqqq](1.58,1.44) node[anchor=north west] {\scriptsize $
F_\gamma(t) $};
\draw [color=qqqqqq](-0.22,3.42) node[anchor=north west] {\scriptsize $
s_1\mathbf{N}_1(t)+s_2\mathbf{N}_2(t) $};
\draw [color=qqqqqq](2.66,3.04) node[anchor=north west] {\scriptsize $ \gamma(t)
$};
\draw [color=qqqqqq](2.42,2.22) node[anchor=north west] {\scriptsize $
F_\gamma(t)\cap \Omega $};
\draw [color=qqqqqq] (3.18,2.895)-- (4.902,5.141);
\draw [line width=0.2pt,dash pattern=on 3pt off 3pt,color=qqqqqq]
(1.517,4.385)-- (4.907,4.427);
\draw [color=qqqqqq](3.22,3.58) node[anchor=north west] {\scriptsize $
S_s^\gamma(t) $};
\draw [color=qqqqqq](3.9,4.04) node[anchor=north west] {\scriptsize $
L_s^\gamma(t) $};
\draw [color=qqqqqq] (3.944,3.891)-- (4.13,3.75);
\draw [color=qqqqqq] (3.18,2.895)-- (3.74,2.46);
\draw [color=qqqqqq] (3.18,2.895)-- (2.84,3.16);
\draw [color=qqqqqq] (4.349,4.42)-- (4.79,4.07);
\draw [color=qqqqqq] (4.531,4.657)-- (4.18,4.91);
\draw [color=qqqqqq](2.76,4.04) node[anchor=north west] {\scriptsize $
L_s(t,\tilde{t})$};
\draw [->,>=stealth',color=qqqqqq] (3.8,3.52) -- (4.036,3.821);
\draw [->,>=stealth',color=qqqqqq] (3.53,3.1) -- (3.309,2.795);
\draw [->,>=stealth',color=qqqqqq] (4.42,3.94) -- (4.64,4.189);
\draw [->,>=stealth',color=qqqqqq] (4.17,3.57) -- (3.525,2.627);
\draw [->,>=stealth',color=qqqqqq] (3.39,3.52) -- (3.033,3.01);
\draw [->,>=stealth',color=qqqqqq] (3.71,3.94) -- (4.364,4.778);
\draw [->,>=stealth',color=qqqqqq] (0.54,4.52) -- (1.93,4.708);
\draw [shift={(-1.402,5.944)},color=qqqqqq] 
plot[domain=4.03:4.172,variable=\t]({1*2.916*cos(\t r)+0*2.916*sin(\t
r)},{0*2.916*cos(\t r)+1*2.916*sin(\t r)});
\begin{scriptsize}
\fill [color=qqqqqq] (-1.76,3.05) circle (1.0pt);
\fill [color=qqqqqq] (-2.224,3.146) circle (1.0pt);
\fill [color=qqqqqq] (-1.121,3.042) circle (1.0pt);
\fill [color=qqqqqq] (3.18,2.895) circle (1.0pt);
\fill [color=qqqqqq] (4.531,4.657) circle (1.0pt);
\fill [color=qqqqqq] (4.349,4.42) circle (1.0pt);
\end{scriptsize}
\end{tikzpicture}
\caption{}
\label{fig14}
\end{figure}

\begin{lemma}\label{0distance2}
$ \ds L_s^{\gamma}(t)\Big ( \sum_{i=1}^{n-1} s_i\kappa_i(t) \Big )\leq 1$
for all $t\in [0,\ell]$ and $s=(s_1,\cdots,s_{n-1})\in \mathbb{S}^{n-2}$.
\end{lemma}
{\em Proof.}
Suppose $F_\gamma(t)$ and $F_\gamma(\tilde{t})$ are not
parallel. Solving for their intersection yields
\begin{equation*}
\ds \gamma(t)+ \sum_{i=1}^{n-1}  s_i\mathbf{N}_i(t) =\gamma(\tilde{t})
+ \sum_{i=1}^{n-1} r_i\mathbf{N}_i(\tilde{t}).
\end{equation*}
This is a linear system of $n$ equations and $2n-2$ unknowns
$(s_i)_{i=1}^{n-1}$ and $(r_i)_{i=1}^{n-1}$. A solution for this system of equations
exists because the two leading front are not parallel. Then
direct computation using Cramer's rule gives the formula for $F(t,\tilde{t})$
explicitly,
\begin{equation*}
\ds F(t,\tilde{t})=\{x\in F_\gamma(t): (x-\gamma(t)) \cdot
\Big (- \sum_{i=1}^{n-1} \frac{h_i(t,\tilde{t})}{H(t,\tilde{t})}\mathbf{N}_i(t) \Big )=1\},
\end{equation*}
where
$$h_i(t,\tilde{t}):=\det[\mathbf{N}_1(\tilde{t}),\cdots,\mathbf{N}_{n-1}(\tilde{t})
,\mathbf{N}_i(t)] $$
for $1\leq i \leq {n-1}$, 
and
$$H(t,\tilde{t})=\det[\mathbf{N}_1(\tilde{t}),\cdots,\mathbf{N}_{n-1}(\tilde{t}),
\gamma(t)-\gamma(\tilde{t})].$$
Note that $H(t,\tilde{t})\neq 0$ since $\gamma(t)-\gamma(\tilde{t})$ is not
parallel to $F_\gamma(\tilde{t})$.

We firstly claim that
\begin{equation}\label{0distance7}
\ds L_s(t,\tilde{t})\Big ( - \sum_{i=1}^{n-1} \frac{h_i(t,\tilde{t})}{H(t,\tilde{t})}s_i  \Big )\leq 1.
\end{equation}

Indeed, we divide the situation into two cases. In the first case,
suppose we travel from $\gamma(t)$ along a given direction
$\ds \sum_{i=1}^{n-1} s_i\mathbf{N}_i(t)$ and hit $F(t,\tilde{t})$, then for $x\in F(t,\tilde{t})$,
$$x-\gamma(t)=L_s(t,\tilde{t})\Big (\sum_{i=1}^{n-1} s_i\mathbf{N}_i(t) \Big ).$$
Therefore,
\begin{equation}\label{0distance5}
\begin{array}{l}
\ds L_s(t,\tilde{t})\Big ( \sum_{i=1}^{n-1} s_i\mathbf{N}_i(t) \Big )
\cdot
\Big (- \sum_{i=1}^{n-1} \frac{h_i(t,\tilde{t})}{H(t,\tilde{t})}\mathbf{N}_i(t) \Big ) 
\\ \ds =L_s(t,\tilde{t})\Big (-\sum_{i=1}^{n-1} \frac{h_i(t,\tilde{t})}{H(t,\tilde{t})}s_i  \Big )=1.
\end{array}
\end{equation}
Suppose for a certain direction $\ds \sum_{i=1}^{n-1} s_i\mathbf{N}_i(t)$ we do not hit $F(t,\tilde{t})$, in which case we set
$L_{s}(t,\tilde{t})=+\infty$, then we must hit $F(t,\tilde{t})$ through
the direction $\ds -\sum_{i=1}^{n-1} s_i\mathbf{N}_i(t)$, therefore,
by (\ref{0distance5}),
$$
\ds L_{-s}(t,\tilde{t})\Big (\sum_{i=1}^{n-1} \frac{h_i(t,\tilde{t})}{H(t,\tilde{t})} s_i \Big )=1. 
$$
In particular, since $L_{-s}(t,\tilde{t})>0$,
$$
\ds \sum_{i=1}^{n-1} \frac{h_i(t,\tilde{t})}{H(t,\tilde{t})}s_i >0.
$$
We then must have,
\begin{equation}\label{0distance4}
L_{s}(t,\tilde{t})\Big (- \sum_{i=1}^{n-1} \frac{h_i(t,\tilde{t})}{H(t,\tilde{t})}s_i  \Big )<0.
\end{equation}
(\ref{0distance5}) and (\ref{0distance4}) together gives that in either case
(\ref{0distance7})
holds true, which proves our claim.

We secondly claim that,
\begin{equation}\label{0distance3}
\ds L_s^{\gamma}(t) \Big (- \sum_{i=1}^{n-1} \frac{h_i(t,\tilde{t})}{H(t,\tilde{t})}s_i \Big )\leq 1
\end{equation}
for all $t, \tilde{t}\in [0,\ell]$ and $s\in \mathbb{S}^{n-2}$.
Indeed, if for a given $t,\tilde{t}$ and $s\in \mathbb{S}^{n-2}$,
$F_\gamma(t)$ and $F_\gamma(\tilde{t})$ are not parallel, and
\begin{equation}\label{bigger}
 \ds - \sum_{i=1}^{n-1} \frac{h_i(t,\tilde{t})}{H(t,\tilde{t})}s_i \geq 0,
\end{equation}
then,
\begin{equation*}
\ds L_s^{\gamma}(t)\Big (- \sum_{i=1}^{n-1} \frac{h_i(t,\tilde{t})}{H(t,\tilde{t})}s_i \Big ) \leq
L_s(t,\tilde{t})\Big (- \sum_{i=1}^{n-1} \frac{h_i (t,\tilde{t})}{H(t,\tilde{t})}s_i  \Big )=1
\end{equation*}
which gives (\ref{0distance3}) for this case. If for a certain $t,\tilde{t}$ and
$s\in \mathbb{S}^{n-2}$, (\ref{bigger}) fails to hold, then (\ref{0distance3}) is obviously satisfied. Finally, if $F_\gamma(t)$ and
$F_\gamma(\tilde{t})$ are parallel, then $h_i(t,\tilde{t})=0 $ for all $1\leq
i\leq n-1$, hence the (\ref{0distance3}) is again satisfied. 

We thirdly claim that
\begin{equation}\label{0kappa}
-\frac{h_i(t,\tilde{t})}{H(t,\tilde{t})}\rightarrow \kappa_i(t), \quad 1\leq i \leq
{n-1}.
\end{equation}
as $\tilde{t}\rightarrow t$. Indeed, since $\det[\gamma'(t), \mathbf{N}_1(t),\cdots,\mathbf{N}_{n-1}(t)
]=1$ for all $t\in[0,\ell]$, 
$$H(t,\tilde{t})\approx \det[\mathbf{N}_1(\tilde{t}),\cdots,\mathbf{N}_{n-1}(\tilde{t}),
\gamma'(\tilde{t})(t-\tilde{t})] =(-1)^{n-1}(t-\tilde{t})$$
as
$\tilde{t}\rightarrow t$. Moreover,
$$
h_i(t,t)=\det[\mathbf{N}_1(t),\cdots,\mathbf{N}_{n-1}(t),\mathbf{N}_i(t)]=0.
$$
Then,
\begin{multline}\label{0det}
-\frac{h_i(t,\tilde{t})}{H(t,\tilde{t})}\approx-\frac{h_i(t,\tilde{t})-h_i(t,t)}
{(-1)^{n-1}(t-\tilde{t})}\rightarrow \\
(-1)^{n-1}\Bigl(\det[\mathbf{N}'_1(t),\cdots,\mathbf{N}_{n-1}(t),\mathbf{N}_i(t)]+\cdots +
\det[\mathbf{N}_1(t),\cdots,\mathbf{N}'_{n-1}(t),\mathbf{N}_i(t)]\Bigr).
\end{multline}
Recalling (\ref{N}) and pluging this expression into (\ref{0det}) and it is easy to see that all other
terms vanish except
\begin{equation*}
\begin{array}{ll}
\ds \det[\mathbf{N}_1(t),\cdots,\mathbf{N}'_{i}(t),
\cdots,\mathbf{N}_{n-1}(t),\mathbf{N}_i(t)] 
\\ \\  & \ds \hspace{-3cm}  =-\kappa_i \det[\mathbf{N}_1(t),\cdots,\gamma'(t),
\cdots,\mathbf{N}_{n-1}(t),\mathbf{N}_i(t)]\\ 
\\ & \ds  \hspace{-3cm} = \kappa_i \det[\mathbf{N}_1(t),\cdots,\mathbf{N}_{n-1}(t),
\gamma'(t)]=(-1)^{n-1}\kappa_i
\end{array}
\end{equation*}
because $\det[\gamma'(t), \mathbf{N}_1(t),\cdots,\mathbf{N}_{n-1}(t)
]=1$. This proves (\ref{0kappa}).

Passing in (\ref{0distance3}) to the limit $\tilde t \to t$ we obtain the lemma. The proof is
complete.
\hspace*{\fill}$\Box$

Recall that $S_s^\gamma(t)$ as defined in (\ref{0length}) satisfies
$0\leq S_s^\gamma(t)\leq L_s^{\gamma}(t)$ for all $s\in \mathbb{S}^{n-2}$ due to the
fact that $F_\gamma(t)\cap F_\gamma(\tilde{t}) \cap \Omega=\emptyset$ for all
$t, \tilde{t}\in [0,\ell], \tilde{t} \neq t$. We then have,
\begin{corollary}\label{0length2}
$ \ds S_s^{\gamma}(t)\Big ( \sum_{i=1}^{n-1} s_i\kappa_i(t) \Big )\leq 1$
for all $t\in [0,\ell]$ and $s=(s_1,\cdots,s_{n-1})\in \mathbb{S}^{n-2}$.
\end{corollary}
{\em Proof.}
If $\ds \sum_{i=1}^{n-1} s_i\kappa_i(t) \geq 0$, then
$ \displaystyle
S_s^{\gamma}(t)\Big (  \sum_{i=1}^{n-1} s_i\kappa_i(t) \Big )\leq 
L_s^{\gamma}(t)\Big (  \sum_{i=1}^{n-1} s_i\kappa_i(t)  \Big )\leq 1.
$
If $\ds \sum_{i=1}^{n-1} s_i\kappa_i(t) <0$, then the result is
obviously true.  
\hspace*{\fill} $\Box$

From the definition of $\Phi$ in (\ref{0phi}), $\Phi$ is Lipschitz, hence its
Jacobian $J_\Phi=\det D\Phi$ exists a.e. on $\Sigma^\gamma$, where
$\Sigma^\gamma$ has two equivalent representations (\ref{0domain3}) and
(\ref{0domain6}). We will show the Corollary \ref{0length2} implies $J_\Phi>0$
a.e. on $\Sigma^\gamma$, namely,
\begin{lemma}\label{0jacobian}
$ \ds J_\Phi(t,s)=1- \sum_{i=1}^{n-1} s_i\kappa_i(t) >0$ for all $(t,s)\in
\Sigma^\gamma$.
\end{lemma}
{\em Proof.}
Differentiating $\Phi(t,s)$ with respect to $(t,s_1,\cdots,s_{n-1})$ gives,
\begin{equation}\label{0jaco}
J_\Phi(t,s)=\det[\gamma'(t)+ \sum_{i=1}^{n-1} s_i\mathbf{N}'_i(t), \mathbf{N}_1(t),\cdots,\mathbf{N}_{n-1}(t)].
\end{equation}
Substituting (\ref{N}) into (\ref{0jaco}), we
obtain, after Gaussian elimination, that,
\begin{equation}\label{0jaco1}
J_\Phi(t,s)=1- \sum_{i=1}^{n-1} s_i\kappa_i(t).
\end{equation}
If $\ds \sum_{i=1}^{n-1} s_i\kappa_i(t) \leq 0$, then obviously
$J_\Phi(t,s)>0$. Suppose now $\ds \sum_{i=1}^{n-1} s_i\kappa_i(t) > 0$.
By (\ref{0length}) and (\ref{0domain6}) we have
\begin{equation*}
\ds \sum_{i=1}^{n-1} s_i\kappa_i(t) =|s|\Big ( \sum_{i=1}^{n-1} \frac{s_i}{|s|}\kappa_{i}(t) \Big ) 
<S_s^{\gamma}(t)\Big (\sum_{i=1}^{n-1} \frac{s_i}{|s|}\kappa_{i}(t) \Big )\leq 1
\end{equation*}
by Corollary \ref{0length2}. Therefore, $J_\Phi(t,s)>0$ for all $(t,s)\in
\Sigma^\gamma$. The proof is complete.
\hspace*{\fill} $\Box$

\subsection{Moving Frames in the target space.}\label{smoving}
We are now in a position to define the moving frame in the target space
$\bbbr^{n+1}$. Let $\mathbf{N}_i(t), 1\leq i \leq {n-1}$ be as in subsection
\ref{sleading}. Define the leading curve corresponding to $\gamma$ in
$u(\Omega(\gamma))$ to be
$$\tilde{\gamma}:=u\circ \gamma. $$
We also recall from subsection \ref{sfoliation} the definitions (\ref{0phi}), (\ref{0leadingplane2}), and that $\nabla u$ is
constant on $P_\gamma (t)$ for each $t\in [0,\ell]$. Hence for each $t\in
[0,\ell]$, $\nabla u \circ \Phi$ is constant on $\Sigma^\gamma(t)$.

Consider the Darboux frame $(\tilde{\gamma}', \mathbf{v}_1,\cdots,\mathbf{v}_{n-1},
\mathbf{n})$ where $\mathbf{v}_i(t) =\nabla u(\gamma(t))\mathbf{N}_i(t)$, $i=1,...,n-1$ and $\mathbf{n}(t)= \tilde{\gamma}'(t)\times \mathbf{v}_1(t) \times \cdots \times
\mathbf{v}_{n-1}(t)$. Since $u$ is an isometric affine map along $P_\gamma (t)$ for each $t\in
[0,\ell]$ we obtain
\begin{equation}\label{0affine}
\ds u(\Phi(t,s))=\tilde{\gamma}(t)+ \sum_{i=1}^{n-1} s_i\mathbf{v}_i(t)
\end{equation}
for all $t\in [0,\ell]$ and $s\in \Sigma^\gamma(t)$.
Differentiating with respect to $t$, by (\ref{0phi}) we get
\begin{equation}\label{0gamma}
  \nabla u(\Phi (t,s))
\Big (\gamma'(t)+\sum_{i=1}^{n-1} s_i\mathbf{N}_i'(t) \Big) 
  = \tilde{\gamma}'(t)+ \sum_{i=1}^{n-1} s_i\mathbf{v}_i'(t), 
\end{equation}
and differentiating with respect to $s_i, 1\leq i\leq n-1$ we obtain for each $i$,
\begin{equation}\label{0v}
\nabla u(\Phi (t,s))\mathbf{N}_i(t)=\mathbf{v}_i(t).
\end{equation}
By the linear expansion of $N_i'$ in (\ref{N}) and (\ref{R^TR'}),
together with (\ref{0gamma}) and (\ref{0v}) we get
\begin{equation}\label{0iso}
\begin{array}{l}
\ds \tilde{\gamma}'(t)+\sum_{i=1}^{n-1} s_i\mathbf{v}_i'(t)  \\ = \ds \nabla u(\Phi (t,s)) \Big 
(1- \sum_{i=1}^{n-1} s_i \kappa_i (t) \Big )\gamma'(t) 
 +\sum_{i=1}^{n-1}s_i \Bigl(\sum_{j=1}^{n-1}\kappa_{i_j}(t)\mathbf{v}_j(t)\Bigr)
\end{array}
\end{equation} 
with $\kappa_{i_i}=0$ and $\kappa_{i_j}=-\kappa_{j_i}$. Also,  by (\ref{0gamma}), for $s=0$ we have 
$$\nabla u(\Phi (t,0))\gamma'(t)=\tilde{\gamma}'(t). $$
Since $\nabla u \circ \Phi$ is constant on $\Sigma^\gamma(t)$ for each $t\in
[0,\ell]$, we obtain
\begin{equation}\label{0t}
\nabla u(\Phi (t,s))\gamma'(t)=\nabla u(\Phi (t,0))\gamma'(t)
=\tilde{\gamma}'(t) \quad \textrm{for all } s\in \Sigma^\gamma(t).
\end{equation} 
Alongside (\ref{0v}), this shows that at each point in $\Omega(\gamma)$, $\nabla u$
maps an orthonormal frame to another orthonormal frame and this orthonormal
frame only depends on $t$. Finally, using (\ref{0iso}) and matching coefficients yields for all $1\leq i\leq n-1$,
\begin{equation}\label{vi}
\mathbf{v}_i'=-\kappa_i \tilde{\gamma}'+\sum_{j=1}^{n-1}\kappa_{i_{j}}\mathbf{v}_{j}, \quad  \kappa_{i_i}=0 \textrm{  and  } \kappa_{i_j}=-\kappa_{j_i} .
\end{equation} 
In other words, the following system of ODEs is satisfied by the Darboux frame of $\tilde \gamma$:
\begin{equation}\label{R^TR'2}
 \hspace{1cm}   \left ( \begin{array}{c}  \tilde{\gamma}' \\ \mathbf{v}_1  \\ \mathbf{v}_2\\
\vdots  \\ \mathbf{v}_{n-1}\\ \mathbf{n} \end{array} \right )'
= \mathcal{K}
\left ( \begin{array}{c}  \tilde{\gamma}' \\ \mathbf{v}_1  \\ \mathbf{v}_2\\
\vdots  \\ \mathbf{v}_{n-1}\\ \mathbf{n} \end{array} \right ),
\end{equation} where the skew-symmetric curvature matrix $\mathcal{K}$ is given by 
 \begin{equation*}
\mathcal{K}= \left ( \begin{array}{cccccc}
 0 & \kappa_1 & \kappa_2 & \cdots & \kappa_{n-1} &\kappa_{\mathbf{n}} \\
 -\kappa_1 & 0 & \kappa_{1_2} & \cdots & \kappa_{1_{n-1}} & 0  \\
 -\kappa_2 & -\kappa_{1_2} & 0 & \cdots & \kappa_{2_{n-1}} & 0 \\
 \vdots & \vdots & \vdots & \cdots &\vdots & \vdots \\
 -\kappa_{n-1} & -\kappa_{1_{n-1}} & -\kappa_{2_{n-1}} & \cdots & 0 & 0\\
 -\kappa_{\mathbf{n}} & 0 & 0 & \cdots & 0 & 0 
\end{array} \right ). 
\end{equation*} 

\subsection{Change of variable formula.}\label{schangeofvariable}
Recall that $\Phi:\Sigma^\gamma\rightarrow \Omega(\gamma)$ is one-to-one and
onto, where $\Sigma^\gamma$ was defined in (\ref{0domain3}), For $(t,s)\in
\Sigma^\gamma$, let $u_i(t,s):=(\frac{\partial }{\partial x_i}u)\circ\Phi
(t,s)$, note that $u_i$ is the $i$th column of $\nabla u \circ \Phi$. The
following holds for all $(t,s)\in \Sigma^\gamma$: 
since $\nabla u^T\mathbf{n} \cdot \gamma'=\mathbf{n}\cdot \nabla
u\gamma'=\mathbf{n}\cdot \tilde{\gamma}' =0$ and $\nabla u^T\mathbf{n} \cdot
\mathbf{N}_j=\mathbf{n}\cdot \nabla u\mathbf{N}_j=\mathbf{n}\cdot
\mathbf{v}_j=0$ for all $1\leq j \leq n-1$, we have $\nabla u^T\mathbf{n}=0$,
i.e. $u_i\cdot \mathbf{n}=0$ for all $1\leq i \leq n$. Thus,
\begin{equation}\label{0ui}
\begin{array}{ll}
\ds 
u_i=\bigl(u_i\cdot \tilde{\gamma}'\bigr)\tilde{\gamma}'
+\sum_{j}\bigl(u_i\cdot \mathbf{v}_j\bigr)\mathbf{v}_j
+\bigl(u_i\cdot \mathbf{n}\bigr)\mathbf{n} & \\ & \ds \hspace{-4.5cm}  
= (u_i\cdot \tilde{\gamma}')\tilde{\gamma}'+\sum_{j}(u_i\cdot
\mathbf{v}_j)\mathbf{v}_j \\ & \ds \hspace{-4.5cm} 
   = (u_i\cdot \nabla u \gamma' )\tilde{\gamma}'+\sum_{j}(u_i\cdot
\nabla u\mathbf{N}_j )\mathbf{v}_j \\ & \ds \hspace{-4.5cm} 
   = (\nabla u^T u_i \cdot \gamma')\tilde{\gamma}'+\sum_{j}(\nabla u^T
u_i \cdot \mathbf{N}_j )\mathbf{v}_j \\ & \ds \hspace{-4.5cm} 
   = (\mathbf{e}_i\cdot
\gamma')\tilde{\gamma}'+\sum_{j}(\mathbf{e}_i\cdot \mathbf{N}_j
)\mathbf{v}_j,
\end{array}
\end{equation}
where $\mathbf{e}_i=(0,..,1,...,0)$. Note that the right hand side of
(\ref{0ui}) is independent of $s$. Differentiating with respect to $s_j, 0\leq j \leq n-1$, by (\ref{0phi})
we get for all $1\leq i \leq n$ and $1\leq j \leq n-1$,
\begin{equation}\label{0s}
(\nabla \frac{\partial }{\partial x_i}u)(\Phi(t,s)) \mathbf{N}_j(t)=0 .
\end{equation}
Differentiating $u_i$ with respect to $t$ we obtain,
\begin{equation}\label{0w}
\begin{array}{rl}
\ds  (\nabla \frac{\partial }{\partial x_i}u)(\Phi(t,s))
(\gamma'(t)+ \sum_{j=1}^{n-1} s_j\mathbf{N}_j'(t)) & \\ &
\hspace{-2.5cm} \ds 
  =(\mathbf{e}_i\cdot
\gamma^{\prime\prime}(t))\tilde{\gamma}'(t)+(\mathbf{e}_i\cdot
\gamma'(t))\tilde{\gamma}''(t) \\ & \hspace{-2.5cm} \ds +\sum_{j=1}^{n-1}(\mathbf{e}_i \cdot
\mathbf{N}_j'(t) )\mathbf{v}_j (t) +\sum_{j=1}^{n-1}(\mathbf{e}_i \cdot
\mathbf{N}_j(t) )\mathbf{v}_j'(t).
\end{array}
\end{equation}
If we write out $\mathbf{N}_i'$ as a linear combination of $\gamma'$ and
$\mathbf{N}_j,j=1,\cdots n-1$ as in (\ref{N}) and (\ref{R^TR'}), the left hand side of (\ref{0w}) becomes
$$ \ds (1- \sum_{j=1}^{N-1}  s_j\kappa_j(t)) (\nabla \frac{\partial
}{\partial x_i}u)(\Phi(t,s))\gamma'(t). $$
For the right hand side of (\ref{0w}), if we write out $\gamma''$, $\tilde{\gamma}''$, $\mathbf{N}_j'$ and $\mathbf{v}_j'$ as linear
combinations of $\gamma'$, $\tilde{\gamma}'$, $\mathbf{N}_\ell$ and $\mathbf{v}_\ell, \ell=1,\cdots,n-1$ and
$\mathbf{n}$ as in (\ref{N}) and (\ref{vi}), we obtain, 
\begin{equation*}
\begin{array}{l}
\ds (\mathbf{e}_i\cdot
\gamma^{\prime\prime})\tilde{\gamma}'+(\mathbf{e}_i\cdot
\gamma')\tilde{\gamma}''  +\sum_{j}(\mathbf{e}_i \cdot
\mathbf{N}_j' )\mathbf{v}_j  +\sum_{j}(\mathbf{e}_i \cdot
\mathbf{N}_j )\mathbf{v}_j' 
 \\ = \ds (\mathbf{e}_i\cdot \sum_{j}\kappa_j\mathbf{N}_j)\tilde{\gamma}'  +(\mathbf{e}_i\cdot \gamma')\bigl(\sum_{j}\kappa_j\mathbf{v}_j
+\kappa_n\mathbf{n}\bigr) 
+\sum_{j}\bigl(\mathbf{e}_i\cdot (-\kappa_j\gamma'+\sum_{\ell}\kappa_{j_\ell}\mathbf{N}_\ell)\bigr)\mathbf{v}_j
\\  \ds +\sum_{j}(\mathbf{e}_i\cdot \mathbf{N}_j)\bigl(-\kappa_j\tilde{\gamma}'+\sum_{\ell }\kappa_{j_\ell}\mathbf{v}_\ell\bigr)
 \\  =  \ds \sum_{j}(\mathbf{e}_i\cdot \kappa_j\mathbf{N}_j)\tilde{\gamma}'+(\mathbf{e}_i\cdot \gamma')\sum_{j}\kappa_j\mathbf{v}_j
+(\mathbf{e}_i\cdot \gamma')\kappa_n\mathbf{n}-(\mathbf{e}_i\cdot \gamma')\sum_{j}\kappa_j\mathbf{v}_j \\ \ds 
+\sum_{j}\sum_{\ell}(\mathbf{e}_i\cdot \kappa_{j_\ell}\mathbf{N}_\ell)\mathbf{v}_j
-\sum_{j}(\mathbf{e}_i\cdot \kappa_j\mathbf{N}_j)\tilde{\gamma}'+\sum_{\ell}\sum_{j}(\mathbf{e}_i\cdot\mathbf{N}_\ell)\kappa_{\ell_j}\mathbf{v}_j
 = \ds (\mathbf{e}_i\cdot \gamma')\kappa_n\mathbf{n}
\end{array}
\end{equation*}
where we used  the fact that $\kappa_{i_j}=-\kappa_{j_i}$. By
Lemma \ref{0jacobian}, $\ds 1- \sum_{j=1}^{n-1} s_j\kappa_j(t) >0$ for
all $(t,s)\in\Sigma^\gamma$. Therefore,
\begin{equation}\label{0second}
\ds (\nabla \frac{\partial }{\partial x_i}u)(\Phi(t,s))\gamma'(t)=\frac{(\mathbf{e}_i
\cdot \gamma'(t))\kappa_{\mathbf{n}}(t)\mathbf{n}(t)}{\ds 1- \sum_{j=1}^{n-1} s_j\kappa_j(t)}.
\end{equation}
Since $\Phi$ is Lipschitz with
$ \ds J_\Phi(t,s)=1- \sum_{j=1}^{n-1} s_j\kappa_j(t) >0 $,
the change of variable $x=\Phi(t,s)$ with (\ref{0affine}) and (\ref{0second}) yields,
\begin{equation}\label{0changeofvariable}
\begin{array}{l}
 \ds \int_{\Omega(\gamma)}|u(x)|^2dx \\  \ds
=
 \int_0^\ell\int_{\Sigma^\gamma(t)}|\tilde{\gamma}(t)+ \sum_{i=1}^{n-1} s_i\mathbf{v}
_i(t)|^2
 \cdot
\bigl(1- \sum_{j=1}^{n-1} s_j\kappa_j(t) \bigr)d\mathcal{H}^{n-1}
(s)dt,
\end{array}
\end{equation}
\begin{equation}\label{0changeofvariable01}
\int_{\Omega(\gamma)}|\nabla u(x)|^2dx =n |\Omega(\gamma)|,
\end{equation}
\begin{equation}\label{0changeofvariable02}
\begin{array}{ll}
\ds \int_{\Omega(\gamma)}|\nabla^2 u(x)|^2dx &  = \ds   
\int_0^\ell\int_{\Sigma^\gamma(t)}\frac{\ds \sum_i(\mathbf{e}_i\cdot
\gamma'(t))^2\kappa_{\mathbf{n}}^2(t)}
{\ds \bigl(1- \sum_{j=1}^{n-1} s_j\kappa_j(t) \bigr)}d\mathcal{H}^{n-1}
(s)dt 
\\  & \ds =\int_0^\ell\int_{\Sigma^\gamma(t)}\frac{\kappa_{\mathbf{n}}^2(t)}
{\ds \bigl(1- \sum_{j=1}^{n-1} s_j\kappa_j(t) \bigr)}d\mathcal{H}^{n-1}
(s)dt.
\end{array}
\end{equation}

\subsection{Approximation process for
$u|_{\Omega(\gamma)}$.}\label{sapproximation}
Recall $L_s^\gamma(t)$ and $S_s^\gamma(t)$ defined in (\ref{0distance}) and (\ref{0length}) respectively. Since all leading fronts meet outside $\Omega$,
we must have $L_s^\gamma(t)\geq S_s^\gamma(t)$ for all $s\in \mathbb{S}^{n-2}$ and $t\in
[0,\ell]$.

\begin{lemma}\label{0dilation}
There exists a sequence of isometries $u_m \in
W^{2,2}(\Omega(\gamma),\bbbr^{n+1})$ converging strongly to $u$ with the
property that each $u_m$ has a suitable leading curve
$\gamma_m:[0,\ell_m]\rightarrow \bbbr^n$ for which $L_s^{\gamma_m}(t)-
S_s^{\gamma_m}(t)>\rho_m>0$ for all $s\in \mathbb{S}^{n-2}$ and $t\in
[0,\ell_m]$.
\end{lemma}
{\em Proof.}
The proof is exactly the same as the 2-dimensional case,
\cite[Prop.3.2]{Pak}. For this reason, it is omitted. \hspace*{\fill} $\Box$ 
 
\begin{remark}\label{0rho}
By the above Lemma, we can just assume $u$ has a suitable leading curve $\gamma$
that satisfies $L_s^{\gamma}(t)- S_s^{\gamma}(t)>\rho>0$ for all $s\in
\mathbb{S}^{n-2}$ and $t\in [0,\ell]$. 
\end{remark}

\begin{lemma}\label{0construction}
Suppose $L_s^{\gamma}(t)- S_s^{\gamma}(t)>\rho>0$ for all $s\in
\mathbb{S}^{n-2}$ and $t\in [0,\ell]$. Then there is a sequence of
smooth maps in $I^{2,2}(\Omega(\gamma),\bbbr^{n+1})$ converging strongly to $u$.
\end{lemma}
{\em Proof.}
The idea is to construct a smooth curve $\gamma_m$ approximating
$\gamma$. We do not know yet this curve is a leading curve of $u_m$ or not, so
we cannot call the $(n-2)$-dimensional hyperplane orthogonal to $\gamma_m$ at
$t$ leading fronts. Instead we call them \textit{orthogonal fronts} and denote
them by $F_{\gamma_m}(t)$. If we manage to show all such orthogonal fronts meet
outside $\Omega(\gamma_m)$, $\gamma_m$ becomes a leading curve for $u_m$ and
$F_{\gamma_m}(t)$ are actually the leading fronts. We then define $u_m$ to be
isometric affine mapping along each leading front $F_{\gamma_m}(t)$. Since all
the leading fronts intersect outside $\Omega$, $u_m$ is well-defined.

We first need the following lemma,
\begin{lemma}\label{0difficult}
There exists smooth curve $\gamma_m$ such that $\gamma_m(t)\rightarrow
\gamma(t)$ strongly in $W^{2,p}([0,\ell],\bbbr^n)$ for all $1\leq p<\infty$ and satisfies
$F_{\gamma_m}(t)\cap F_{\gamma_m}(\tilde{t})\cap \overline{\Omega}=\emptyset$
for all $t,\tilde{t}\in [0,\ell]$.
\end{lemma}

{\em Proof.}
The construction is long and technical so we postpone the proof to Appendix \ref{AA}. 
\hspace*{\fill} $\Box$

We also need to define the curves $\tilde{\gamma}_m$ in the target space
$u(\Omega(\gamma))$
corresponding to $\gamma_m$. Recall that the normal curvature $\kappa_{\mathbf{n}}$ defined in  
(\ref{R^TR'2})  is bounded. We choose a sequence of uniformly bounded
smooth function $\tilde{\kappa}_{\mathbf{n},m}$ such that
$\tilde{\kappa}_{\mathbf{n},m}\rightarrow \kappa_{\mathbf{n}}$ a.e. in
$[0,\ell]$, (and hence in $L^p$ for all $1\leq p<\infty$).

We need to flatten $\tilde{\kappa}_{\mathbf{n},m}$ around the end points $0$ and
$\ell$ for two reasons: first, it might happen
that $\Omega(\gamma)\nsubseteq
\Omega(\gamma_m)$ so we need to extend the isometric immersion defined on
$\Omega(\gamma_m)$ smoothly to the region of $\Omega(\gamma)$ outside
$\Omega(\gamma_m)$. Second, so far all the construction is on one covered
domain $\Omega(\gamma)$ and our final goal is to glue all the different covered
domains together smoothly. By flattening
$\tilde{\kappa}_{\mathbf{n},m}$ around the end point $0$ and $\ell$,
$u_m$ constructed later is affine near the leading planes $P_\gamma(0)$ and
$P_\gamma(\ell)$ (for definition of leading planes see Definition
\ref{0leadingplane}) so that we can join all the pieces smoothly. The
modification goes as follows: by (\ref{0second}), the second derivative of $u$
vanishes
whenever
$\kappa_{\mathbf{n}}=0$. Put
\begin{eqnarray*}
\ell_m^*= \biggl\{ \begin{array} {l}
\ell \quad \textrm { if} \quad \Omega(\gamma)\subset
\Omega(\gamma_m) \textrm{  and},\\
\sup\{t\in [0,\ell], F_{\gamma_m}(t) \cap F_\gamma(\ell) \cap
\overline{\Omega}(\gamma) =\emptyset\} \textrm {  otherwise}.
\end{array}
\end{eqnarray*}
By step 1 of Lemma \ref{0difficult} in the Appendix, $F_{\gamma_m}(t)\rightarrow F_\gamma(t)$
uniformly, hence $\ell_m^*\rightarrow \ell$ as $m\rightarrow \infty$.

Let $\psi_1$ be any smooth non-negative function which is $0$ on $[-1,\infty)$ and
$1$ on $(-\infty, -2)$. Let $\psi_2$ be any smooth positive function which is
$0$ on $(-\infty,1]$ and $1$ on $(2,\infty)$. We put,
$$
\kappa_{\mathbf{n},m}(t):=\psi_1(m(t-\ell_m^*))\psi_2(mt)\tilde{\kappa}_{\mathbf
{n},m}(t), \quad t\in [0,\ell]$$
and we solve the following linear system for initial values
$\tilde{\gamma}'_m(0)=\tilde{\gamma}'(0)$,
$\mathbf{v}_{i,m}(0)=\mathbf{v}_{i}(0)$, and $\mathbf{n}_m(0)=\mathbf{n}(0)$:
\begin{align*}
   \left ( \begin{array}{c}  \tilde{\gamma}'_m \\ \mathbf{v}_{1,m}  \\
\mathbf{v}_{2,m}\\ \vdots  \\ \mathbf{v}_{n-1,m}\\ \mathbf{n}_m \end{array}
\right )'
=  \mathcal{K}_m
\left ( \begin{array}{c}  \tilde{\gamma}'_m \\ \mathbf{v}_{1,m}  \\
\mathbf{v}_{2,m}\\ \vdots  \\ \mathbf{v}_{n-1,m}\\ \mathbf{n}_m  \end{array}
\right ), 
\end{align*} where the matrix $\mathcal {K}_m$ is given by 
\begin{align*} 
\mathcal{K}_m=\begin{pmatrix}
 0 & \kappa_{1,m} & \kappa_{2,m} & \cdots & \kappa_{n-1,m} & \kappa_{\mathbf{n},m}
\\
 -\kappa_{1,m} & 0 & \kappa_{1_2,m} & \cdots & \kappa_{1_{n-1},m} & 0\\
 -\kappa_{2,m} & -\kappa_{1_2,m} & 0 & \cdots & \kappa_{2_{n-1},m} & 0\\
  \vdots & \vdots & \vdots & \cdots  & \vdots & \vdots \\
 -\kappa_{n-1,m} & -\kappa_{1_{n-1},m} & -\kappa_{2_{n-1},m} & \cdots & 0 & 0\\
 -\kappa_{\mathbf{n},m} & 0 & 0 & \cdots & 0 & 0\\
 \end{pmatrix}.
\end{align*}

We define
$$\tilde{\gamma}_m(t)=\tilde{\gamma}(0)+\int_0^t \tilde{\gamma}'_m(\tau)d\tau.$$
By the same argument as in step 1 in the proof of Lemma \ref{0difficult},
$\tilde{\gamma}_m\rightarrow \tilde{\gamma}$ in $W^{2,p}([0,\ell],\bbbr^{n+1})$
and the moving frame $(\tilde{\gamma}'_m,
\mathbf{v}_{1,m},\cdots,\mathbf{v}_{n-1,m}, \mathbf{n}_m)$ converges
to $(\tilde{\gamma}', \mathbf{v}_1,\cdots,\mathbf{v}_{n-1}, \mathbf{n})$ uniformly.

Eventually, we define our approximating sequence $u_m$ on $\Omega(\gamma_m)$:
\begin{equation}\label{0affine1}
\ds u_m(\gamma_m(t)+ \sum_{i=1}^{n-1} s_i\mathbf{N}_{i,m}(t))\\
=\ds \tilde{\gamma}_m(t)+ \sum_{i=1}^{n-1} s_i \mathbf{v}_{i,m}(t),  
\end{equation}
where $\gamma_m$ is defined in Lemma \ref{0difficult}. Such $\gamma_m$ assures
that all its leading fronts intersect outside $\overline{\Omega}$, hence
$u_m$ is well-defined and smooth over $\Omega(\gamma)\cap \Omega(\gamma_m)$.

As before, let $\Phi_m:[0,\ell]\times \bbbr^{n-1}\rightarrow \bbbr^n$ be defined
as
$$
\ds \Phi_m(t,s)=\gamma_m(t)+ \sum_{i=1}^{n-1} s_i\mathbf{N}_{i,m}(t),
$$
and let $\Delta^{\gamma_m}=\{(t,s):\Phi_m(t,s)\in \Omega(\gamma)\}$. Same
argument as
Step 6 in Lemma \ref{0difficult} gives that $\Phi_m(t,s)$ is a bi-Lipschitz mapping
of $\Delta^{\gamma_m}$ onto $\Omega(\gamma)\cap \Omega(\gamma_m)$. By
differentiating with respect to
$t,s_1,\cdots, s_{n-1}$, as in (\ref{0t}) and (\ref{0v}), we see that at each point
of $x$, $\nabla u_m(x)$ maps an orthonormal frame to an orthonormal frame. Hence
$\nabla u_m(x)^T \nabla u_m(x) ={\rm I}$. Moreover, $u_m$ is affine near
$P_{\gamma_m}(\ell)$ and can be extended by
an affine isometry over $\Omega(\gamma)$. Therefore, $u_m\in
I^{2,2}(\Omega(\gamma), \bbbr^n)$.
Everything we have proved for isometric immersions of course applies, in
particular, by (\ref{0ui}), (\ref{0second}), and (\ref{0s}) we have for all $1\leq i\leq n-1$,
\begin{equation}\label{0ui1}
\frac{\partial}{\partial x_i}u_m\circ\Phi_m (t,s)=
   (\mathbf{e}_i\cdot
\gamma'_m(t))\tilde{\gamma}_m'(t)+\sum_{j=1}^{n-1}(\mathbf{e}_i\cdot
\mathbf{N}_{j,m}(t) )\mathbf{v}_{j,m}(t), 
\end{equation}
\begin{equation}\label{0second1}
(\nabla \frac{\partial }{\partial
x_i}u_m)(\Phi_m(t,s))\gamma_m'(t)=\frac{(\mathbf{e}_i\cdot
\gamma'_m(t))\kappa_{\mathbf{n},m}(t)\mathbf{n}(t)}{\ds 1
- \sum_{i=1}^{n-1} s_i\kappa_{i,m} (t)}, \quad \textrm{and}
\end{equation}
\begin{equation}\label{0s1}
(\nabla \frac{\partial }{\partial x_i}u_m)(\Phi_m(t,s))\mathbf{N}_{j,m}(t)=0,
\quad \textrm{for all  }\, 1\leq j\leq n-1.
\end{equation}
for all $t\in [0,\ell]$ and $s=(s_1,\cdots,s_{n-1})\in
\Delta^{\gamma_m}(t)$. Moreover, by (\ref{0changeofvariable}), (\ref{0changeofvariable01}), and (\ref{0changeofvariable02}) we compute,
\begin{equation}\label{0changeofvariable1}
\begin{array}{l}
\ds \int_{\Omega(\gamma)}|u_m(x)|^2dx \\ \ds =\int_{\Omega(\gamma)\cap
\Omega(\gamma_m)}|u_m(x)|^2dx+\int_{\Omega(\gamma)\setminus
\Omega(\gamma_m)}|u_m(x)|^2dx \\  \ds 
=\ds \int_0^\ell\int_{\Delta^{\gamma_m}(t)}|\tilde{\gamma}_m(t)+ \sum_{i=1}^{n-1} s_i \mathbf{v}_{i,m}
(t)|^2 
 \cdot
\Big (1- \sum_{i=1}^{n-1} s_j\kappa_{j,m}(t) \Big )d\mathcal{H}^{
n-1}(s)dt\\  \ds 
+\int_{\Omega(\gamma)\setminus \Omega(\gamma_m)}|u_m(\ell)+\nabla
u_m(\ell)(x-\gamma_m(\ell))|^2dx,
\end{array}
\end{equation}
\begin{equation}
\int_{\Omega(\gamma)}|\nabla u_m(x)|^2dx =n
|\Omega(\gamma)|,
\end{equation}
\begin{equation}
\begin{array}{l}
\ds \int_{\Omega(\gamma)}|\nabla^2 u_m(x)|^2dx  \\ \ds =  \int_{\Omega(\gamma)\cap
\Omega(\gamma_m)}|\nabla^2 u_m(x)|^2dx  \ds +
\int_{\Omega(\gamma)\setminus
\Omega(\gamma_m)}|\nabla^2 u_m(x)|^2dx\\  \ds 
\ds =\int_0^\ell\int_{\Delta^{\gamma_m}(t)}\frac{\kappa_{\mathbf{n},m}^2(t)}
{\Big (1- \sum_{i=1}^{n-1} s_i \kappa_{i,m}(t) \Big )}d\mathcal{H}^{
n-1}(s)dt+0.
\end{array}
\end{equation}

It is easy to see that $u_m\rightarrow u$ in $W^{2,2}(\Omega(\gamma),\bbbr^{n+1})$
because $(\gamma_m',
\mathbf{N}_{1,m},\cdots\mathbf{N}_{n-1,m})$ converges to $ (\gamma',
\mathbf{N}_{1},\cdots\mathbf{N}_{n-1})$ uniformly, $(\tilde{\gamma}'_m,
\mathbf{v}_{1,m},\cdots,\mathbf{v}_{n-1,m},
\mathbf{n}_m)$ converges to $(\tilde{\gamma}', \mathbf{v}_1,\cdots,\mathbf{v}_{n-1},
\mathbf{n})$ uniformly, $\kappa_{\mathbf{n},m}\rightarrow
\kappa_{\mathbf{n}}$, $\kappa_{i,m}\rightarrow
\kappa_{i}, 1\leq i\leq n-1$ in $L^p([0,\ell])$ for all $1\leq p<\infty$,
$\ds 1-\sum_{i=1}^{n-1} s_i\kappa_{i,m}(t) \geq
\min\{\rho/16d,1/2\}$, $\Delta^{\gamma_m}(t)\rightarrow \Sigma^\gamma(t)$ for
all $t\in [0,\ell]$
and $|\Omega(\gamma)\setminus \Omega(\gamma_m)|\rightarrow 0$.
The proof is complete.
\hspace*{\fill} $\Box$

Combining Lemmas \ref{0dilation} and \ref{0construction} we get a smooth
approximation sequence for any isometry $u$ in $\Omega(\gamma)$.

\subsection{Approximation for $u$ in $\Omega$.}\label{sappro}
The proof is exactly the same as
the proof in section 3.3 in \cite{Pak}. Since it is the final part of the argument, we briefly review it 
for the convenience of the reader. 

Recall that we defined a maximal region on which $u$ is affine a \textit{body}
if its
boundary contains more than \textit{two} different $(n-1)$-planes in $\Omega$
(recall Definition \ref{0kplane} for the definition of $(n-1)$-planes in
$\Omega$) and we have shown that we can assume $\Omega$ has only a finite number
of
bodies and is partitioned into bodies and covered domains. We call the maximal
subdomain covered by some leading curve $\gamma$ an \textit{arm}. Similarly to
Lemma \ref{0body} we also have,
\begin{lemma}
It is sufficient to prove Theorem \ref{0density} for a function in
$I^{2,2}(\Omega, \bbbr^{n+1})$ with finite number of arms.
\end{lemma}
{\em Proof.} The proof is the same as the two dimensional case in
\cite[Lemma 3.9]{Pak} and is omitted for brevity. \hspace*{\fill} $\Box$

\begin{figure}[ht]
\centering
\begin{tikzpicture}[line cap=round,line join=round,>=triangle
45,x=1.0cm,y=1.0cm]
\clip(-4.28,1.34) rectangle (5.96,6.3);
\fill[line width=0pt,color=qqffff,fill=qqffff,fill opacity=0.5] (-1.844,5.839) -- (-2.571,1.937) -- (-1.213,1.676) -- (0.311,1.982) -- cycle;
\fill[line width=0pt,color=qqffff,fill=qqffff,fill opacity=0.5] (-1.038,5.878) -- (0.641,2.154) -- (1.638,3.472) -- cycle;
\fill[line width=0pt,color=qqffff,fill=qqffff,fill opacity=0.5] (-3.325,5.231) -- (-3.178,5.337) -- (-3.116,5.378) -- (-3.062,5.411) -- (-2.999,5.448) -- (-3.162,2.249) -- (-3.298,2.347) -- (-3.399,2.43) -- (-3.503,2.525) -- (-3.59,2.613) -- (-3.683,2.72) -- (-3.764,2.827) -- (-3.831,2.926) -- (-3.894,3.035) -- (-3.959,3.17) -- (-4.009,3.299) -- cycle;
\fill[line width=0pt,color=qqffff,fill=qqffff,fill opacity=0.5] (-0.508,5.813) -- (-0.34,5.776) -- (-0.158,5.727) -- (-0.006,5.678) -- (0.124,5.63) -- (0.315,5.548) -- (0.48,5.464) -- (0.592,5.401) -- (0.689,5.34) -- (0.782,5.276) -- (0.878,5.204) -- (0.988,5.113) -- (1.071,5.037) -- (1.173,4.932) -- (1.296,4.785) -- (1.388,4.654) -- (1.457,4.538) -- (1.528,4.393) -- (1.582,4.252) -- (1.637,4.04) -- cycle;
\draw [rotate around={-1.166:(-1.215,3.78)},color=qqqqqq] (-1.215,3.78) ellipse (2.879cm and 2.104cm);
\draw [color=qqttzz] (-3.529,5.053)-- (-4.089,3.687);
\draw [color=qqttzz] (-2.774,5.563)-- (-3.039,2.169);
\draw [color=qqttzz] (-2.461,5.688)-- (-2.875,2.077);
\draw [color=qqttzz] (-1.422,5.88)-- (0.508,2.078);
\draw [dash pattern=on 3pt off 3pt,color=qqqqff] (-0.158,5.727)-- (1.528,4.393);
\draw [color=qqttzz] (0.867,2.308)-- (1.594,3.292);
\draw [color=qqttzz] (-0.981,1.681)-- (0.149,1.914);
\draw [color=qqttzz] (-2.723,2.002)-- (-2.153,5.778);
\draw [color=qqttzz] (-3.325,5.231)-- (-4.009,3.299);
\draw [color=qqttzz] (-2.999,5.448)-- (-3.162,2.249);
\draw [color=qqttzz] (-1.844,5.839)-- (-2.571,1.937);
\draw [color=qqttzz] (-1.844,5.839)-- (0.311,1.982);
\draw [color=qqttzz] (-2.571,1.937)-- (-1.213,1.676);
\draw [color=qqttzz] (-1.213,1.676)-- (0.311,1.982);
\draw [color=qqttzz] (-1.038,5.878)-- (0.641,2.154);
\draw [color=qqttzz] (-1.038,5.878)-- (1.638,3.472);
\draw [color=qqttzz] (-2.333,1.852)-- (-1.427,1.684);
\draw [color=qqttzz] (-0.508,5.813)-- (1.637,4.04);
\draw [color=qqttzz] (-0.794,5.857)-- (1.664,3.74);
\draw [color=qqttzz] (1.638,3.472)-- (0.641,2.154);
\draw [dash pattern=on 3pt off 3pt,color=qqqqff] (0.315,5.548)-- (1.296,4.785);
\draw [dash pattern=on 3pt off 3pt,color=qqqqff] (-3.178,5.337)-- (-3.764,2.827);
\draw [dash pattern=on 3pt off 3pt,color=qqqqff] (-3.116,5.378)-- (-3.59,2.613);
\draw [dash pattern=on 3pt off 3pt,color=qqqqff] (-3.062,5.411)-- (-3.399,2.43);
\draw [dash pattern=on 3pt off 3pt,color=qqqqff] (-3.241,5.294)-- (-3.894,3.035);
\draw [color=qqqqqq](-1.88,3.7) node[anchor=north west] {$B_1$};
\draw [color=qqqqqq](0.14,4.28) node[anchor=north west] {$B_2$};
\draw [shift={(-2.97,6.23)},color=qqqqqq]  plot[domain=4.211:4.667,variable=\t]({1*2.112*cos(\t r)+0*2.112*sin(\t r)},{0*2.112*cos(\t r)+1*2.112*sin(\t r)});
\draw [color=qqqqqq] (-1.955,1.819)-- (-1.97,1.757);
\draw [color=qqqqqq] (-0.484,1.822)-- (-0.47,1.741);
\draw [color=qqqqqq] (1.121,2.789)-- (1.263,2.686);
\draw [shift={(0.817,1.729)},color=qqqqqq]  plot[domain=1.989:2.193,variable=\t]({1*2.827*cos(\t r)+0*2.827*sin(\t r)},{0*2.827*cos(\t r)+1*2.827*sin(\t r)});
\draw [shift={(-5.213,9.59)},color=qqqqqq]  plot[domain=5.555:5.659,variable=\t]({1*7.528*cos(\t r)+0*7.528*sin(\t r)},{0*7.528*cos(\t r)+1*7.528*sin(\t r)});
\draw [color=qqqqqq](2.62,4.52) node[anchor=north west] {$ B_1$};
\draw [color=qqqqqq](4.42,4.56) node[anchor=north west] {$B_2$};
\draw [->,>=stealth',color=qqqqqq] (2.8,3.72) -- (5.02,3.73);
\draw [->,>=stealth',color=qqqqqq] (5.02,3.73) -- (5.72,4.43);
\draw [->,>=stealth',color=qqqqqq] (5.02,3.73) -- (5.72,3.03);
\draw [->,>=stealth',color=qqqqqq] (2.8,3.72) -- (2.1,4.42);
\draw [->,>=stealth',color=qqqqqq] (2.8,3.72) -- (3.3,2.86);
\draw [->,>=stealth',color=qqqqqq] (2.8,3.72) -- (3.66,3.22);
\draw [shift={(-3.407,-2.554)},color=qqqqqq]  plot[domain=1.387:1.52,variable=\t]({1*6.689*cos(\t r)+0*6.689*sin(\t r)},{0*6.689*cos(\t r)+1*6.689*sin(\t r)});
\begin{scriptsize}
\fill [color=qqqqqq] (2.8,3.72) circle (1.0pt);
\fill [color=qqqqqq] (5.02,3.73) circle (1.0pt);
\end{scriptsize}
\end{tikzpicture}
\caption{Graph of retraction of $\Omega$.}
\label{fig18}
\end{figure}
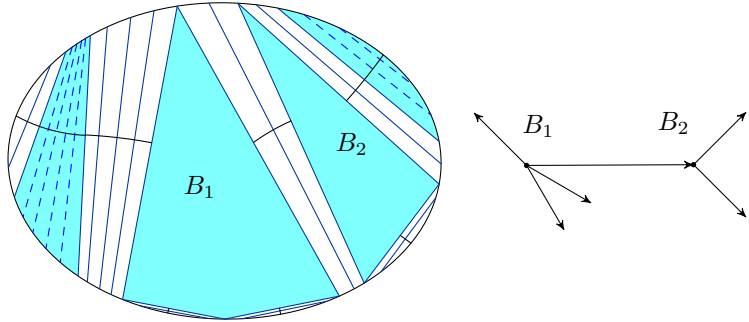

Now since $\Omega$ is convex and simply-connected, we claim that two bodies are
connected through one chain of bodies and arms: it suffices to consider the
graph obtained by retracting bodies to vertices and arms to edges. This graph is
simply connected because it is a deformation retract of $\Omega$. Therefore
every two vertices are connected through only one chain of edges, which proves
the claim (Figure \ref{fig18}).

We begin by a central body $B_1$ and define our approximating sequence on each
arm as in subsection \ref{sapproximation}. Note that for this final purpose, we
have
constructed our approximating smooth isometric immersion to be affine near both
ends, this allows us to apply an affine transformation to the target space of each
arm so that the affine regions near its ends
join together smoothly all the way till we reach $B_2$. Meanwhile, we also
apply an affine transformation to $u(B_2)$ so
that it joins the last arm smoothly. It is easy to see from the uniform
convergence of each term in representation (\ref{0ui1}) that such affine
transformation converges to identity as $m\rightarrow 0$. Now we continue our
construction using
$B_2$ as a new starting point. Note that we will never come back to $B_1$
because they are connected through only one chain of arms. The construction of the approximating sequence on 
the entire domain $\Omega$ is complete.
\hspace*{\fill} $\Box$

\appendix
\section{Proof of Lemma \ref{0difficult}}\label{AA}

\textbf{Step 1.} Recall from the
matrix of moving frame defined in Subsection \ref{sleading}. that 
$\ds \gamma
^{\prime \prime
}(t)= \sum_{i=1}^{n-1} \kappa_i(t)\mathbf{N}_i(t))$,
with $\kappa_i$ bounded. We can choose uniformly bounded smooth functions
$\tilde{\kappa}_{i,m}\rightarrow \kappa_i$ a.e. on $[0,\ell]$, and hence in
measure due to the fact that $[0,\ell]$ is bounded. Since the sequence $\tilde{\kappa}_{i,m}$
are uniformly bounded, it follows $\tilde{\kappa}_{i,m}\rightarrow \kappa_i$ in
$L^p$ for all $1\leq p <\infty$. Similarly we can find uniformly bounded smooth
functions $\kappa_{i_j,m}\rightarrow \kappa_{i_j}$ a.e. on $[0,\ell]$ (hence in
$L^p$ for all $1\leq p<\infty$) for $\kappa_{i_j}, 1\leq i,j\leq n-1$. By
solving the following system of ODEs:  
\begin{align*}
   \left ( \begin{array}{c}  \Gamma_m' \\ \mathbf{N}_{1,m}  \\ \mathbf{N}_{2,m}
\\ \vdots  \\ \mathbf{N}_{n-1,m} \end{array} \right )'
= \widetilde {\mathcal{K}}_m
\left ( \begin{array}{c}   \Gamma_m' \\ \mathbf{N}_{1,m}  \\ \mathbf{N}_{2,m} \\
\vdots  \\ \mathbf{N}_{n-1,m} \end{array} \right ),
\end{align*} where the matric $\widetilde{\mathcal{K}}_m$ is given by 
\begin{align*} 
\widetilde{\mathcal{K}}_m=
\begin{pmatrix}
 0 & \tilde{\kappa}_{1,m} & \tilde{\kappa}_{2,m} & \cdots &
\tilde{\kappa}_{n-1,m}\\
 -\tilde{\kappa}_{1,m} & 0 & \kappa_{1_2,m} & \cdots & \kappa_{1_{n-1},m} \\
 -\tilde{\kappa}_{2,m} & -\kappa_{1_2,m} & 0 & \cdots & \kappa_{2_{n-1},m} \\
  \vdots & \vdots & \vdots & \cdots & \vdots \\
 -\tilde{\kappa}_{n-1,m} & -\kappa_{1_{n-1},m} & -\kappa_{2_{n-1},m} & \cdots & 0
\\
\end{pmatrix},
\end{align*} we obtain a unique orthogonal frame $(\Gamma_m'(t),
\mathbf{N}_{1,m}(t),\cdots,\mathbf{N}_{n-1,m}(t))$ with initial condition
$\Gamma_m'(0)=\gamma'(0)$, and $\mathbf{N}_{i,m}(0)=\mathbf{N}_i(0)$. We can
then define
$$\Gamma_m(t)=\Gamma(0)+\int_0^t \Gamma'_m(\tau) d\tau.$$
We want to show that $(\Gamma_m', \mathbf{N}_{1,m},\cdots,\mathbf{N}_{n-1,m})\rightarrow
(\gamma', \mathbf{N}_{1},\cdots,\mathbf{N}_{n-1})$ uniformly. This result is given
by the following theorem due to Opial, \cite{opial}, Theorem 1.
\begin{lemma}[Opial]\label{0opial}
Suppose the linear system of differential equations,
\begin{equation}\label{0diffeq}
x'(t)=A_k(t)x(t), \quad x(0)=a_k, \quad k=0,1,2,\cdots
\end{equation}
admits a solution $x_k(t)$ in $[0,\ell]$ for all $k$. Suppose $a_k\rightarrow
a_0$,
$$
\int_0^t A_k(s)\,ds\rightarrow \int_0^t A_0(s)\,ds
$$
uniformly for all $t\in [0,\ell]$ and $A_k$ is a bounded sequence in $L^1$, i.e.
$$
\sup_k\|A_k\|_{L^1([0,\ell])}<\infty,
$$ then the solutions
$x_k(t)$ converge to $x_0(t$) uniformly. 
\end{lemma}
Since $\tilde{\kappa}_{i,m}\rightarrow \kappa_i$ and $\kappa_{i_j,m}\rightarrow
\kappa_{i_j}$ in $L^p$ for all $1\leq p<\infty$, in particular for $p=1$, the
conditions in Lemma \ref{0opial} are satisfied, hence $(\Gamma_m',
\mathbf{N}_{1,m},\cdots,\mathbf{N}_{n-1,m})$ converges to $(\gamma',
\mathbf{N}_{1},\cdots,\mathbf{N}_{n-1})$ uniformly. Since
$\ds \Gamma_m^{''}= \sum_{i=1}^{n-1} \tilde{\kappa}_{i,m}\mathbf{N}_{i,m}$, $\Gamma_m^{''}$ are uniformly bounded, and
$\Gamma_m^{''}\rightarrow \gamma^{''}$ a.e. (and hence in $L^p$ for all $1\leq p
<\infty$), Poincar\'{e} inequality for intervals implies that $\Gamma_m\rightarrow
\gamma$ in $W^{2,p}([0,\ell],\bbbr^n)$ for all $1\leq p<\infty$.

However $\Gamma_m$ is not our desired curve since we cannot guarantee that all its
leading fronts intersect outside $\overline{\Omega}$. This happens if $\Gamma_m$
is too ``curvy''. We need to ``flatten'' its curvature continuously. This needs
to be done in several steps:

\textbf{Step 2}. We construct
$\tilde{\tilde{\kappa}}_m=(\tilde{\tilde{\kappa}}_{1,m},\cdots,\tilde{\tilde{\kappa
}}_{n-1,m})$ continuous on $t\in [0,\ell]$ and
for each $t\in [0,\ell]$ and $s=(s_1,\cdots,s_{n-1})\in \mathbb{S}^{n-2}$,
\begin{equation}\label{0curvature}
\ds (S^{\Gamma_m}_s(t)+\frac{\rho}{2})(\sum_{i=1}^{n-1} s_i\tilde{\tilde{\kappa}}_{i,m}(t))\leq 1
\end{equation}
where
$$
\ds S^{\Gamma_m}_s(t)=\sup \{S\geq 0:
\Gamma_m(t)+S\bigl(\sum_{i=1}^{n-1} s_i\mathbf{N}_{i,m}(t) \bigr)\in \Omega\}. 
$$

We first need the following lemma using implicit function theorem for $C^1$ functions. 
\begin{lemma}\label{0implicit}
$S^{\Gamma_m}_s(t)$ is uniformly continuous on $(s,t)\in \mathbb{S}^{n-2}\times [0,\ell]$ and 
$S^{\Gamma_m}_s(t)$ converges to  $S^{\gamma}_s(t)$
uniformly on $(s,t)\in \mathbb{S}^{n-2}\times [0,\ell]$. 
\end{lemma}
{\em Proof.}
Let $t_0\in [0,\ell]$ and $s^0=(s_1^0,\cdots,s_{n-1}^0)\in \mathbb{S}^{n-1}$ be
arbitrary. We parametrize locally $\mathbb{S}^{n-2}$ by the polar coordinates:
$s_i=s_i(\theta)$ where $\theta=(\theta_1,\cdots,\theta_{n-2})\in U_1\subset
[0,\pi)^{n-3}\times [0,2\pi)$. Let $\theta^0\in U_1$ be such that
$s^0_i=s_i(\theta^0)$.

Let $\gamma^0=\gamma(t_0)$ and $\mathbf{N}_i^0=\mathbf{N}_i(t_0)$. Let $x_0$ be
the intersection of the line segment $\ds L=\{
\gamma^0+S\bigl( \sum_{i=1}^{n-1} s_i^0\mathbf{N}^0_i \bigr), 0\leq
S\}$ and $\partial \Omega$. Then
$\ds x_0=\gamma^0+S_0\bigl( \sum_{i=1}^{n-1} s^0_{i}\mathbf{N}^0_i \bigr)$ for some $S_0>0$.

Since $\Omega$ is a $C^1$ domain, there exits an open subset of $U_2\subset
\bbbr^{n-1}$ and a $C^1$ function $\alpha:U_2\rightarrow \partial \Omega$ and
$\alpha(\eta^0_1,\cdots,\eta^0_{n-1})=x_0$ for some $(\eta^0_1,\cdots,\eta^0_{n-1})\in
U_2$.

Consider $F:\bbbr^n\times \bbbr^{n\times (n-1)}\times U_1 \times \bbbr  \times
U_2\rightarrow \bbbr^n$
\begin{equation*}
F(\gamma, \mathbf{N}_1,\cdots\mathbf{N}_{n-1},\theta, S, \eta_1,\cdots,\eta_{n-1})=\\
\ds \gamma+S\bigl(\sum_{i=1}^{n-1} s_i(\theta)\mathbf{N}_i \bigr)-\alpha(\eta_1,\cdots,\eta_{n-1}).
\end{equation*}
Since $x_0\in \partial \Omega\cap L$,
$$F(\gamma^0, \mathbf{N}_1^0,\cdots,\mathbf{N}_{n-1}^0, \theta^0, S^0,
\eta_1^0,\cdots,\eta_{n-1}^0)=0. $$

Let
$$\mathbf{x}= (\gamma, \mathbf{N}_1,\cdots,\mathbf{N}_{n-1}, \theta), \quad \mathbf{y}=(S,\eta_1,\cdots,\eta_{n-1}) $$ and 
$$\alpha_{k}:=\frac{\partial \alpha }{\partial \eta_k}, \quad 1\leq k\leq n-1. $$
Then
\begin{multline*}
\det\bigl[(\frac{\partial F}{\partial \mathbf{y}}(\gamma^0,
\mathbf{N}_1^0,\cdots,\mathbf{N}_{n-1}^0,\theta^0, S^0,
\eta_1^0,\cdots,\eta_{n-1}^0)\bigr]=\\
\ds \det\bigl[\sum_{i=1}^{n-1}  s_i(\theta^0)\mathbf{N}^0_i 
,\alpha_{1}(\eta_1^0,\cdots,\eta_{n-1}^0),\cdots,
\alpha_{n-1}(\eta_1^0,\cdots,\eta_{n-1}^0) \bigr]\neq 0.
\end{multline*}
Otherwise, the line segment $L$ would be parallel to the tangent plane of
$\partial \Omega$ at $x_0$, which is not possible since $\Omega$ is convex.

By implicit function theorem, there is an open neighborhood 
$$
V_1\subset
\bbbr^n\times \bbbr^{n\times (n-1)}\times U_1
$$ of $\mathbf{x}_0=(\gamma^0,
\mathbf{N}^0_1,\cdots,\mathbf{N}^0_{n-1}, \theta^0)$, $V_2\subset \bbbr\times U_2$
of $\mathbf{y}_0=(S^0, \eta^0_1,\cdots,\eta^0_{n-1})$, and a $C^1$ map
$\mathbf{y}: V_1\rightarrow V_2$ such that
$$
F(\mathbf{x},\mathbf{y}(\mathbf{x}))=F(\mathbf{x},S(\mathbf{x}),\eta_1(\mathbf{x
}),\cdots,\eta_{n-1}(\mathbf{x}))=0.
$$
for all $\mathbf{x}\in V_1$.

Since $\gamma$, $\mathbf{N}_i, 1\leq i\leq n-1$ are Lipschitz on $[0,\ell]$ and
$\Gamma_m\rightarrow \gamma$ uniformly and $\mathbf{N}_{i,m}\rightarrow
\mathbf{N}_i$ uniformly on $[0,\ell]$ for all $1\leq i\leq n-1$, there exists an
open interval $O\subset \bbbr$ containing $t_0$, an open subset $\Delta \subset
U_1$ containing $\theta_0$ and an integer $M$ such that for all $t\in
[0,\ell]\cap O$, $\theta\in \Delta$ and $m\geq M$,
$$\mathbf{x}(t,\theta)=(\gamma(t), \mathbf{N}_1(t),\cdots,\mathbf{N}_{n-1}(t),
\theta)\in V_1 \quad \textrm{and}, $$
$$ \mathbf{x}_m(t,\theta)=(\Gamma_m(t),
\mathbf{N}_{1,m}(t),\cdots,\mathbf{N}_{n-1,m}(t), \theta)\in V_1. $$
Apparently $\mathbf{x}_m(t,\theta)\rightarrow \mathbf{x}(t,\theta)$ uniformly
for all $t\in [0,\ell]\cap O$ and $\theta\in \Delta$. Since $S$ is $C^1$ on $\mathbf{x}\in V_1$,
\begin{equation}\label{0uniform3}
S(\mathbf{x}_m(t,\theta))\rightarrow S(\mathbf{x}(t,\theta)) \textrm{ uniformly
on } t\in [0,\ell]\cap O \textrm{ and } \theta\in \Delta.
\end{equation}
Moreover, since $S$ is $C^1$ on $\mathbf{x}\in V_1$ and $\mathbf{x}$ is uniformly continuous on $t\in
[0,\ell]\cap O$ and $s\in s(\Delta)$, $S$ is uniformly continuous on $t\in
[0,\ell]\cap O$ and $s\in s(\Delta)$.

Now note that since $F\bigl(\mathbf{x}(t,\theta),
\mathbf{y}(\mathbf{x}(t,\theta))\bigr)=0$ and $F\bigl(\mathbf{x}_m(t,\theta),
\mathbf{y}(\mathbf{x}_m(t,\theta))\bigr)=0$, for each $s=
s(\theta) \in s(\Delta)\subset \mathbb{S}^{n-2}$, we have
$S_s^\gamma(t)=S(\mathbf{x}(t,\theta))$ and
$S_s^{\Gamma_m}(t)=S(\mathbf{x}_m(t,\theta))$. Thus by (\ref{0uniform3}),
$$S_s^{\Gamma_m}(t)\rightarrow S_s^\gamma(t) \textrm{ uniformly on } t\in
[0,\ell]\cap O \textrm{ and } s\in s(\Delta),$$
and $S_s^{\Gamma_m}(t)$ is uniformly continuous on $t\in
[0,\ell]\cap O$ and $s\in s(\Delta)$. 

It remains to observe that since $[0,\ell]$ and $\mathbb{S}^{n-2}$ are both
compact they  can be covered by a finite union of neighborhoods on which
(\ref{0uniform3}) holds. The proof is complete.
\hspace*{\fill} $\Box$

Define
$$
\ds \lambda_m(t):=\min\Bigl\{1, 1/  ( {\sup_{|s|=1}\{(S^{\Gamma_m}_s(t)+\frac{
\rho}{2})(\ds \sum_{i=1}^{n-1} s_i\tilde{\kappa}_{i,m}(t)) \} ) }
\Bigr\}$$
where $\tilde{\kappa}_{i,m}, 1\leq i\leq n-1$ are those found in Step 1. A first
observation is that $0< \lambda_m\leq 1$. Indeed, there must exist
$s\in\mathbb{S}^{n-2}$ such that
$\ds \sum_{i=1}^{n-1} s_i\tilde{\kappa}_{i,m}(t) \geq0$ so the
supreme over all $s\in\mathbb{S}^{n-2}$ must be nonnegative. On the
other hand, $S^{\Gamma_m}_s$ as well as all
$\tilde{\kappa}_{i,m}$ are bounded so $\lambda_m$ is bounded below by a positive number. 

Second, we observe that  $\lambda_m$ is continuous. Indeed, by Lemma \ref{0implicit}, 
$$(S^{\Gamma_m}_s(t)+\rho/2)(\ds \sum_{i=1}^{n-1}  s_i\tilde{\kappa}_{i,m}(t))$$ is uniformly 
continuous on $(s,t)\in \mathbb{S}^{n-2}\times [0,\ell]$. 
Hence the supreme over $\mathbb{S}^{n-2}$ is attained and a simple argument gives
$$h(t):=\sup_{|s|=1}\{(S^{\Gamma_m}_s(t)+\frac{\rho}{2})(\ds \sum_{i=1}^{n-1} s_i\tilde{\kappa}_{i,m}(t))\}$$ is continuous. 

We then define a vector valued function
$\tilde{\tilde{\kappa}}_{m}=(\tilde{\tilde{\kappa}}_{1,m},\cdots,\tilde{\tilde{
\kappa}}_{n-1,m})$ as
$$(\tilde{\tilde{\kappa}}_{1,m}(t),\cdots,\tilde{\tilde{\kappa}}_{n-1,m}
(t)):=\lambda_{m}(t)(\tilde{\kappa}_{1,m}(t),\cdots,\tilde{\kappa}_{n-1,m}
(t))$$
$\tilde{\tilde{\kappa}}_m$ is obviously continuous. It remains to show
$\tilde{\tilde{\kappa}}_m$ satisfies (\ref{0curvature}). Indeed, for any
$s=(s_1,\cdots,s_{n-1})\in \mathbb{S}^{n-2}$,
\begin{equation*}
\ds (S^{\Gamma_m}_s(t)+\frac{\rho}{2})(\sum_{i=1}^{n-1} s_i\tilde{\tilde{\kappa}}_{i,m}(t))\\
\ds =\lambda_{m}(S^{\Gamma_m}_s(t)+\frac{\rho}{2})(\sum_{i=1}^{n-1} s_i\tilde{\kappa}_{i,m}
(t)).
\end{equation*}

If $\ds \sum_{i=1}^{n-1}  s_i\tilde{\kappa}_{i,m}(t) \geq 0$,
then by
the definition of $\lambda_m$,
\begin{equation*}
\lambda_{m}(t)(S^{\Gamma_m}_s(t)+\frac{\rho}{2})
\ds \bigl( \sum_{i=1}^{n-1} s_i\tilde{\kappa}_{i,m}(t) \bigr)\\
\ds  \leq
\min\bigl\{(S^{\Gamma_m}_s(t)+\frac{\rho}{2})
\bigl( \sum_{i=1}^{n-1} s_i\tilde{\kappa}_{i,m} (t) \bigr),\, 1\bigr\}\leq 1.
\end{equation*}
If $\ds \sum_{i=1}^{n-1} s_i\tilde{\kappa}_{i,m}(t) < 0$, then
\begin{equation*}
\ds \lambda_{m}(t)(S^{\Gamma_m}_s(t)+\frac{\rho}{2})(\sum_{i=1}^{n-1} s_i\tilde{\kappa}_{i,m} (t))<0\leq 1.
\end{equation*}
Thus (\ref{0curvature}) is satisfied.

\textbf{Step 3}. We want to show that
$(\tilde{\tilde{\kappa}}_{1,m},\cdots,\tilde{\tilde{\kappa}}_{n-1,m})\rightarrow
(\kappa_{1},\cdots,\kappa_{n-1})$ a.e. Indeed, we know that
$\tilde{\kappa}_m=(\tilde{\kappa}_{1,m},\cdots,\tilde{\kappa}_{n-1,m})\rightarrow
(\kappa_{1},\cdots,\kappa_{n-1})$ a.e.. Therefore, all we need to show is
$\lambda_{m}\rightarrow 1 \quad \textrm{a.e.}$.

By possibly replacing $\lambda_m$ by a subsequence, it suffices to prove $\lambda_m \rightarrow 1$ in measure. From the
definition of $\lambda_m$, it is enough to show the Lebesgue measure of the set:
$$
\ds E_m=\{t\in [0,l], \exists s\in
\mathbb{S}^{n-2},(S^{\Gamma_m}_s(t)+\frac{\rho}{2})( \sum_{i=1}^{n-1} s_i\tilde{\kappa}_{i,m} (t))> 1\}
$$
goes to zero. First by assumption, $L_s^{\gamma}(t)- S_s^{\gamma}(t)>\rho>0$ and by Lemma
\ref{0distance2},
$\ds L_s^{\gamma}(t)( \sum_{i=1}^{n-1} \kappa_i(t)s_i )\leq 1$, thus,
\begin{equation}\label{0a1}
\ds (S_s^\gamma(t)+\rho)(\sum_{i=1}^{n-1} s_i\kappa_i(t) )\leq 1, 
\end{equation}
for all $t\in [0,\ell]$ and $s\in \mathbb{S}^{n-2}$.
Indeed, if $\ds \sum_{i=1}^{n-1} s_i\kappa_i(t))\geq 0$,
\begin{equation*}
\ds (S_s^\gamma(t)+\rho)(\sum_{i=1}^{n-1} s_i\kappa_i(t))\leq 
\ds L_s^{\gamma}(t)(\sum_{i=1}^{n-1} \kappa_i(t)s_i)\leq 1
\end{equation*}
which gives (\ref{0a1}).

If $t\in E_m$, there is $s\in \mathbb{S}^{n-2}$ such that
$$
\ds \sum_{i=1}^{n-1} s_i\tilde{\kappa}_{i,m}(t) >\frac{1}{S^
{\Gamma_m}_s(t)+\rho/2}. $$
Therefore all $t\in E_m$ and our choice of $s=s(t)$ as above, we have,
\begin{equation*}
\begin{array}{ll}
\ds |\tilde{\kappa}_m(t)-\kappa(t)| & \ds \geq
\sum_{i=1}^{n-1} s_i\tilde{\kappa}_{i,m}(t) 
  -( \sum_{i=1}^{n-1} s_i\kappa_i(t) ) \\ & \ds    >
\frac{\rho/2+S_s^\gamma(t)-S^{\Gamma_m}_s(t)}{(S^{\Gamma_m}
_s(t)+\rho/2)(S_s^\gamma(t)+\rho)}
   \geq \frac{\rho/2-|S_s^\gamma(t)-S^{\Gamma_m}_s(t)|}{\rho^2/2}.
\end{array}
\end{equation*}
By Lemma \ref{0implicit}, we have
\begin{equation*}
S^{\Gamma_m}_s(t)\rightarrow S^\gamma_s(t) \textrm{ uniformly on } s\in
\mathbb{S}^{n-2} \textrm{ and } t\in [0,\ell],
\end{equation*}
then we can find $m$ sufficiently large so that
$|S_s^\gamma(t)-S^{\Gamma_m}_s(t)|< \rho/4$ for all $s\in \mathbb{S}^{n-2}$ and
$t\in [0,\ell]$. Since $\tilde{\kappa}_{m}\rightarrow \kappa$ a.e.,
$$\lim_{m\rightarrow \infty}|E_m|\leq \lim_{m\rightarrow \infty}|\{t:
|\tilde{\kappa}_m(t)-\kappa(t)|\geq \frac{1}{2\rho}\}|=0$$
which is what we wanted to show.

\textbf{Step 4}. Since
$\tilde{\tilde{\kappa}}_m=(\tilde{\tilde{\kappa}}_{1,m},\cdots,\tilde{\tilde{\kappa
}}_{n-1,m})$ are continuous, for each $m$ we can find $\kappa_{m}$ smooth and
$|\tilde{\tilde{\kappa}}_{m}-\kappa_m|\rightarrow 0$ uniformly on $t\in
[0,\ell]$. Hence for $m$ sufficiently large,
\begin{equation}\label{0curvature2}
\ds (S^{\Gamma_m}_s(t)+\frac{\rho}{4})(\sum_{i=1}^{n-1} s_i\kappa_{i,m}(t) )\leq 1
\end{equation}

\textbf{Step 5}. We now define our desired curve $\gamma_m$. Given
$\kappa_m=(\kappa_{1,m},\cdots\kappa_{n-1,m})$ smooth as found in Step 4, and
$\kappa_{i_j,m}\rightarrow \kappa_{i_j}$ found in step 1, we again solve the system of ODEs
\begin{align*}
   \left ( \begin{array}{c}  \gamma_m' \\ \mathbf{N}_{1,m}  \\ \mathbf{N}_{2,m}
\\ \vdots  \\ \mathbf{N}_{n-1,m} \end{array} \right )'
=  \mathcal{K}^{n\times n}_m
\left ( \begin{array}{c}   \gamma_m' \\ \mathbf{N}_{1,m}  \\ \mathbf{N}_{2,m} \\
\vdots  \\ \mathbf{N}_{n-1,m} \end{array} \right ),
\end{align*} where 
\begin{align*}
\mathcal{K}^{n\times n}_m=
\begin{pmatrix}
 0 & \kappa_{1,m} & \kappa_{2,m} & \cdots & \kappa_{n-1,m} \\
 -\kappa_{1,m} & 0 & \kappa_{1_2,m} & \cdots & \kappa_{1_{n-1},m} \\
 -\kappa_{2,m} & -\kappa_{1_2,m} & 0 & \cdots & \kappa_{2_{n-1},m} \\
  \vdots & \vdots & \vdots & \cdots & \vdots \\
 -\kappa_{n-1,m} & -\kappa_{1_{n-1},m} & -\kappa_{2_{n-1},m} & \cdots & 0 \\
\end{pmatrix},
\end{align*}
 and denote by the orthogonal frame $(\gamma_m'(t),
\mathbf{N}_{1,m}(t),\cdots\mathbf{N}_{n-1,m}(t))$ the unique solution with initial conditions
$\gamma_m'(0)=\gamma'(0)$ and $\mathbf{N}_{i,m}(0)=\mathbf{N}_i(0)$. Moreover,
by Lemma \ref{0opial}, $(\gamma_m'(t),
\mathbf{N}_{1,m}(t),\cdots\mathbf{N}_{n-1,m}(t))\rightarrow (\gamma'(t),
\mathbf{N}_{1}(t),\cdots\mathbf{N}_{n-1}(t))$ uniformly. Let
$$\gamma_m(t)=\gamma(0)+\int_0^t \gamma'_m(\tau) d\tau.$$
We claim $\gamma_m$ satisfies for $m$ sufficiently large,
\begin{equation}\label{0curvature3}
\ds (S^{\gamma_m}_s(t)+\frac{\rho}{8})(\sum_{i=1}^{n-1} s_i\kappa_{i,m}(t) )\leq 1
\end{equation}
Indeed, by the same argument of Lemma \ref{0implicit} using implicit function
theorem, $S^{\gamma_m}_s$ also converges to $S^\gamma_s$
uniformly. Together with Lemma \ref{0implicit} we obtain that 
$|S^{\gamma_m}_s-S^{\Gamma_m}_s|$ converges to $0$ uniformly.
Thus the claim follows from (\ref{0curvature2}).

\textbf{Step 6.} Finally, we claim that orthogonal fronts satisfy
$F_{\gamma_m}(t)\cap F_{\gamma_m}(\tilde{t})\cap \overline{\Omega}=\emptyset$
for all $t,\tilde{t}\in [0,\ell]$.

For $\gamma_m$ and its moving frame $(\gamma_m',
\mathbf{N}_{1,m},\cdots,\mathbf{N}_{n-1,m})$ found in Step 5, let
$$
\Phi_m:[0,\ell]\times \bbbr^{n-1}\rightarrow \bbbr^n
$$ be defined as
$$
\ds \Phi_m(t,s)=\gamma_m(t)+ \sum_{i=1}^{n-1} s_i\mathbf{N}_{i,m}(t).$$
Let $\Sigma^{\gamma_m}=\{(t,s): \Phi_m(t,s)\in \overline{\Omega}\}$. By the same
argument as Lemma \ref{0phi2}, $\Phi_m$ maps $\Sigma^{\gamma_m}$ onto
$\overline{\Omega(\gamma_m)}$ where $\Omega(\gamma_m)$ is the subset of $\Omega$
covered by all orthogonal fronts $F_{\gamma_m}(t),t\in [0,\ell]$. By the same computation as in (\ref{0jaco1}),
$$
\ds J_{\Phi_m}(t,s)= 1-\sum_{i=1}^{n-1} s_i\kappa_{i,m}(t).
$$
Let $d:=\textrm{diam}(\Omega)$, we claim that
\begin{equation*}
\ds 1-\sum_{i=1}^{n-1} s_i\kappa_{i,m}(t) \geq \min\{\rho/16d,1/2\},
\end{equation*}
for all $(t,s)\in \Sigma^{\gamma_m}$. Indeed, if
$\ds \sum_{i=1}^{n-1} (s_i/|s|)\kappa_{i,m}(t) \geq 1/{2d}$,
then by (\ref{0curvature3}),
\begin{equation*}
\ds 1-|s|\Bigl( \sum_{i=1}^{n-1} \frac{s_i}{|s|}\kappa_{i,m}(t) \Bigr) 
\ds \geq 1-S^{\gamma_m}_s(t)
\Bigl(\sum_{i=1}^{n-1} \frac{s_i}{|s|}\kappa_{i,m}(t) \Bigr) 
\ds \geq\frac{\rho}{8}\Bigl(\sum_{i=1}^{n-1} \frac{s_i}{|s|}\kappa_{i,m}(t) \Bigr) 
\geq\frac{\rho}{8}\cdot \frac{1}{2d}.
\end{equation*}
If $\ds \sum_{i=1}^{n-1} (s_i/|s|)\kappa_{i,m}(t) < 1/{2d}$,
then
$$
\ds 1-|s|\Bigl(\sum_{i=1}^{n-1} \frac{s_i}{|s|}\kappa_{i,m}(t) \Bigr)>1-\frac{|s|}{2d}\geq \frac{1}{2}. 
$$
Hence, the claim follows. By Inverse function theorem due to Clarke
\cite{clarke}, $\Phi$ admits a local Lipschitz inverse, actually a
\textit{global} Lipschitz inverse
$\Phi_m^{-1}:\overline{\Omega(\gamma_m)}\rightarrow \Sigma^{\gamma_m}$ since the
Jacobian is everywhere bounded below by a positive constant in
$\Sigma^{\gamma_m}$. In particular, $\Phi_m$ is one-to-one on
$\Sigma^{\gamma_m}$. This implies all orthogonal fronts $F_{\gamma_m}(t),t\in
[0,\ell]$ meet outside $\overline{ \Omega}$. The proof of Lemma
\ref{0difficult} is complete.
\hspace*{\fill} $\Box$

\end{document}